\newcount\Comments  
\documentclass[letterpaper,10pt]{scrartcl}

\usepackage[utf8]{inputenc}
\usepackage[english]{babel}
\usepackage{amsmath,amstext,amsbsy,amsopn,amscd,amsxtra,upref,amssymb}
\usepackage{amsfonts,bm}
\usepackage{amsthm}
\usepackage{color}
\usepackage{graphicx}
\usepackage{sidecap}
\usepackage{multirow}
\usepackage{booktabs}
\usepackage{bm}
\usepackage{fullpage}
\usepackage{wrapfig}

\usepackage{tabularx}
\usepackage[font=small]{caption}
\usepackage{subfig}
 
\usepackage{longtable}
\usepackage{rotating}

\usepackage{algorithm}
\usepackage{algpseudocode}
\extrafloats{100}

\usepackage{lineno,hyperref}

\graphicspath{{figures/}}
\DeclareGraphicsExtensions{.pdf,.eps,.png,.jpg,.jpeg}

\usepackage{grffile}

\newcommand{\R}{\mathbb{R}}

\newcommand{\bfa}{\boldsymbol a}

\newcommand{\bfx}{\boldsymbol x}
\newcommand{\bfX}{\boldsymbol X}

\newcommand{\Mcal}{\mathcal{M}}

\newcommand{\Ucal}{\mathcal{U}}

\newcommand{\Xcal}{\mathcal{X}}

\newcommand{\Zcal}{\mathcal{Z}}

\newcommand{\bfmu}{\boldsymbol \mu}

\newcommand{\bfb}{\boldsymbol b}
\newcommand{\bfc}{\boldsymbol c}
\newcommand{\bfe}{\boldsymbol e}

\newcommand{\bfu}{\boldsymbol u}

\newcommand{\bfW}{\boldsymbol W}

\newcommand{\bfSigma}{\boldsymbol\Sigma}

\newcommand{\<}{\langle}
\renewcommand{\>}{\rangle}

\DeclareMathOperator*{\argmin}{argmin}

   {\begin{trivlist}\item[]{\bfseries\sffamily Keywords:}\ }
   {\end{trivlist}}

\theoremstyle{definition}

\newtheorem{proposition}{Proposition}

\date{}

\ifpdf
\hypersetup{
  pdftitle={Neural Galerkin Schemes with Active Learning for High-Dimensional Evolution Equations}, 
  pdfauthor={Joan Bruna and Benjamin Peherstorfer and Eric Vanden-Eijnden} 
}
\fi

\title{Neural Galerkin Schemes with Active Learning for High-Dimensional Evolution Equations}

\author{Joan Bruna\thanks{Courant Institute of Mathematical Sciences, New York University, New York, NY 10012} \and Benjamin Peherstorfer\footnotemark[1] \and Eric Vanden-Eijnden\footnotemark[1]}

\begin{document}

\maketitle

\begin{abstract}
Deep neural networks have been shown to provide accurate function approximations in high dimensions. However, fitting network parameters requires informative training data that are often challenging to collect in science and engineering applications. This work proposes Neural Galerkin schemes based on deep learning that generate training data with active learning for numerically solving high-dimensional partial differential equations. Neural Galerkin schemes build on the Dirac-Frenkel variational principle to train networks by minimizing the residual sequentially over time, which enables adaptively collecting new training data in a self-informed manner that is guided by the dynamics described by the partial differential equations. This is in contrast to other machine learning methods that aim to fit network parameters globally in time without taking into account training data acquisition. Our finding is that the active form of gathering training data of the proposed Neural Galerkin schemes is key for numerically realizing the expressive power of networks in high dimensions. Numerical experiments demonstrate that Neural Galerkin schemes have the potential to enable simulating phenomena and processes with many variables for which traditional and other deep-learning-based solvers fail, especially when features of the solutions evolve locally such as in high-dimensional wave propagation problems and interacting particle systems described by Fokker-Planck and kinetic equations.
\end{abstract}

\section{Introduction}
Partial differential equations (PDEs) are used to describe the dynamics of systems in a wide variety of science and engineering applications. Many of these equations, however, are challenging to solve, analytically and even numerically. Classical techniques from scientific computing such as the finite element method and other grid-based methods can work well if the solution domain is low dimensional; however, in higher dimensions, grid-based methods fail because their computational costs increase at an exponential rate with the dimension of the domain. This curse of dimensionality leaves simulating many systems of interest out of scope of classical PDE solvers. For example, interacting particles described by Schr\"odinger's and other kinetic equations such as Fokker-Planck and Boltzmann equations lead to high-dimensional systems when the number of particles is more than a few.
Even multiscale techniques as sophisticated as multigrid, fast multipole, and adaptive-mesh refinement methods fail in such instances.
Thus, high-dimensional approximation problems require a fundamentally different notion of discretization than via grids to circumvent the curse of dimensionality. In supervised machine learning applications, such adaptivity is achieved with deep neural networks (DNNs) via feature or representation learning.

The last few years have seen many developments in the direction of using DNNs for numerically solving PDEs, with the introduction of techniques such as physics-informed neural networks and the Deep Ritz method \cite{doi:10.1002/cnm.1640100303,SIRIGNANO20181339,RAISSI2019686,BERG201828}.  In these schemes, the solution of a time-dependent PDE is represented with a DNN and the parameters of the network are optimized by minimizing the residual of the PDE on collocation points from the temporal and spatial domains. These collocation points can lie on a grid, or be sampled randomly in time-space. Either way, their positions are usually not tailored to the (unknown) solution of the PDE. In high dimensions, the lack of adaptivity in the training procedure ultimately suffers from the curse of dimensionality  because all features of the solution have to be discovered without guidance. In addition, the training treats space and time equivalently, which is inefficient for initial value problems if the goal is predicting the PDE solution over long time ranges, and it can lead to a loss of causality in the solution.

In this work we take a different approach. We develop time-integrators for PDEs that use DNNs to represent the solution but update the parameters sequentially from one time slice to the next rather than globally over the whole time-space domain. To achieve this, we follow the well-established Dirac-Frenkel variational principle \cite{dirac_1930,Frenkel1934} that is widely used in computational chemistry \cite{Lubich2008} and scientific computing \cite{doi:10.1137/050639703,SAPSIS20092347,hesthaven_pagliantini_rozza_2022} and that has been considered in the context of PDEs and DNNs in the works~\cite{PhysRevE.104.045303,doi:10.1137/21M1415972}. It updates the network parameters so that the PDE residual is minimized as the numerical solution is evolved in time. The key novelty of our Neural Galerkin schemes with active learning is that we combine the sequential updating of the network parameters given via the Dirac-Frenkel variational principle with adaptive sampling (``data collection'') over time to determine at which sampling points to evaluate the residual to guide the parameter updates. This means that the measure---from which samples are drawn for estimating the expected residual---is adapted over time so that the samples follow where the important features of the solutions are, which is critical in high-dimensional spatial domains and when the solution has local features that evolve in the spatial domain.  
Importantly, our scheme uses the PDEs that govern the solution to inform the adaptive sampling as well, rather than relying on a black-box global optimization strategy. No \textit{a~priori} data about the PDE solutions is used, however.  
Overall, we summarize that our approach leverages adaptivity in \emph{both} function approximation and sampling.

\paragraph{Summary of main contributions} Our main results can be summarized as follows: 1. Given a nonlinear parametric representation (e.g.~via a DNN) of the solution of an initial value problem for a time-dependent PDE, we leverage the Dirac-Frenkel variational principle as in \cite{dirac_1930,Frenkel1934,Lubich2008,doi:10.1137/050639703,SAPSIS20092347,hesthaven_pagliantini_rozza_2022,PhysRevE.104.045303,doi:10.1137/21M1415972} to derive systems of nonlinear evolution equations for the parameters of the nonlinear parametric representation.  The system of equations that govern the parameters can then be integrated using standard solvers with different level of sophistication, for example by dynamically choosing the time-step size based on properties of the numerical solution.
In particular, with suitable time-integration schemes, the proposed approach takes larger time steps when possible and corrects to smaller time-step sizes if the dynamics of the solution require it. This is in contrast to typical collocation methods such as the Deep Ritz method that samples over the whole time-space domain without having the opportunity to exploit the smoothness of solutions over time and thus has to draw samples even when the solution hardly varies.

2. The evolution equations that we derive for the DNN parameters involve operators that require estimation via sampling in space. We show how to perform this sampling adaptatively, by using the property of the current solution itself, which is again in contrast to collocation methods.
We propose a dynamical estimation of the loss, which is shown to be essential in two aspects: First it achieves accuracy in situations where static, non-adaptive sampling methods fail to capture the details of the solution that determine its evolution. 
Second, it is key to the interpretation of the PDE solution itself, since the solution can be quite complicated, and sampling its feature to compute statistics typically is the only way to exploit and interpret it.

3. We illustrate the viability and usefulness of our approach on a series of test cases. First we consider two equations in low dimensions, the Korteweg-de Vries  (KdV) and the Allen-Cahn (AC) equations to show that our method can be used to efficiently capture solutions with localized features that move in the spatial domain. The features are solitons in the context of the KdV equation and domain walls in the context of the AC equation, which can both be modeled with few degrees of freedom using neural networks with specific nonlinear units.  Second, we show that our method can be used to simulate PDEs in unbounded, high-dimensional domains, that are out of reach of standard methods. Specifically, we study advection equations of hyperbolic type  with wave speeds that vary in space and time, and Fokker-Planck equations of parabolic type for the evolution of several interacting particles. In both cases, our results show that our method is able to efficiently capture the evolution of the solution and overcome the sampling limitations of collocation-based methods.

\paragraph{Related works} The need for adaptive data acquisition in the context of machine learning for problems in science and engineering has been emphasized in previous works~\cite{rotskoff2021active}. 
There also is a large body of work on numerically solving PDEs with DNN parametrizations based on collocation over the spatiotemporal domain (space-time discretizations), which makes adaptive sampling challenging; see, e.g., Refs.~\cite{doi:10.1002/cnm.1640100303,SIRIGNANO20181339,RAISSI2019686,BERG201828}. Methods that do not rely on collocation typically exploit specific formulations of elliptic and semilinear parabolic PDEs \cite{E2017,Han8505}, and focus on specific approximation tasks such as learning committor functions \cite{Khoo2018,doi:10.1063/1.5110439,rotskoff2021active}, coarse-graining \cite{Bar-Sinai15344,Kochkove2101784118,WANG2020109402}, or de-noising \cite{RUDY2019483}. Thus, these approaches are limited in their scope and, in particular, are not applicable to pure transport dynamics given by kinetic and hyperbolic equations.
There also is a range of surrogate-modeling methods \cite{RozzaPateraSurvey,SIREVSurvey} based on nonlinear parametrizations such as \cite{WANG2019289,DBLP:journals/corr/abs-2108-08481,SMAI-JCM_2021__7__121_0,KHOO_2020,LEE2020108973,TaddeiShock,OLEARYROSEBERRY2022114199,ehrlacher19,doi:10.1137/17M1140571,https://doi.org/10.48550/arxiv.2203.00360,QIAN2020132401,hesthaven_pagliantini_rozza_2022}; however, these methods require access to training sets of ground truth PDE solutions, which is precisely what is unavailable in our setting because classical PDE solvers are cursed by dimension.
The surrogate-modeling method introduced in \cite{LEE2020108973,KIM2022110841} embeds large state vectors stemming from traditional discretizations of PDEs in low-dimensional manifolds and then propagates forward in time the embedding based on the dynamics given by the traditional discretizations; however, this means that in each time step a lifting from low to high dimensions has to be performed, which is inherently computationally expensive. In contrast, the proposed Neural Galerkin schemes directly integrate the parameters of the PDE solutions, rather than embeddings of high-dimensional states. This approach circumvents the requirement of needing traditional discretizations of the underlying PDEs, which avoids the expensive lifting to higher dimensions that is intractable for problems formulated over high-dimensional spatial domains.

There is a wide range of methods that builds on the Dirac-Frenkel variational principle to update the parameters of nonlinear parametrizations, see, e.g., the works on dynamic low-rank approximations and related methods \cite{doi:10.1137/050639703,doi:10.1137/16M1071493,hesthaven_pagliantini_rozza_2022,doi:10.1137/17M1123286,Cecilia2020,SAPSIS20092347,Peherstorfer15aDEIM,doi:10.1137/19M1257275,doi:10.1137/140967787}. Closest to our approach are the works~\cite{PhysRevE.104.045303,doi:10.1137/21M1415972}, which derive the same dynamics for parameters of nonlinear parametrizations based on the Dirac-Frenkel variational principle, where the work \cite{doi:10.1137/21M1415972} additionally derives schemes that conserve quantities. However, both works \cite{PhysRevE.104.045303,doi:10.1137/21M1415972} focus exclusively either on specific nonlinear parametrizations and PDEs where the objective value can be computed analytically or on low-dimensional problems where data acquisition with static, uniform sampling is tractable.  In contrast, our Neural Galerkin schemes with adaptive sampling applies to high-dimensional problems with generic nonlinear parametrizations where adaptive data collection becomes crucial, as will be illustrated by our case studies.

\section{Setup and problem formulation}
In this section, we formulation the problem of numerically solving time-dependent PDEs, which are sometimes referred to as evolution equations. We are particularly interested in PDEs formulated over potentially high-dimensional spatial domains and with solutions that exhibit local features that are transported over time.

\subsection{Evolution equations}
Given a spatial domain $\Xcal \subseteq \R^d$, we will consider the evolution of a time-dependent field $u:[0,\infty) \times \Xcal \to \R$ which at all times belongs to some function space $\Ucal$ and whose dynamics is governed by the PDE
\begin{equation}
 \begin{aligned}
 \partial_t u(t, \bfx) = & f(t,\bfx,u) \qquad &&\text{for}\ \ (t,\bfx) \in [0,\infty) \times\Xcal\,,\\
 u(0, \bfx) = &u_0(\bfx) &&\text{for}\ \ \bfx \in \Xcal
 \end{aligned}
\label{eq:nonlinear:PDE}
\end{equation}
where  $u_0 \in \Ucal$ is the initial condition. By appropriately choosing  $f: [0,\infty) \times \Xcal\times \Ucal \to \R$, Eq.~\eqref{eq:nonlinear:PDE} can represent different PDEs of interest. For example, we get an advection-diffusion-reaction equation by taking 
\[
f(t,\bfx,u) = b(t,\bfx) \cdot \nabla u + a(t,\bfx) : \nabla \nabla u + G(t,\bfx,u)
\]
for a $b: [0,\infty)\times \Xcal\to \R^d$, $a: [0,\infty)\times \Xcal\to \R^d\times \R^d$ and $G:[0,\infty)\times \Xcal\times \R\to \R$.
In the following, we only consider problems with appropriate boundary conditions for Eq.~\eqref{eq:nonlinear:PDE} that make Eq.~\eqref{eq:nonlinear:PDE} well-posed for all $t\in[0,\infty)$, which means that its solution  exists, is unique, and depends continuously on $u_0 \in \Ucal$.

\subsection{Problem formulation}\label{sec:ProbForm}
Widely used approaches for numerically solving PDEs of the form \eqref{eq:nonlinear:PDE} with deep neural networks rely on globally optimizing network parameters over the time-space domain. Let $\tilde{u}: [0, \infty) \times \Xcal \times \Theta \to \mathbb{R}$ be a neural-network parametric approximation of $u$ so that $\tilde{u}(t, \bfx; \theta)$ approximates $u(t, \bfx)$ for a parameter $\theta \in \Theta$. Notice that time $t$ is an input of the network function $\tilde{u}$. Widely used collocation-based methods train the network parameters $\theta$ of $\tilde{u}$ to globally minimize the expected residual 
\begin{equation}
r(t, \bfx; \theta) = \partial_t \tilde{u}(t, \bfx; \theta) - f(t, \bfx, \tilde{u})
\label{eq:ProbForm:GlobalRes}
\end{equation}
of the PDE,
\begin{equation}
\min_{\theta \in \Theta} \mathbb{E}_{(t, \bfx) \sim \nu} \left[|r(t, \bfx; \theta)|^2\right]\,.
\label{eq:ProbForm:ExpLoss}
\end{equation}
The objective of the optimization problem \eqref{eq:ProbForm:ExpLoss} depends on the expected value of the squared residual over the time-space domain with respect to the distribution $\nu$. To estimate the expected value in \eqref{eq:ProbForm:ExpLoss}, the typical approach is to draw $i = 1, \dots, n$ samples $(t_i, \bfx_i)$ from the distribution $\nu$, to form a Monte Carlo estimator based on the drawn samples, and to minimize the Monte Carlo estimate of the expected value as in
\begin{equation}
\min_{\theta \in \Theta} \frac{1}{n}\sum_{i = 1}^n |r(t_i, \bfx_i; \theta)|^2\,,\qquad (t_i, \bfx_i) \sim \nu.
\label{eq:ProbForm:MCEst}
\end{equation}
However, it can be challenging to ensure that the Monte Carlo estimator in \eqref{eq:ProbForm:MCEst} is a good estimate of the expected value in \eqref{eq:ProbForm:ExpLoss}. In particular, in high-dimensional spatial domains $\Xcal$ and if the solutions of the PDEs~\eqref{eq:nonlinear:PDE} develop spatially localized structures, then uniform (or other un-informed, un-guided, solution-independent) sampling schemes quickly fail to provide accurate Monte Carlo estimators; see also \cite{rotskoff2021active}. 
In fact, even though the deep neural network parametrization might be expressive enough and exploit properties of the solution to circumvent the curse of dimensionality from an approximation-theoretic point of view, minimizing the objective in \eqref{eq:ProbForm:MCEst} can be well affected by the curse of dimensionality through the variance of the estimator, especially with local solution features advecting in high-dimensional spatial domains.

\section{Neural Galerkin schemes via the Dirac-Frenkel variational principle}
In this section, we break with the traditional time-space approach of training DNNs in the context of PDEs and instead formulate the training of deep neural network parametrizations via the Dirac-Frenkel variational principle \cite{dirac_1930,Frenkel1934,Lubich2008,PhysRevE.104.045303}, which prescribes a system of ordinary differential equations (ODEs) for the network parameters. This means that the network parameters are ultimately trained by solving a system of ODEs over time, rather than solving a global optimization problem over the time-space domain as in traditional collocation-based methods. Carrying over the dynamics described by the PDEs of interest to the network parameters is of critical importance to our Neural Galerkin approach because it will allow us to use information from previous time steps to inform finding the solution at the current time step, including the measure against which to estimate the residual. 

We emphasize that using the Dirac-Frenkel variational principle to train nonlinear parametrizations has a long history in molecular dynamics and computational chemistry; see \cite{Lubich2008,hesthaven_pagliantini_rozza_2022} for the history and early works. More recent works \cite{PhysRevE.104.045303,doi:10.1137/21M1415972} in the context of deep neural networks derive the same systems of ODEs that we obtain in the following.

\subsection{Nonlinear parametrizations that depend on time} We parametrically represent the solution $u(t)$ at time $t$ as  
$U(\theta(t)) \in \Ucal$ with parameters $\theta(t) \in \Theta$ and $U: \Theta \times \Xcal \to \mathbb{R}$, i.e. we use the ansatz
\begin{equation}
\label{eq:rep:sol}
u(t,\bfx) = U(\theta(t),\bfx)  \qquad \text{for}\ \ (t,\bfx)\in[0,\infty)\times\Xcal\,.
\end{equation}
Two remarks are in order. First, time $t$ enters the parametrization \eqref{eq:rep:sol} via the parameter $\theta(t)$, instead of being an input to the network function as in the global, time-space approaches discussed in Section~\ref{sec:ProbForm}. This is of critical importance to our approach because it will allow us to use information from previous times to inform approximations of the PDE solution at future times. Second, the parametrization $U$ may depend nonlinearly on $\theta(t)$, which is in stark contrast to the majority of classical approximations in scientific computing that have a linear dependence on the parameter \cite{FEMBook}. And it is in stark contrast to traditional, linear model reduction methods that seek linear approximations in subspaces; see \cite{P22AMS}. 

For the representation in Eq.~\eqref{eq:rep:sol} to be complete, i.e. such that for any $u(t)\in \Ucal$ there exists at least one $\theta(t) \in \Theta$  such that $U(\theta(t))\equiv u(t)$, $\Theta$ should itself be an infinite-dimensional function space in general---for more discussion on this point, see the Appendix. In practice, however, we will be interested in situations where the parametric representation $U$ depends on a finite-dimensional parameter $\theta(t)$ such as in a DNN, where $\theta(t)$ denotes the vector of adjustable parameters/weights in this network whose number can be large but is finite---specific DNN representations are discussed below. In this case, using Eq.~\eqref{eq:rep:sol} as ansatz for the solution of Eq.~\eqref{eq:nonlinear:PDE} introduces approximation errors, the magnitude of which we aim to keep low at all times.

\subsection{Training nonlinear parametrization via the Dirac-Frenkel variational principle}\label{sec:NG:DiracFrenkel} We emphasize that we consider the situation where we do not have access to training data to learn the parameter $\theta(t)$ so that $U(\theta(t))$ approximates well the PDE solution. Instead, we use the structure of the governing equation~\eqref{eq:nonlinear:PDE} to find $\theta(t)$ so that the error is low. To this end, note that inserting the ansatz solution~$U(\theta(t))$ in~Eq.~\eqref{eq:nonlinear:PDE}, assuming differentiability of $\theta(t)$ and using $\partial_t U(\theta(t)) = \nabla_{\theta}U(\theta)\cdot \dot{\theta}(t)$, leads to the residual function $r:\Theta\times \dot\Theta \times \Xcal\to\R$ defined as
\begin{equation}
    r_t(\theta,\eta,\bfx)  \equiv  \nabla_\theta U(\theta,\bfx)\cdot \eta  - f(t,\bfx,U(\theta))\label{eq:res}\,,
\end{equation}
where $\dot\Theta$ is the set of time derivatives of $\theta(t)$. Notice that the residual $r_t$ defined in \eqref{eq:res} depends on time $t$, whereas the residual $r$ defined in \eqref{eq:ProbForm:GlobalRes} is global over the time-space domain.

The goal is now to keep the error made by the parametrization described in Eq.~\eqref{eq:rep:sol} low by minimizing the magnitude of the residual $r_t$ with respect to $\eta$, given an appropriate metric. One possibility is to do so globally over some time interval $t\in[0,T]$ with $T>0$; this option, which is discussed in~\ref{sec:comp:coll}, 
has the disadvantage that it leads to a complicated boundary value problem that requires the storage of the entire path $\theta(t)$ for $t\in[0,T]$. Instead, we follow the Dirac-Frenkel variational principle \cite{dirac_1930,Frenkel1934,Lubich2008}, see also \cite{PhysRevE.104.045303} for the principle applied in the context of deep neural network, which leads to an initial value problem that can be solved over arbitrary long times. Specifically, we seek $\theta(t)$ such that for all $t>0$ it holds
\begin{equation}
    \label{eq:min}
    \dot \theta(t) \in \argmin_{\eta\in \dot \Theta} J_t(\theta(t),\eta)
\end{equation}
where we define the objective function $J_t: \Theta\times\dot \Theta \to \mathbb{R}$ as
\begin{equation}
J_t(\theta,\eta) = \frac12\int_\Xcal |r_t(\theta,\eta,\bfx)|^2 d\nu(\bfx).
\label{eq:NG:ContTime:MinResObj}
\end{equation}
In Eq.~\eqref{eq:NG:ContTime:MinResObj}, $\nu$ is a positive measure with full support on $\Xcal$, which we will discuss in detail later in Section~\ref{sec:Active}. 

The initial condition~$\theta_0$ can be obtained via e.g.~minimization of the least-squares loss between $u_0$ and $U(\theta_0)$:
\begin{equation}
\label{eq:loss:init}
    \theta_0 \in \argmin_{\theta\in\Theta} \int_\Xcal |u_0(\bfx)-U(\theta,\bfx)|^2 d\nu(\bfx).
\end{equation}
Similarly, boundary conditions can be enforced either naturally by choosing $U(\theta(t))$ such that these conditions are satisfied for all $\theta(t)\in\Theta$, or by adding appropriate penalty terms to the residual~\eqref{eq:res} and the objective~\eqref{eq:NG:ContTime:MinResObj}. 
For more justification of Eq.~\eqref{eq:NG:ODE}  note that this equation also arises as a consequence of the proposition introduced in~\ref{sec:deriv:NGODE}.  

\subsection{System of ODEs for network parameters}\label{sec:NG:SystemODEs} The minimizers of $J_t(\theta(t),\eta)$ solve the Euler-Lagrange equation
\begin{equation}
\nabla_{\eta} J_t(\theta(t),\eta) = 0\,.
\label{eq:NG:ContTime:ResMinEquation}
\end{equation}
Written explicitly, Eq.~\eqref{eq:NG:ContTime:ResMinEquation} is a system of ODEs for $\theta(t)$:
\begin{equation}
M(\theta)\dot\theta = F(t,\theta)\,, \qquad \theta(0) = \theta_0\,,
\label{eq:NG:ODE}
\end{equation}
where we defined 
\begin{equation}
\begin{aligned}
     M(\theta) = &   \int_\Xcal \nabla_\theta U(\theta,\bfx)\otimes \nabla_\theta U(\theta,\bfx)d\nu(\bfx)\,, 
    \\
    F(t,\theta) =  & \int_\Xcal \nabla_\theta U(\theta,\bfx) f(t,\bfx,U(\theta))d\nu(\bfx)\,,
\end{aligned}
\label{eq:M:F:def}
\end{equation}
in which $\otimes$ denotes the outer product. Following the same principle, analogous equations as given by \eqref{eq:M:F:def} have been derived in other settings of nonlinear parametrizations as well; see, e.g., the works on dynamic low-rank approximations and related methods \cite{doi:10.1137/050639703,SAPSIS20092347} and the works \cite{UngerTransformModes2020,PhysRevE.104.045303,doi:10.1137/21M1415972}.

The system of ODEs defined in~\eqref{eq:NG:ODE} for $\theta(t)$ is the basis for the Neural Galerkin schemes we introduce below. We refer to $U(\theta(t))$ as Neural Galerkin solution because Eq.~\eqref{eq:NG:ODE} can be derived via testing the residual $r_t(\theta,\dot \theta,\bfx)$ in Eq.~\eqref{eq:res} using  $\nabla_{\theta}U(\theta, x)$ as test function and Galerkin projection with respect to the inner product of $L^2(\Xcal; \nu)$. The objective function~\eqref{eq:NG:ContTime:MinResObj} also arises if we take the inner product between $r_t(\theta,\dot \theta,\bfx)$ and a test function $g(\bfx)$ and maximize it over all $g(\bfx)$ with norm one in $L^2(\Xcal; \nu_{\theta})$.

\subsection{Discretization in time} Eq.~\eqref{eq:NG:ODE} is a system ODEs that can be numerically integrated using standard methods; thus training the network parameters becomes solving a system of ODEs in time.  Let us denote by $\theta^k$ the numerical approximation to $\theta(t_k)$ where the times $t_k$ with $k=0,1,2,\ldots$ are defined recursively via $t_0=0$, $t_{k+1} = t_k + \delta t_k$ using the time steps $\delta t_k>0$, possibly nonuniform and chosen adaptively. To update $\theta^k$, we can either use explicit or implicit integrators.

\subsubsection{Explicit integrators.} The forward Euler scheme is an explicit scheme that leads to 
\begin{equation}
M(\theta^{k})\theta^{k+1} = M(\theta^{k})\theta^{k} + \delta t_k F(t_{k}, \theta^{k})\,.
\label{eq:ODE:DiscExp}
\end{equation}
This is an explicit equation for $\theta^{k+1}$ which is solved at every time step---to this end it may be useful to add a regularizing term $\lambda \, \text{Id}$ to $M$, where $\text{Id}$ is the identity.
 Other integrators can be used as well, such as the adaptive Runge–Kutta–Fehlberg (RK45) method  \cite{DORMAND198019}, which leads to linear systems similar to~\eqref{eq:ODE:DiscExp} at every stage---RK45 will be used in the examples below.

\subsubsection{Implicit integrators.} We consider the backward Euler scheme as an example for implicit integrators, which leads to the system of equations
\begin{equation}
M(\theta^{k+1})\theta^{k+1} =  M(\theta^{k+1})\theta^{k} + \delta t_k  F(t_{k+1}, \theta^{k+1})\,.
\label{eq:ODE:DiscImp}
\end{equation}
Implicit integrators are more costly per time step since they lead to nonlinear equations such as~\eqref{eq:ODE:DiscImp} for $\theta^{k+1}$ that have to be solved at each time step.

We stress that the adaptive nature of our Neural Galerkin approach is different from time-space collocation methods based on neural networks \cite{doi:10.1002/cnm.1640100303,SIRIGNANO20181339,RAISSI2019686} and other DNN-based PDE solution methods that are formulated over equidistant time steps \cite{Han8505}.

Instead of deriving the variational formulation \eqref{eq:min} in continuous time and then discretizing the Euler-Lagrange equations \eqref{eq:NG:ContTime:ResMinEquation}, one could also first discretize in time and then derive the Euler-Lagrange equation of the discrete-in-time optimization problem; this is similar to the distinction between discretize-then-optimize versus optimize-then-discretize approaches used e.g. in inverse problems and PDE-constrained optimization: see e.g. Ref.~\cite{doi:10.1137/1.9780898718720}.

\section{Active learning with Neural Galerkin schemes}\label{sec:Active}
We now introduce active learning for Neural Galerkin schemes to efficiently estimate $M(\theta)$ and $F(t, \theta)$. To achieve this, we make the measure in the objective \eqref{eq:NG:ContTime:MinResObj} and the definition of $M$ and $F$ depend on the parameter $\theta(t)$ so that we can adapt the measure $\nu_{\theta(t)}$ over time; see Section~\ref{sec:Active:Formulation}. We then formulate adaptive Monte Carlo estimators of $M$ and $F$ with the aim of achieving lower mean-squared errors than regular Monte Carlo estimations that are based on un-informed, e.g., uniform, sampling. In Section~\ref{sec:Active:AdaptMeasure}, we propose to directly use a time-adaptive measure and in Section~\ref{sec:Active:IS} we show that Neural Galerkin schemes can be combined with importance sampling. Section~\ref{sec:Active:Arch} specifically considers shallow networks with Gaussian units that gives rise to an adaptive sampling scheme that leverages the network architecture.

\subsection{Active learning via time-dependent measures}\label{sec:Active:Formulation}
For generic choices of parametrization~$U(\theta)$, the integrals in Eq.~\eqref{eq:M:F:def} do not admit a closed-form solution and so
will need to be numerically estimated.  In low dimensions, this calculation can be performed by quadrature on a grid. When $\Xcal$ is high dimensional, however, quadrature becomes computationally intractable quickly and we must proceed differently. If the measure $\nu$ against which is integrated in Eq.~\eqref{eq:M:F:def} is a probability measure, we can consider using a vanilla Monte Carlo estimator for each term, by drawing $n$ samples $\{\bfx_i\}_{i=1}^n$ from $\nu$ and replacing the expectations by empirical averages over these samples.  This estimator is efficient to approximate certain kernels uniformly over high-dimensional spaces \cite{NIPS2007_013a006f}, but not necessarily if the solution to the PDE~\eqref{eq:nonlinear:PDE} develops spatially localized structures, as is the case in many physical systems of interest.

Instead, we consider the objective $J_t$ defined in \eqref{eq:NG:ContTime:MinResObj} and reformulate it with a parameter-dependent measure $\nu_{\theta}$ as
\begin{equation}\label{eq:Active:NewJt}
J_{\theta}(\theta, \eta) = \frac{1}{2}\int_{\Xcal} |r_t(\theta, \eta, \bfx)|^2\mathrm d\nu_{\theta}(\bfx)\,.
\end{equation}
The key difference to $J_t$ used in Section~\ref{sec:NG:DiracFrenkel}, is that now the objective defined in \eqref{eq:Active:NewJt} integrates the squared residual against a measure $\nu_{\theta}$ that depends on the parameter $\theta(t)$ and thus it depends on time $t$. Following an analogous derivation as in Section~\ref{sec:NG:SystemODEs}, we obtain the operators
\begin{equation}
\begin{aligned}
     M(\theta) = &   \int_\Xcal \nabla_\theta U(\theta,\bfx)\otimes \nabla_\theta U(\theta,\bfx)d\nu_{\theta}(\bfx) 
    \\
    F(t,\theta) =  & \int_\Xcal \nabla_\theta U(\theta,\bfx) f(t,\bfx,U(\theta))d\nu_{\theta}(\bfx)\,,
\end{aligned}
\label{eq:M:F:defAdaptMeasure}
\end{equation}
with measure $\nu_{\theta}$ depending on $\theta$. This new formulation with measure $\nu_{\theta}$ allows us to adapt the measure so that Monte Carlo estimators of $M$ and $F$ have a low mean-squared error.

\subsection{Monte Carlo estimators with adaptive samples}\label{sec:Active:AdaptMeasure} Let $\{\bfx_i(t)\}_{i = 1}^n$ be $n$ samples of $\nu_{\theta(t)}$ and notice that the samples $\bfx_1(t), \dots, \bfx_n(t)$ depend on time $t$. Monte Carlo estimators $\tilde{M}$ and $\tilde{F}$ of $M$ and $F$, respectively, are then given by
\begin{equation}
    \begin{aligned}
    \tilde M(\theta(t)) = &  \frac1n \sum_{i=1}^n\nabla_\theta U(\theta(t),\bfx_i(t))\otimes \nabla_\theta U(\theta(t),\bfx_i(t)),\\
 \tilde F(t,\theta(t)) = &  \frac1n \sum_{i=1}^n \nabla_\theta U(\theta(t),\bfx_i(t)) f(t,\bfx_i(t),U(\theta(t)))\,.
\end{aligned}
\label{eq:F:M:emp:2}
\end{equation}
For these estimators to be effective, their variances must be small, which requires choosing the measure $\nu_{\theta(t)}$ in a way that is informed by the solution~$U(\theta(t))$. This choice must be made on a case-by-case basis. For example, below we consider Fokker-Planck equations whose solutions are normalized probability density functions: in this case, it is reasonable to take $\nu_{\theta} = U(\theta)$---this choice leads to low variance estimators if the density has mass localized in a few relatively small regions that move in $\Xcal$. If in addition the parametrization $U(\theta)$ involves Gaussian units, as in the Fokker-Planck examples below, sampling from it is straightforward---with this specific choice we can can in fact evaluate analytically the integrals of the components of $M(\theta)$ and sampling is only required to estimate $F(t,\theta)$. Other choices for $\nu_{\theta}$ are possible too, e.g. by taking $\nu_{\theta(t)}$ proportional to $|\nabla_{\bfx} U(\theta(t),\bfx)|$ in situations where the solution has localized fronts that propagate in $\Xcal$.

\subsection{Importance sampling with respect to a nominal measure}\label{sec:Active:IS} Consider a nominal measure $\nu$. Then, an adaptive measure $\nu_{\theta}$ can be realized via importance sampling by defining positive weights $\omega: \Xcal\times \Theta \to(0,\infty)$ such that $Z_\theta =\int_\Xcal\omega(\bfx,\theta)\mathrm d\nu(\bfx) < \infty$ for all $\theta\in \Theta$. Notice that the weights given by $\omega$ depend on the parameter $\theta(t)$ and thus on time $t$. One then draws a set of samples $\{\bfx_i(t)\}_{i=1}^n$ from $d\nu_{\theta(t)}(\bfx) = Z_\theta^{-1} \omega(\bfx,\theta(t)) d\nu(\bfx)$ and uses the
reweighted estimators $\tilde{M}$ and $\tilde{F}$ given by
\begin{equation}
    \begin{aligned}
    \tilde M(\theta) = &  \frac1n \sum_{i=1}^n\frac{\nabla_\theta U(\theta,\bfx_i)\otimes \nabla_\theta U(\theta,\bfx_i)}{ \omega(\bfx_i,\theta)},\\
 \tilde F(t,\theta) = &  \frac1n \sum_{i=1}^n \frac{\nabla_\theta U(\theta,\bfx_i) f(t,\bfx_i,U(\theta))}{\omega(\bfx_i,\theta)}\,.
\end{aligned}
\label{eq:F:M:emp:1}
\end{equation}

\subsection{Neural architectures and sampling}\label{sec:Active:Arch} Our method offers the possibility to use DNNs of arbitrary complexity and sophistication for the parametrization $U(\theta,\bfx)$, in  a way that can be tailored to known properties or symmetries of the solution of Eq.~\eqref{eq:nonlinear:PDE}, as well as its boundary conditions.
However, the network architecture can also be chosen such that it aligns well with the sampling from the adaptive measure $\nu_{\theta(t)}$. In later numerical experiments, we will consider shallow (one-hidden-layer) network with $m$ nodes given by
\begin{equation}
U(\theta, \bfx) = \sum\nolimits_{i = 1}^m c_i\varphi(\bfx, w_i, \bfb_i)\,,
\label{eq:NumExp:Exp:ShallowNetwork}
\end{equation}
where $\theta$ is defined as $\theta = (c_i,w_i,\bfb_i)_{i=1}^m$ and $c_i,w_i \in \mathbb{R}$ and $\bfb_i \in \mathbb{R}^{d}$ are the parameters of the DNN, and $\varphi:\Xcal\times\R^{d+1} \to \mathbb{R}$ is the following nonlinear unit (activation function) given by the Gaussian kernel
\begin{equation}
\varphi_{\text{G}}(\bfx, w, \bfb) = \exp\left(-w^2|\bfx - \bfb|^2\right)\,.
\label{eq:NumExp:Exp:ExpUnit}
\end{equation}
Here and below the exponential function acts componentwise, and we  refer to the $c_i$ as the coefficients of the network and to the $w_i, \bfb_i$ as its features.

The network architecture \eqref{eq:NumExp:Exp:ShallowNetwork} is especially useful for approximating solutions that are probability density functions, such as for the Fokker-Planck equation considered below. In such a case, as discussed in Section~\ref{sec:Active:AdaptMeasure}, it is appropriate to set $\nu_{\theta(t)} = U(\theta(t))$. Then, the set of samples $\{(\bfx_i(t)\}_{i = 1}^n$ at time $t$ is obtained via sampling from $U(\theta)$ directly, which can be done efficiently by exploiting that the units \eqref{eq:NumExp:Exp:ExpUnit} are Gaussians.

\section{Numerical experiments}
We demonstrate Neural Galerkin schemes with active learning on several numerical examples. We first discuss network architectures in Section~\ref{sec:NumExp:NN} and then consider numerical examples reaching from low-dimensional benchmark problems in Section~\ref{sec:NumExp:KdV} and Section~\ref{sec:NumExp:AC} to PDEs whose solutions develop local features in moderately high-dimensional spatial domains in Section~\ref{sec:NumExp:Adv} and Section~\ref{sec:NumExp:Particle}.

\subsection{Network architectures}\label{sec:NumExp:NN}
We use two specific neural architectures. The first is a shallow (one-hidden-layer) network with $m$ nodes given by \eqref{eq:NumExp:Exp:ShallowNetwork}. The unit $\varphi$ is one of the following two nonlinear units (activation function): The first is the Gaussian kernel
defined in \eqref{eq:NumExp:Exp:ExpUnit} and the second one is
\begin{equation}
\varphi_{\text{G}}^L(\bfx, w, \bfb) = \exp\left(-w^2|\sin\left(\pi(\bfx - \bfb)/L\right)|^2\right)\,,
\label{eq:NumExp:Exp:ExpLUnit}
\end{equation}
which we use when $\Xcal=L\mathbb{T}^d$ with $L>0$ and we need to enforce periodicity. The sine function acts componentwise on its argument. 
We  refer to the $c_i$ as the coefficients of the network and to the $w_i, \bfb_i$ as its features.

The other neural architecture that we use is a feedforward neural network with $\ell \in \mathbb{N}$ hidden layers and $m$ nodes per layer:
\begin{multline}
U(\theta, \bfx) = \bfc(t)^T\tanh(\bfW_{\ell}\tanh(\bfW_{\ell - 1}(\cdots
\varphi^L_{\tanh}(\bfx, \bfW_1, \bfb_1) \cdots ) + \bfb_{\ell-1}) + \bfb_\ell)\,,
\label{eq:NumExp:Tanh:DeepNetwork}
\end{multline}
where $\theta = (\bfc, \bfW_1, \dots, \bfW_{\ell}, \bfb_1, \dots, \bfb_{\ell})$ with $\bfc\in \R^m$, $\bfW_1 \in \mathbb{R}^{m \times d}, \bfW_i \in \mathbb{R}^{m \times m}, i > 1$, and $\bfb_i \in \mathbb{R}^d$ are the parameters, and $\varphi^L_{\tanh}$ is the nonlinear unit
\begin{equation}
\varphi_{\tanh}^L(\bfx, \bfW, \bfb) = \tanh\left(\bfW\sin\left(2\pi\left(\bfx - \bfb\right)/L\right)\right)\,.
\end{equation}
The tanh function acts compontwise on its argument.

\begin{figure}[t]
\begin{tabular}{ll}
\includegraphics[width=0.48735\columnwidth]{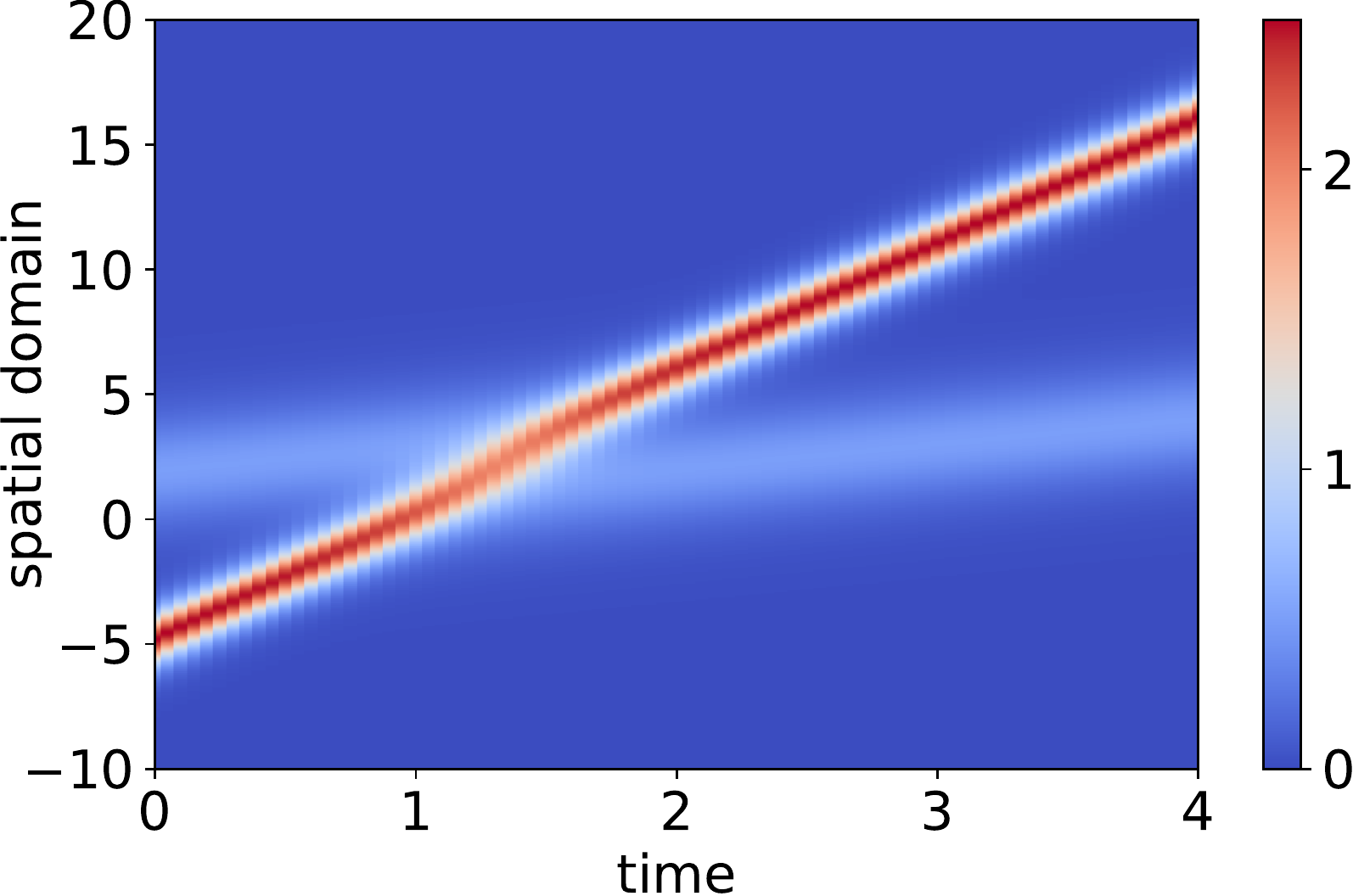}& \includegraphics[width=0.4275\columnwidth]{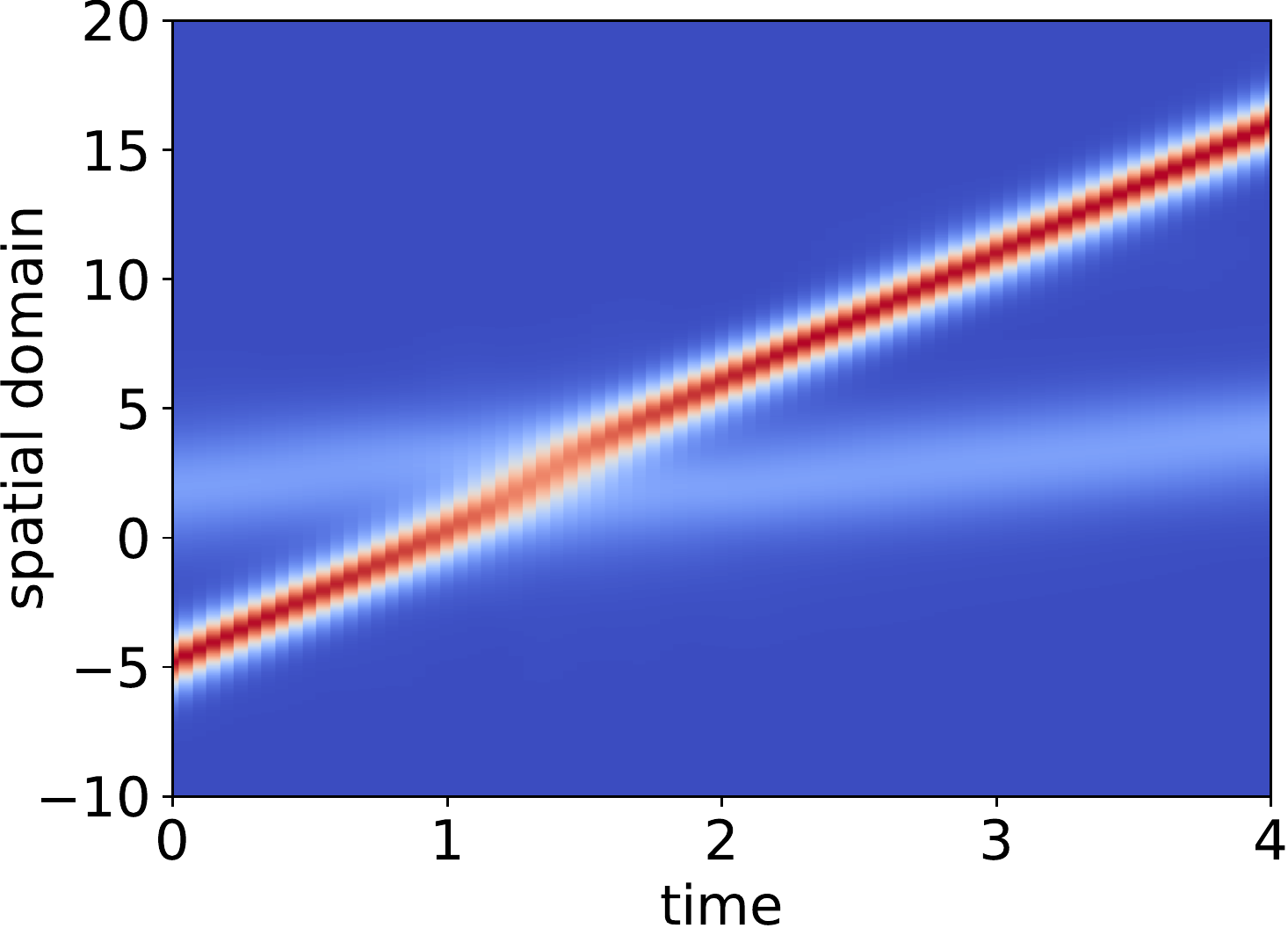}\\
\multicolumn{1}{c}{\footnotesize (a) truth} & \multicolumn{1}{c}{\footnotesize (b) Neural Galerkin} \\
\includegraphics[width=0.4275\columnwidth]{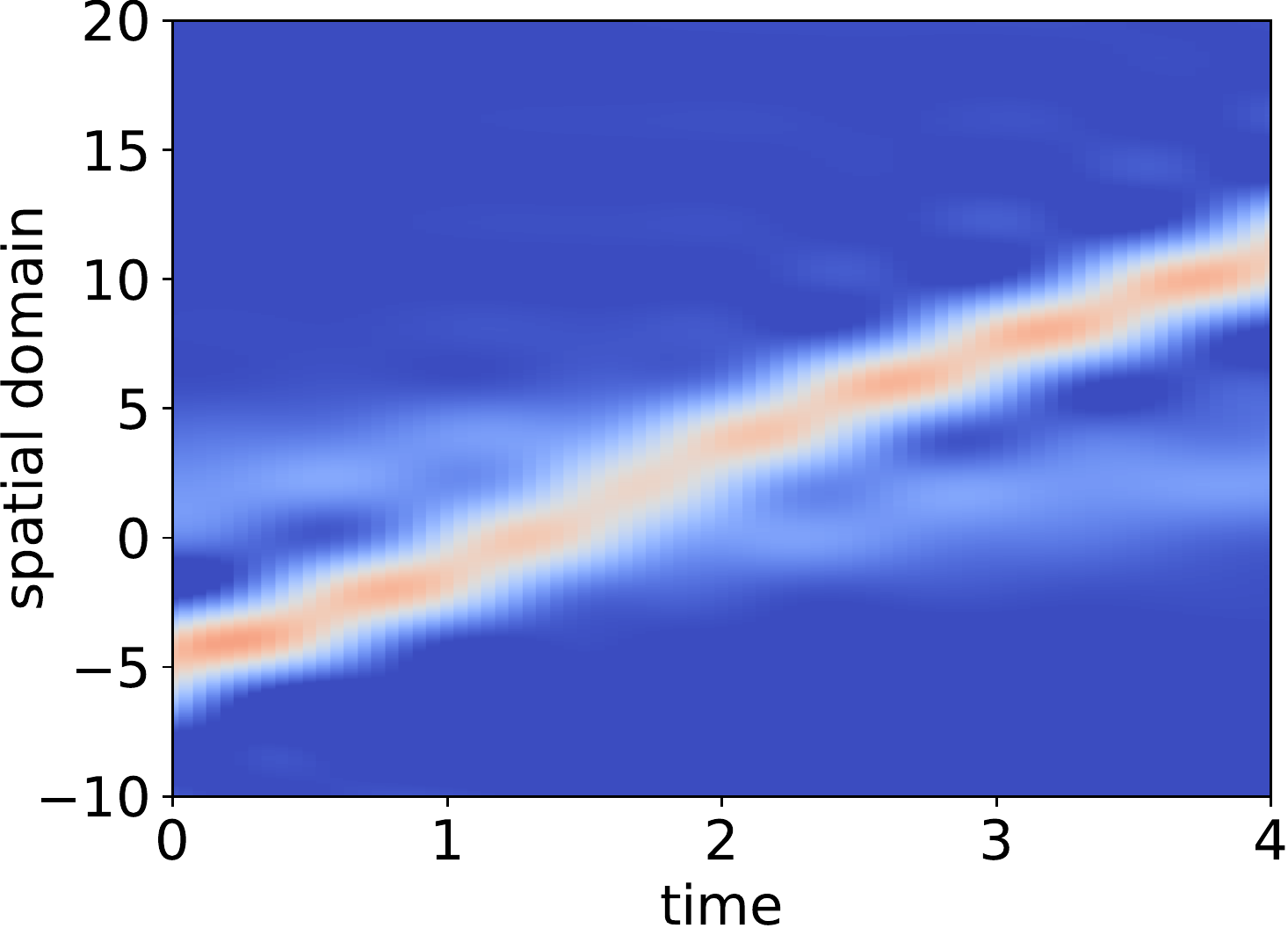} &
\includegraphics[width=0.4275\columnwidth]{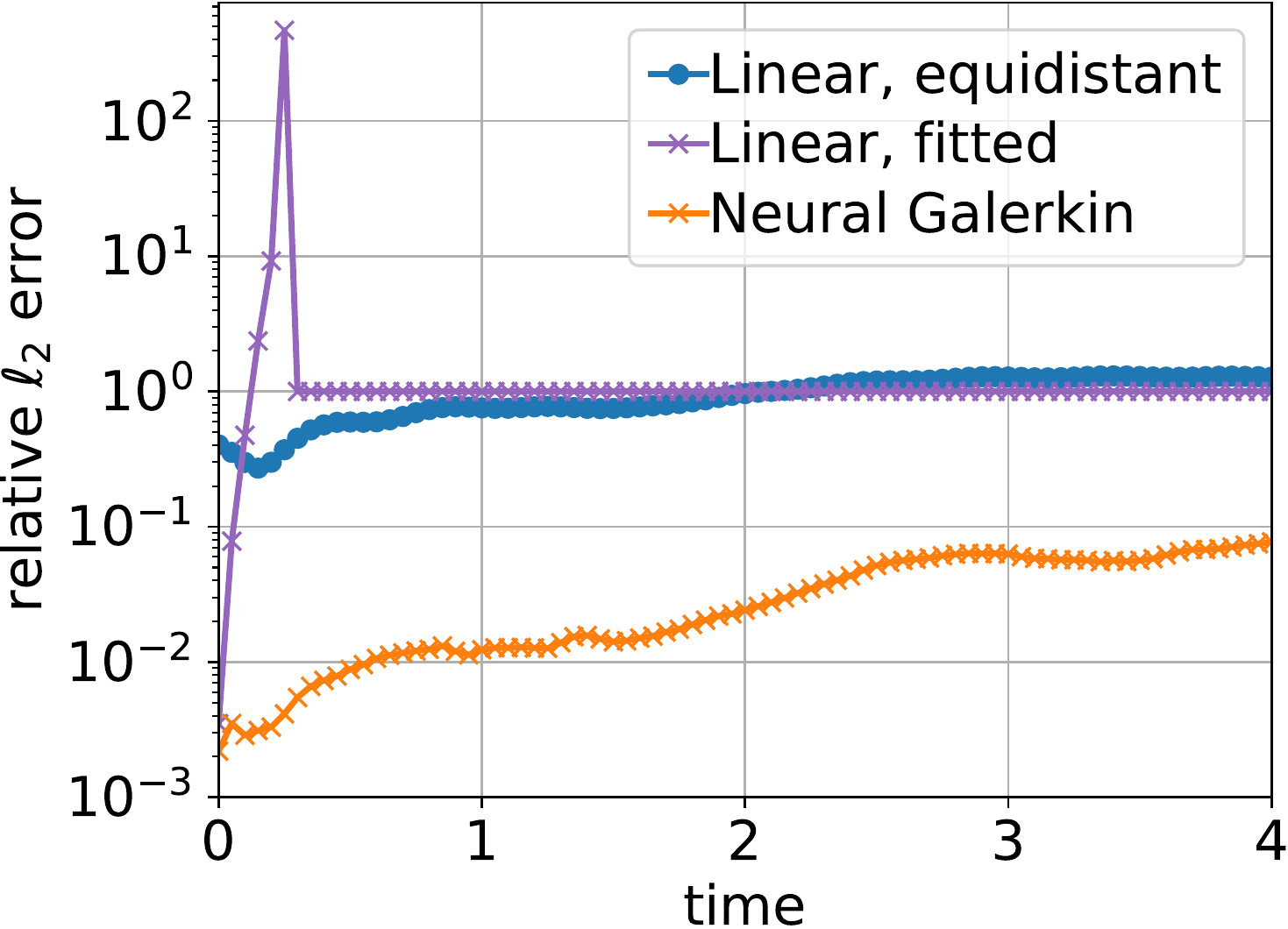}\\
\multicolumn{1}{c}{\footnotesize (c) Linear Galerkin, equidistant}  & \multicolumn{1}{c}{\footnotesize (d) error over time}
\end{tabular}
\caption{Korteweg-De Vries: Neural Galerkin schemes integrate in time a nonlinear parametrization of the solution and so obtain an accurate approximation with few degrees of freedom in this example. In contrast, linear (standard) Galerkin that derives approximations in fixed spaces with fixed bases without feature adaptation leads to poor approximations when there are local dynamics in the spatial domain.} 
\label{fig:KdVSpaceTimePlots}
\end{figure}

\subsection{Korteweg-de Vries equation}\label{sec:NumExp:KdV}
Consider the KdV equation
$
\partial_t u + \partial_x^3 u + 6 u \partial_x u = 0\,
$
on the one-dimensional spatial domain $\Xcal = [-20, 40) \subset \mathbb{R}$ with periodic boundary condition.  We follow the setup described in Ref.~\cite{TAHA1984231}, which for the initial conditions described in their Section~I-a-(ii) leads to a solution of the KdV equation that consists of two interacting solitons that approach each other, collide, and then separate---this solution, shown in Figure~\ref{fig:KdVSpaceTimePlots}(a), is available analytically, thereby providing us with a benchmark of the numerical results.

For the parametric solution, we use the shallow neural network defined in Eq.~\eqref{eq:NumExp:Exp:ShallowNetwork}  with the Gaussian units defined in Eq.~\eqref{eq:NumExp:Exp:ExpLUnit} and $m = 10$ nodes. Because the dimension of the spatial domain is one in this example, we simply take $d\nu_\theta(x) = dx$, and sample $n = 1000$ points uniformily in $[0,1]$ to estimate $M$ and $F$.  
The initial parameter $\theta_0$ is obtained via a least-squares fit of the initial condition. 
The batch size is $10^5$ and number of iterations is $10^5$. Samples are drawn uniformly in the spatial domain and the learning rate is $10^{-1}$. We take five replicates with randomized initialization of parameters and then use the fit with lowest error on test samples.
As integrator, we use RK45 so that the time-step size is adaptively chosen based on the dynamics of the solution.

The time-space plot of the Neural Galerkin approximation shown in Figure~\ref{fig:KdVSpaceTimePlots}(b) is in close agreement with the analytic solution shown in Figure~\ref{fig:KdVSpaceTimePlots}(a). The relative $\ell_2$ error remains low over time, as shown in Figure~\ref{fig:KdVSpaceTimePlots}(d). The relative $\ell_2$ error is computed as follows: Consider the $w = 2048$ equidistant grid points $x_1, \dots, x_w$ in $\Xcal$ and define
\[
\bfu(t) = [u(t, x_1), \dots, u(t, x_w)]^T \in \mathbb{R}^w\,,
\]
to be the vector of the analytic solution $u$ at time $t$ and at the grid points $x_1, \dots, x_w$. Similarly, define $\tilde{\bfu}(t) = [\tilde{u}(t, x_1), \dots, \tilde{u}(t, x_w)]^T$ as the solution vector corresponding to an approximation $\tilde{u}$. Then, the reported relative $\ell_2$ error is
\[
e_{\ell_2} = \frac{\sum_{k = 0}^K\|\bfu(t_k) - \tilde{\bfu}(t_k)\|_2^2}{\sum_{k = 0}^K\|\bfu(t_k)\|_2^2}\,.
\]

\begin{figure}
\centering\begin{tabular}{cc}
\includegraphics[width=0.40\columnwidth]{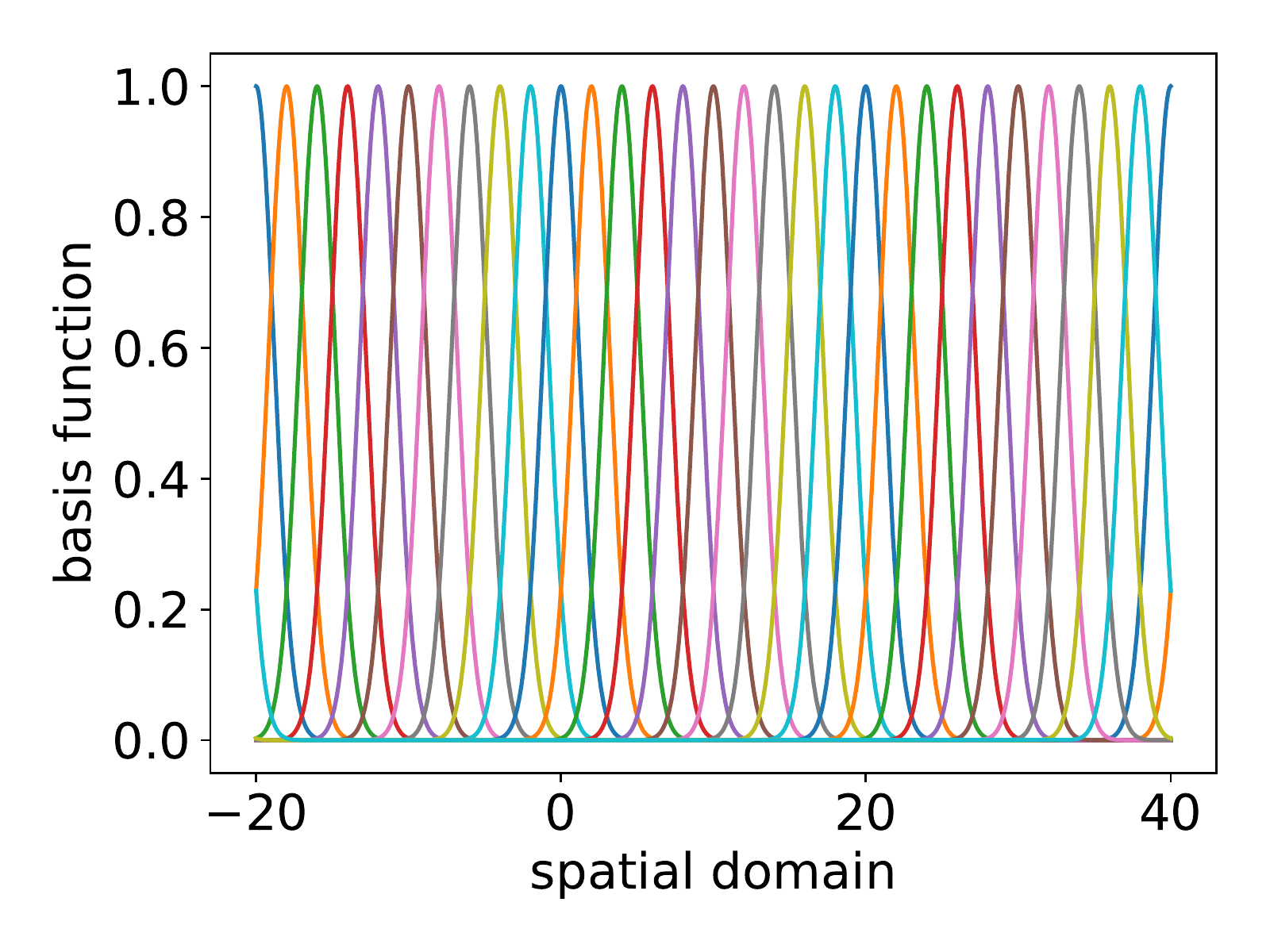} & \includegraphics[width=0.40\columnwidth]{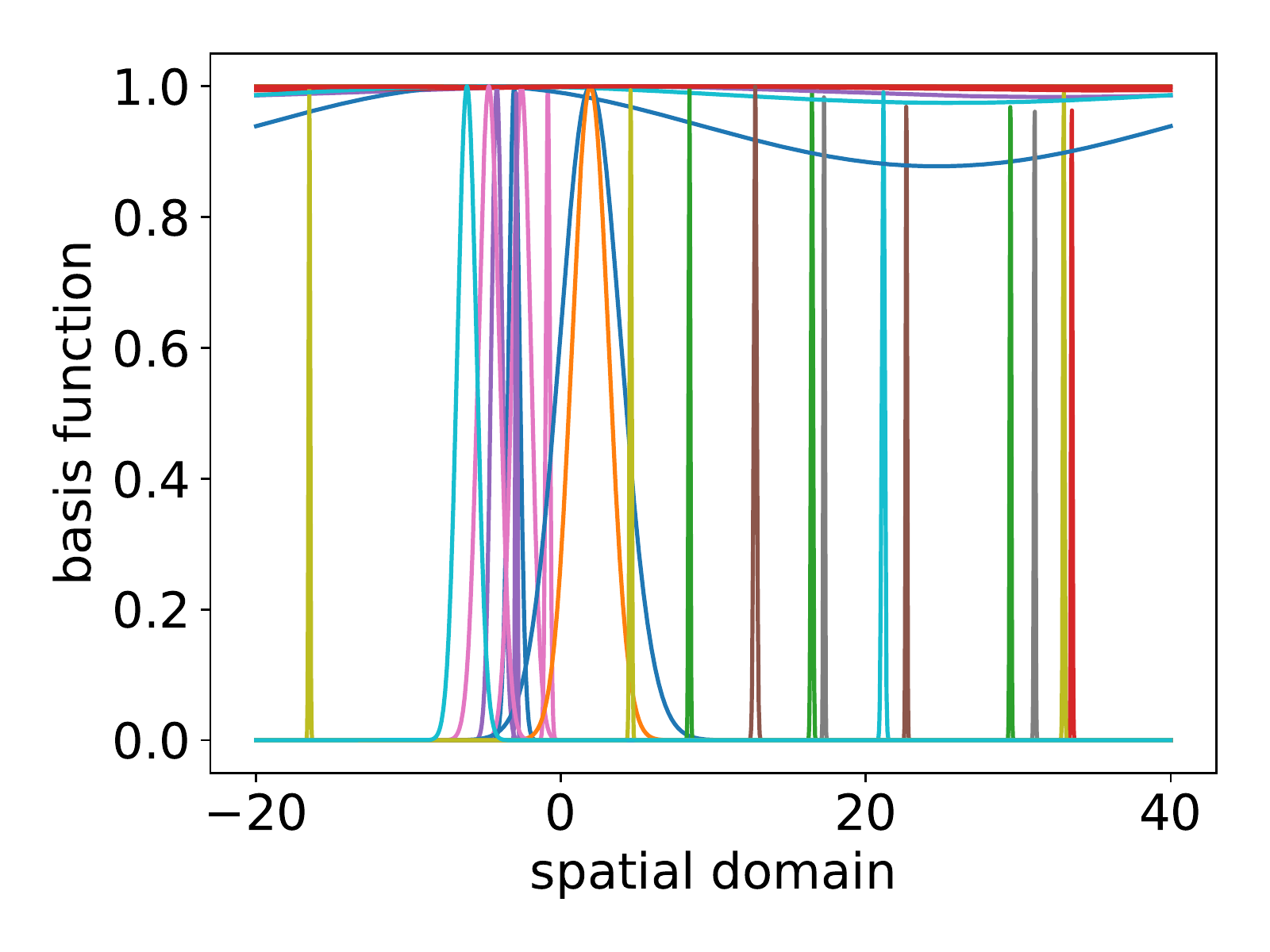}\\ 
\scriptsize (a) KdV, equidistant basis functions & \scriptsize (b) KdV, fitted basis functions 
\end{tabular}
\caption{KdV: Basis functions used for linear Galerkin approximations.}
\label{supp:fig:BasisFun}
\end{figure}

To emphasize the importance of the nonlinear feature adaptation, we also consider linear approximations with network Eq.~\eqref{eq:NumExp:Exp:ShallowNetwork} with $m=30$ nodes but with fixed features $w_i$, and $\bfb_i$ that are either located equidistantly in $\Xcal$ or fitted to the initial condition and then kept fixed. For linear Galerkin with equidistant basis functions, the basis functions are Gaussians units $\varphi_{G}^L$ with a fixed bandwidth and equidistantly located in the domain $\Omega$, see Figure~\ref{supp:fig:BasisFun}a. For the comparison with linear Galerkin with basis functions fitted to the initial condition, we use the units fitted to the initial condition as described above and keep the features fixed, see Figure~\ref{supp:fig:BasisFun}b. In both cases, the linear approximation has the same number of degrees of freedom as the nonlinear Neural Galerkin approximation.

The linear approximation leads to large errors that develop after only a few time steps, as apparent from the solution shown in Figure~\ref{fig:KdVSpaceTimePlots}(c). This is also in agreement with Figure~\ref{fig:KdVSpaceTimePlots}(d) which shows that the linear Galerkin approximations based on fixed features lead to errors that are orders of magnitude larger than the Neural Galerkin approximation with the same number of degrees of freedom.

\begin{figure}[t]
\begin{tabular}{lc}
\includegraphics[width=0.48735\columnwidth]{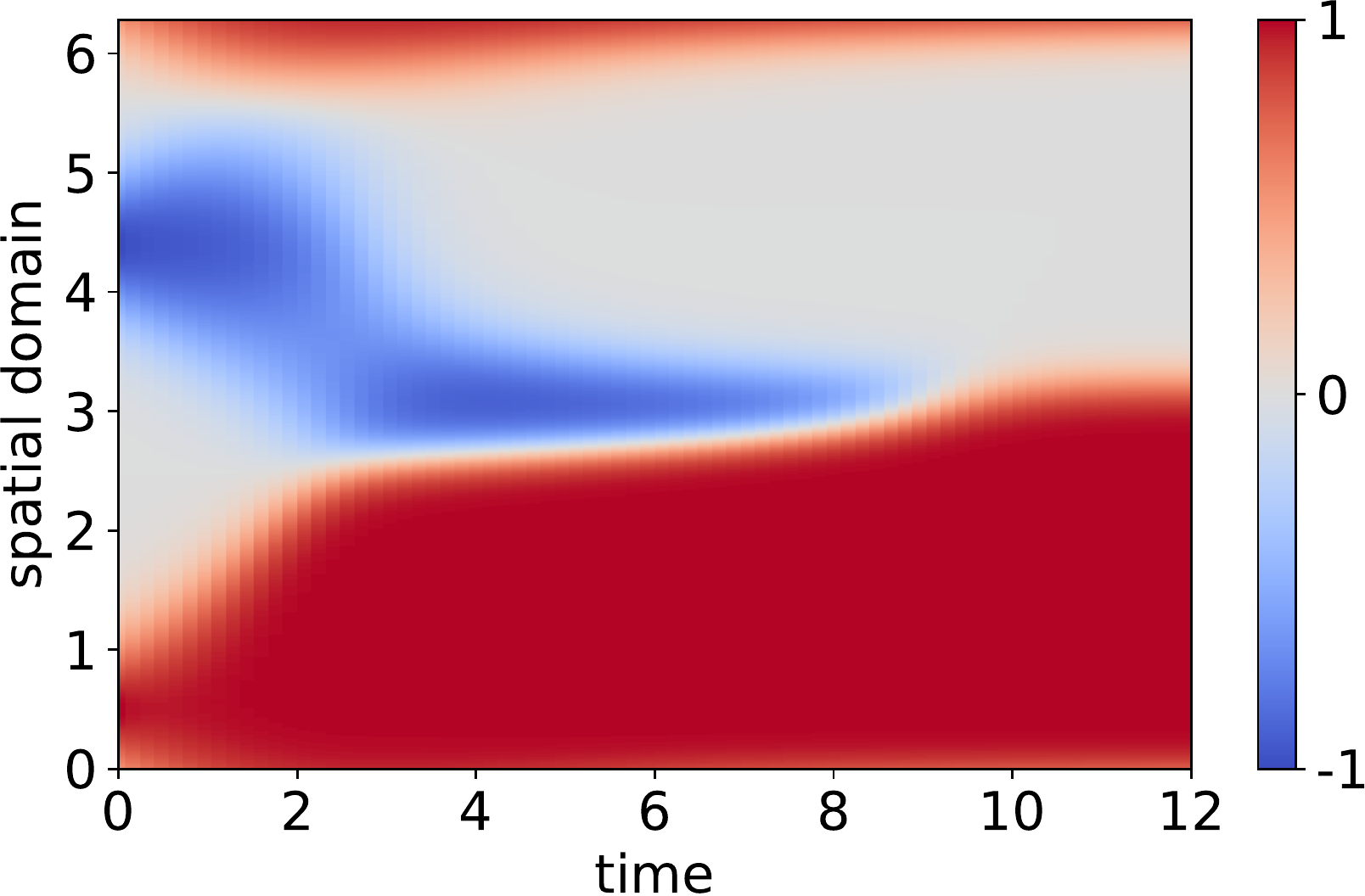} &
\includegraphics[width=0.4275\columnwidth]{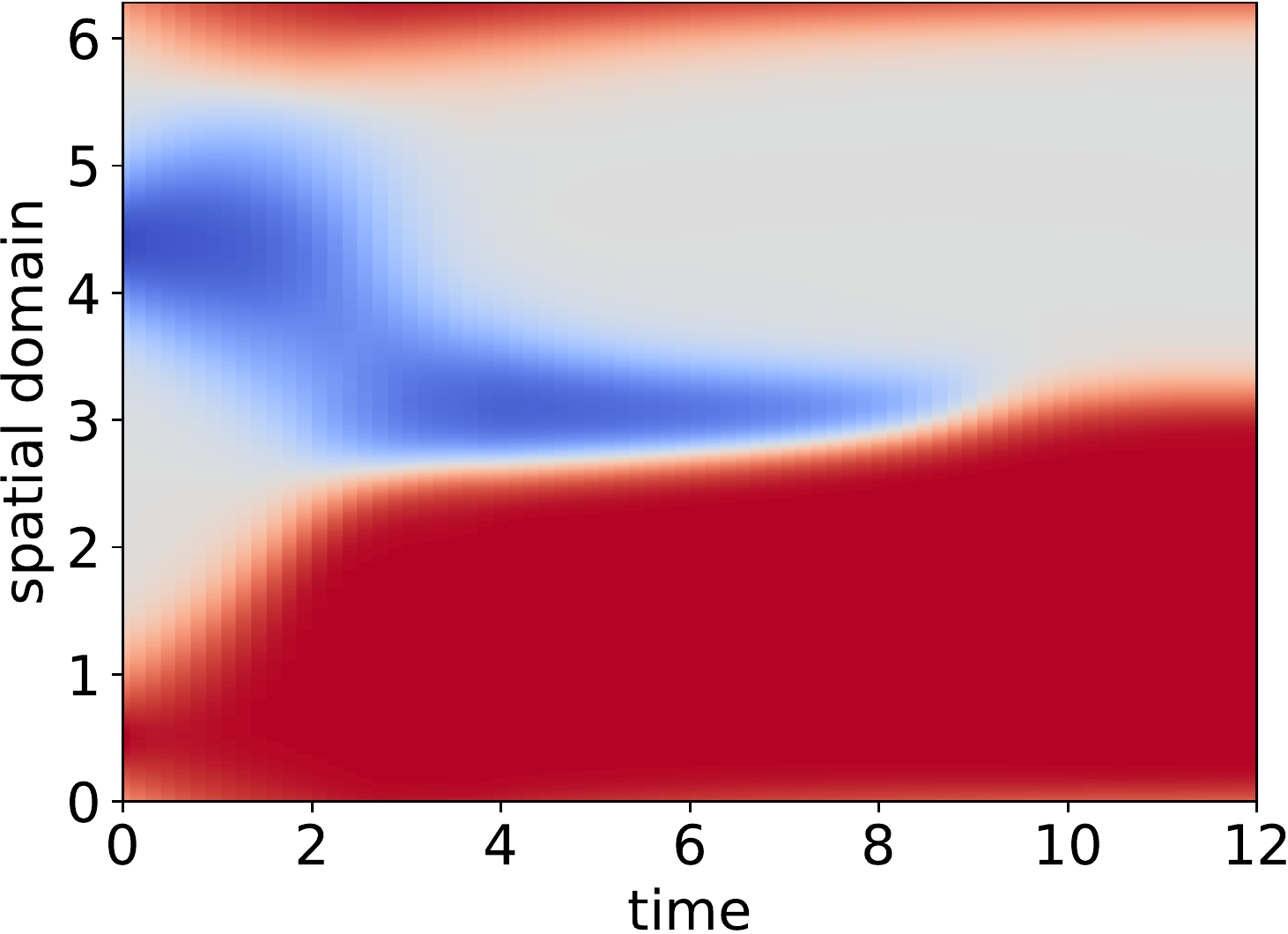} \\
\multicolumn{1}{c}{\footnotesize (a) benchmark} & \multicolumn{1}{c}{\footnotesize (b) Neural Galerkin}\\
\includegraphics[width=0.4275\columnwidth]{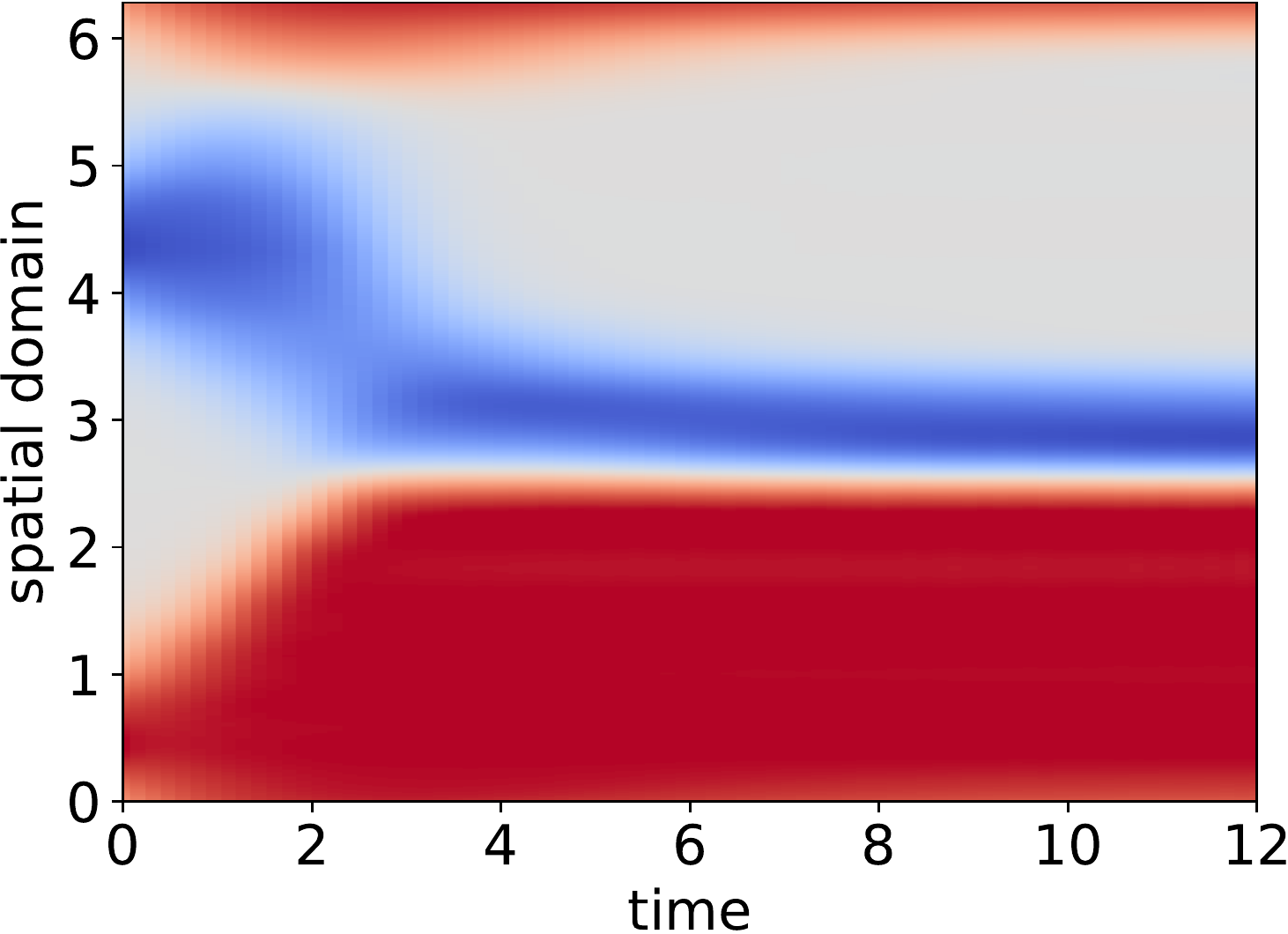}
& \includegraphics[width=0.4275\columnwidth]{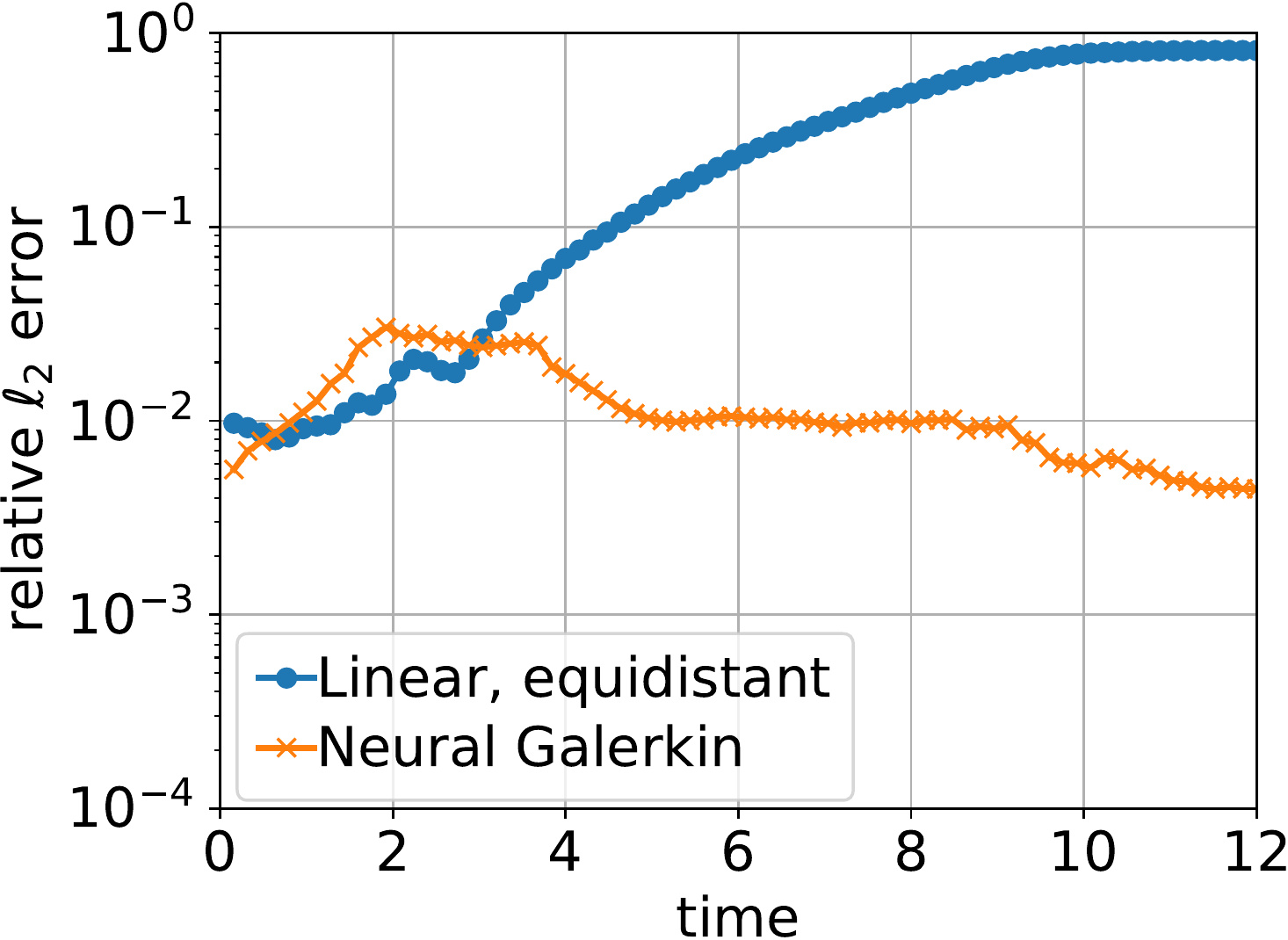}\\
\multicolumn{1}{c}{\footnotesize (c) Linear Galerkin, equidistant} & \multicolumn{1}{c}{\footnotesize (d) error over time}
\end{tabular}
\caption{Allen-Cahn: The proposed Neural Galerkin approach with a three-hidden-layer network correctly predicts the dynamics of the sharp walls that eventually separate the solution into two flat pieces at $0$ and $+1$ in this experiment. In contrast, linear Galerkin with the same number of degrees of freedom wrongly predicts a solution with three flat pieces at $-1, 0$, and $+1$ because of the lack of expressiveness of the linear approximation.}
\label{fig:ACSpaceTimePlots}
\end{figure}

\begin{figure}[!h]
\begin{tabular}{lc}
\includegraphics[width=0.48735\columnwidth]{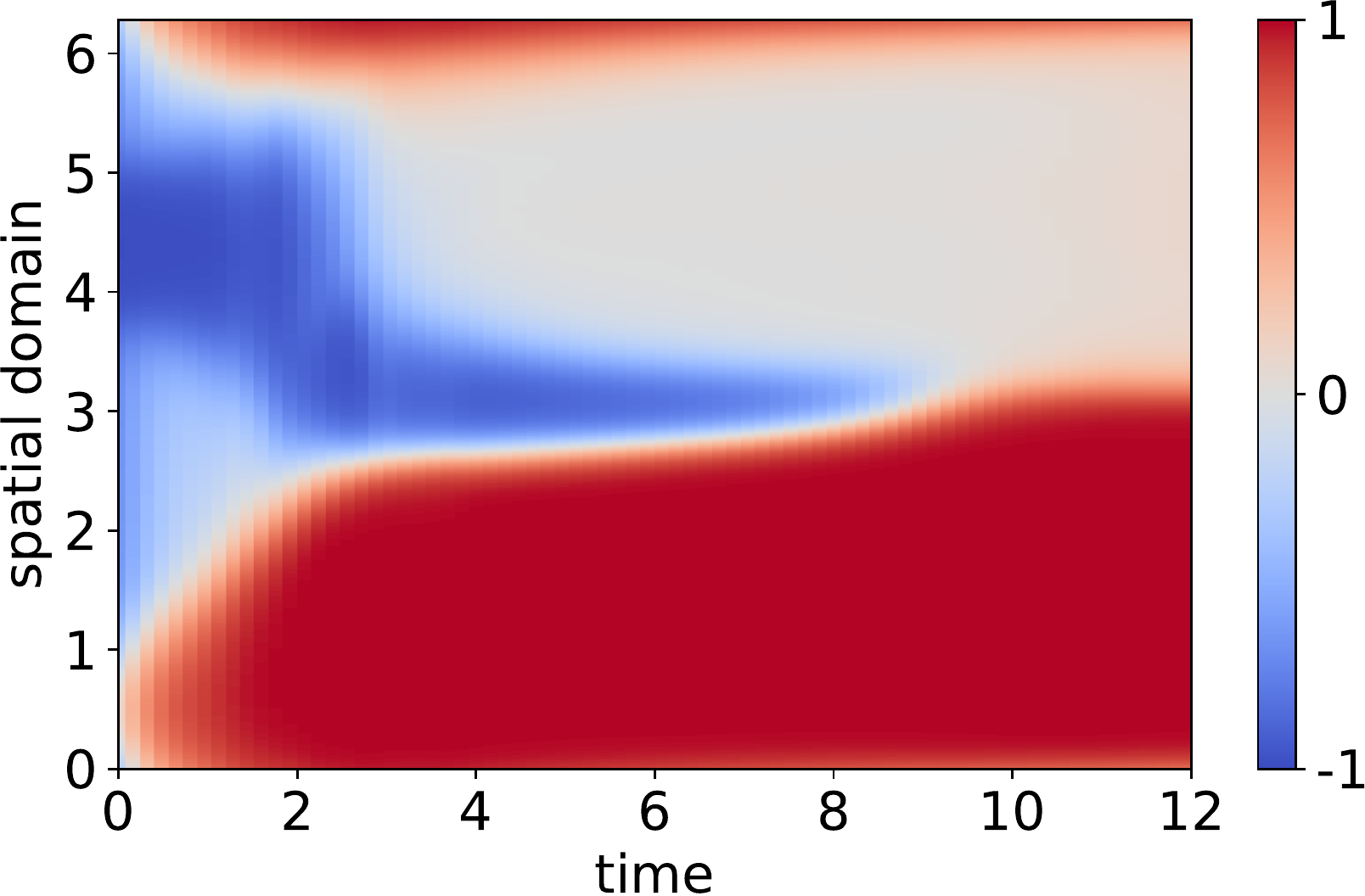} & \includegraphics[width=0.48735\columnwidth]{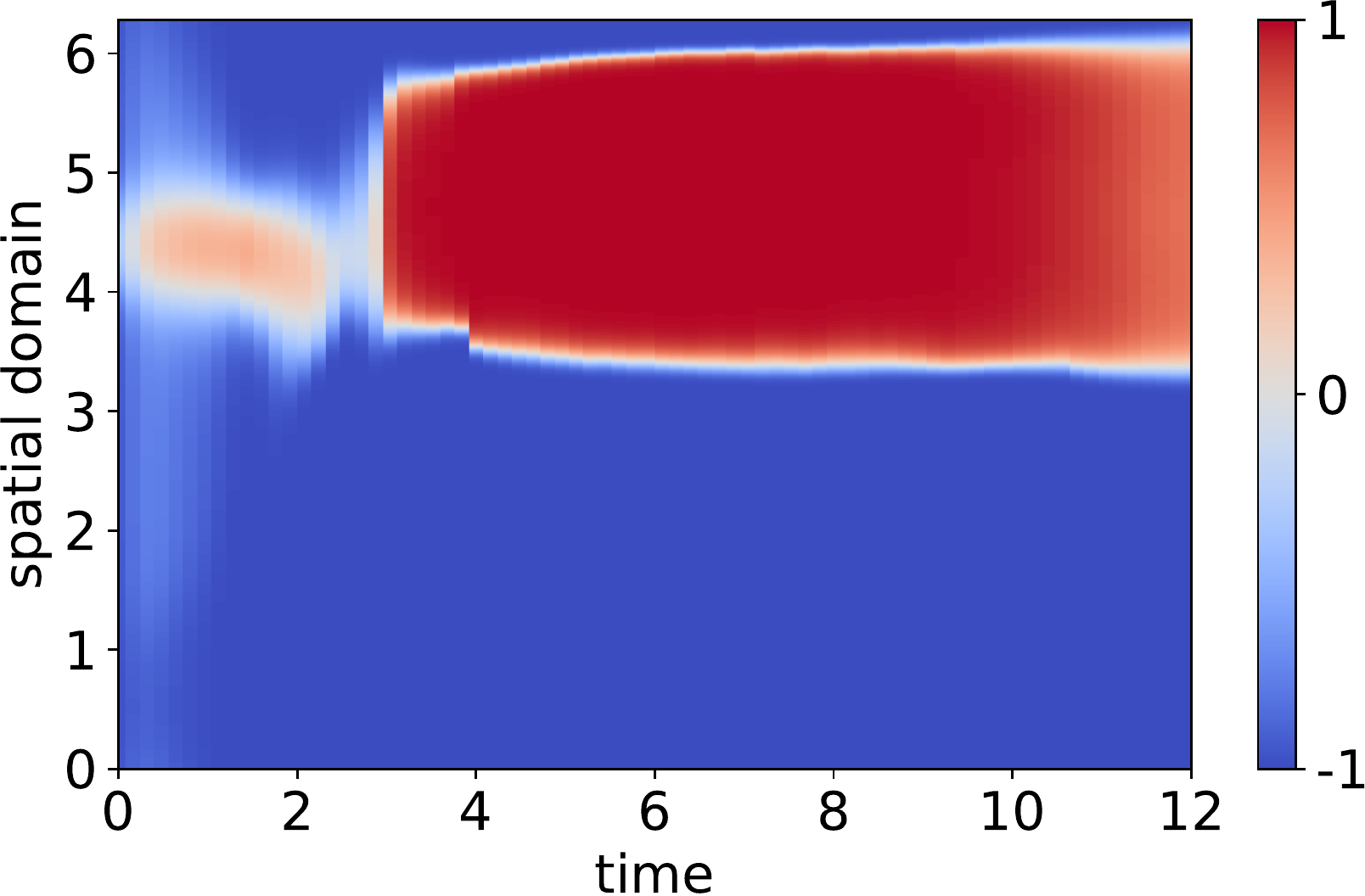}\\
\multicolumn{1}{c}{\footnotesize (a) output component one} & \footnotesize (b) output component two\\
\includegraphics[width=0.44\columnwidth]{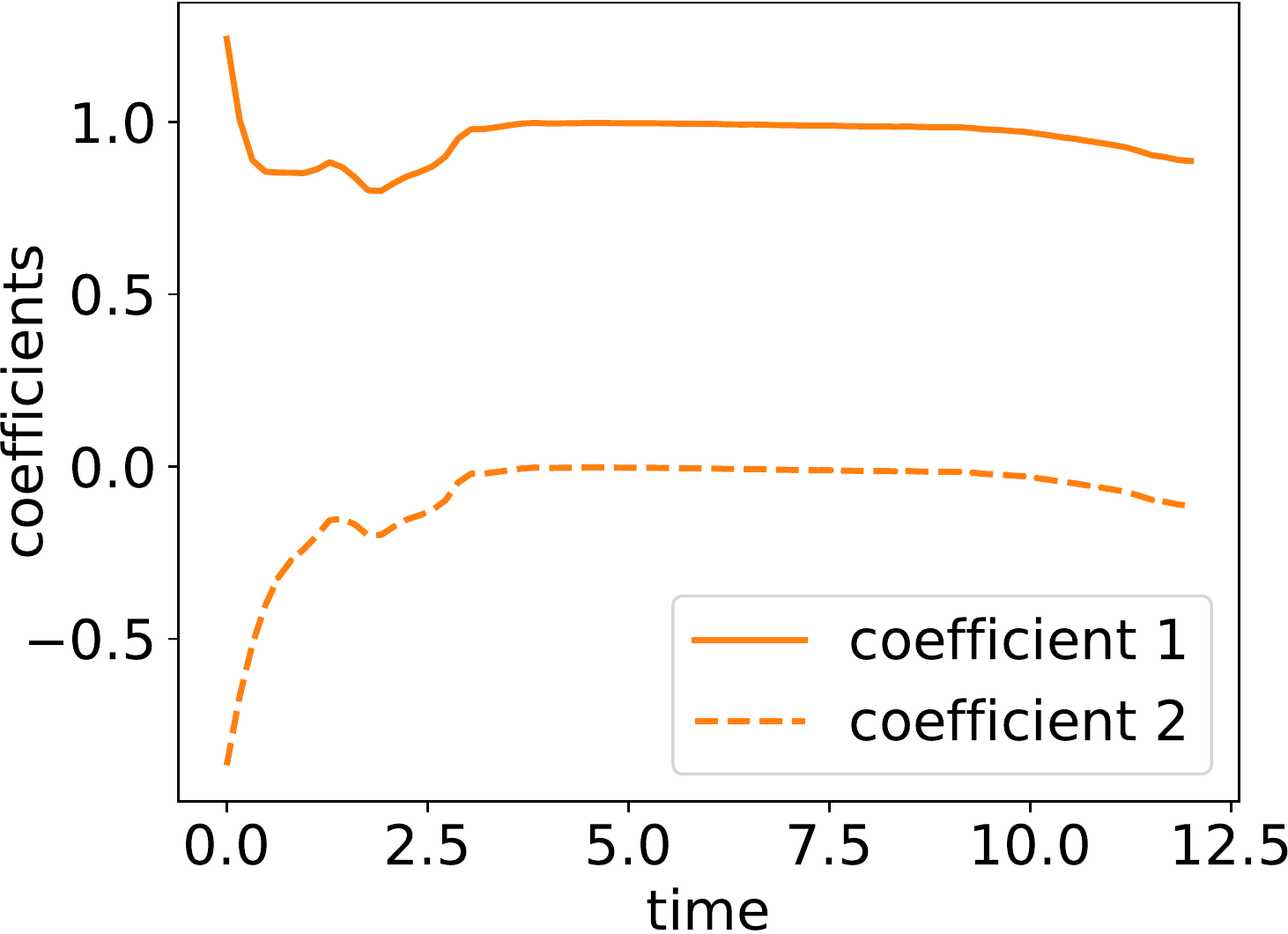} & \includegraphics[width=0.48735\columnwidth]{figures/F2_ChafeeInfante_RBFp_BwdE_adamSGD_uni_zero_ZIfixedFeatures-uni-1e+00-1e+00-1e+00-1_N16_nrL1_dt0.01_bs1000-1000_nrIter10000_h0.1_Ubenchmark-crop}  \\
\multicolumn{1}{c}{\footnotesize (c) coefficients} & \footnotesize (d) benchmark\\
\end{tabular}
\caption{Allen-Cahn: The plots visualize how the used Neural Galerkin scheme propagates forward in time the coefficients and features of a DNN parametrization by following the dynamics prescribed by the PDE. Notice that the output components, defined in Eq.~\eqref{eq:ACOutputComponent} and shown in (a) and (b), reflect the transition of the sharp walls between the flat pieces of the solution and that the dynamics rapidly change from time $t = 1.05$ to about $t = 2.5$ when the potential changes sign in parts of the spatial domain.}
\label{fig:ACOutputsCoeff}
\end{figure}

\subsection{Allen-Cahn (AC) equation}\label{sec:NumExp:AC}
Consider the prototypical reaction diffusion equation known as the AC equation
\[
\partial_t u = \epsilon\partial_x^2 u - a(t, x)(u - u^3)
\]
in the one-dimensional domain $\Xcal= [0, 2\pi)$
with periodic boundary conditions. We set $\epsilon = 5 \times 10^{-2}$ in the diffusion term and let the coefficient in the reaction term vary in time and space using $a(t, x) = 1.05 + t\sin(x)$. This coefficient drives $u$ towards $u=\pm1$ for $t<1.05$, wheres for $t>1.05$ the solution $u$ is driven towards $\pm1$ for $x$ in the  region where $\sin(x)<1.05/t$ and $0$ in the  region where $\sin(x)>1.05/t$; as $t\to\infty$ these two regions converge to $(0,\pi)$ and $(\pi,2\pi)$, respectively. The initial condition is
\[
u(0, x) = \varphi_{\text{G}}^L(x, \sqrt{10}, 1/2) - \varphi_{\text{G}}^L(x, \sqrt{10}, 4.4)\,,
\]
with $\varphi_{\text{G}}^L$ defined in Eq.~\eqref{eq:NumExp:Exp:ExpLUnit} and $L = 1/2$. 

Figure~\ref{fig:ACSpaceTimePlots}(a) shows a benchmark solution computed with a grid-based method. The benchmark solution is computed with a finite-difference discretization on 2048 equidistant grid points in the spatial domain. Time is discretized with semi-implicit Euler where the linear operators are treated implicitly and the nonlinear operators are treated explicitly. The time-step size is $10^{-5}$. 
The solution consists of relatively flat pieces separated by sharp walls, which evolve in a way that is consistent with the evolving sign structure of the coefficient $a(t,x)$ described above. Figure~\ref{fig:ACSpaceTimePlots}(b)  shows the Neural Galerkin approximation obtained with the three-hidden-layer network defined in \eqref{eq:NumExp:Tanh:DeepNetwork}, with $m = 2$ tanh units (16 degrees of freedom) and fixed uniform measure $d\nu_{\theta}(x) = dx$. The initial condition for the Neural Galerkin scheme is fitted analogously to fitting of the initial condition in the experiment with the KdV equation. Neural Galerkin uses a backward Euler discretization in time and then takes the gradient with automatic differentiation implemented in JAX (\url{https://github.com/google/jax}). The time-step size is $\delta t = 10^{-2}$. The nonlinear system is solved with ADAM and 10,000 iterations.
The corresponding optimization problem in each time step is approximately solved with SGD and $n = 1000$ samples to estimate $M$ and $F$. As can be seen, the Neural Galerkin solution agrees well with the benchmark solution.  In contrast, the approximation obtained with linear Galerkin with 16 equidistantly located Gaussian basis functions defined in \eqref{eq:NumExp:Exp:ExpLUnit} and shown in Figure~\ref{fig:ACSpaceTimePlots}(c) leads to a poorer approximation. This is further quantified in Figure~\ref{fig:ACSpaceTimePlots}(d) which shows that the Neural Galerkin approximation with a three-hidden-layer network achieves  an $\ell_2$ error  that is two orders of magnitude lower than that of the linear Galerkin approximation. 

Figure~\ref{fig:ACOutputsCoeff} visualizes the adaptation of the coefficients and features of the Neural Galerkin approximation. 
Consider the network \eqref{eq:NumExp:Tanh:DeepNetwork} with $\bfc(t) = (c_1(t), c_2(t))^T \in \mathbb{R}^{2}$, $\ell=3$, and set
\begin{multline}
\psi_i(t, \bfx) = \bfe_i^T \tanh(\bfW_{3}(t)\tanh(\bfW_{2 - 1}(t)(
\varphi^L_{\tanh}(\bfx, \bfW_1(t), \bfb_1(t)) ) + \bfb_{2}(t)) + \bfb_3(t))\,,
\label{eq:ACOutputComponent}
\end{multline}
where $\bfe_i \in \mathbb{R}^2$ with $i=1,2$ are the two canonical unit vectors. Thus the output components $\psi_1(t, \bfx)$ and $ \psi_2(t, \bfx)$ are weighted by the coefficients $c_1(t)$ and $c_2(t)$ to obtain the Neural Galerkin approximation; these two components are shown in Figure~\ref{fig:ACOutputsCoeff}(a) and (b). The corresponding coefficients are visualized in panel (c) and the benchmark solution is plotted again in panel (d). Even though the Neural Galerkin approximation uses only two $\psi_i(t, \bfx)$, because these functions depends nonlinearly on parameters that evolve in time, they adapt their shape to accurately approximate the solution.

\begin{figure}[!h]
\begin{center}
\begin{tabular}{ccc}
\hspace*{-0.25cm}\includegraphics[width=0.33\columnwidth]{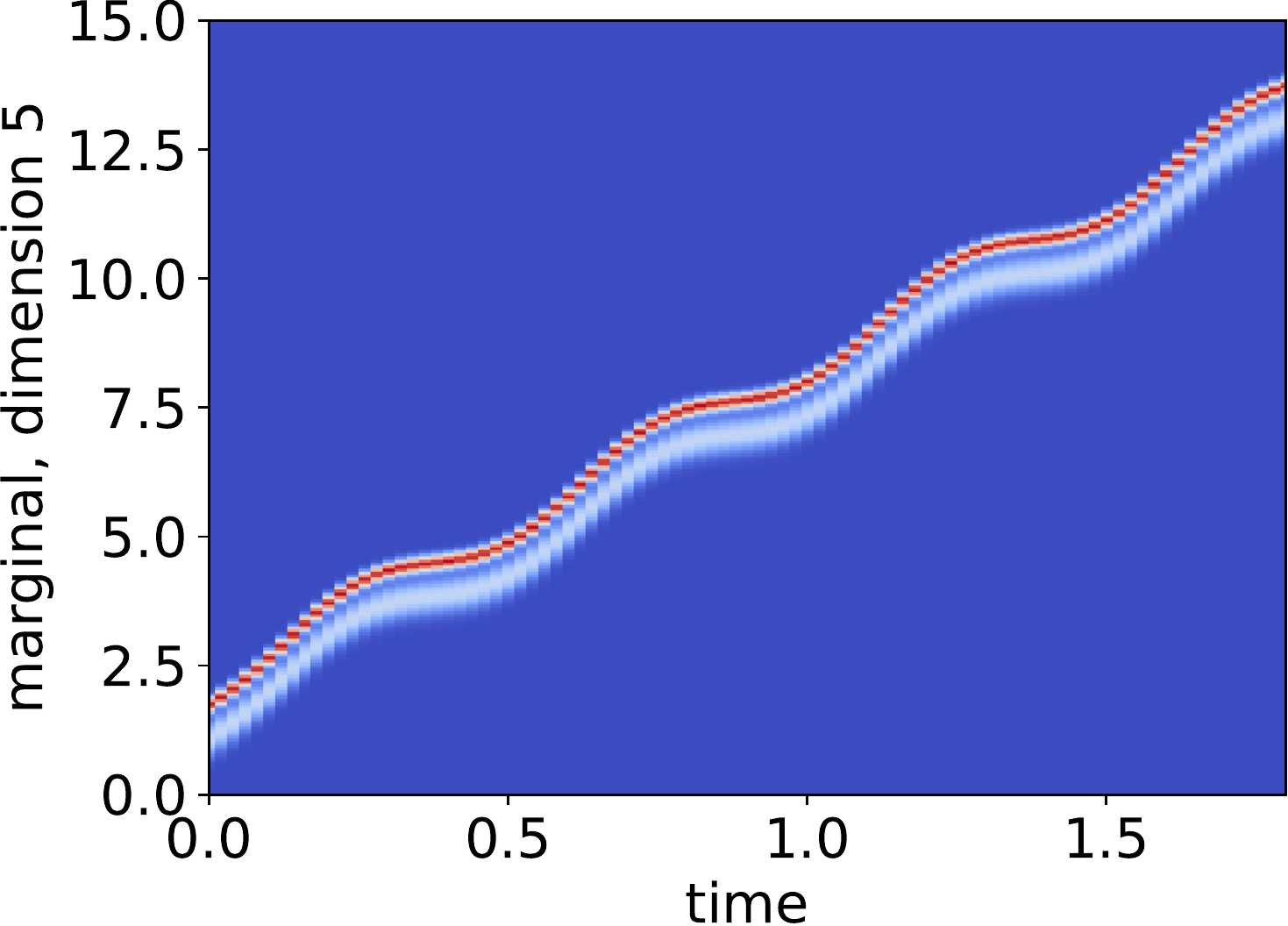} & \hspace*{-0.25cm}\includegraphics[width=0.33\columnwidth]{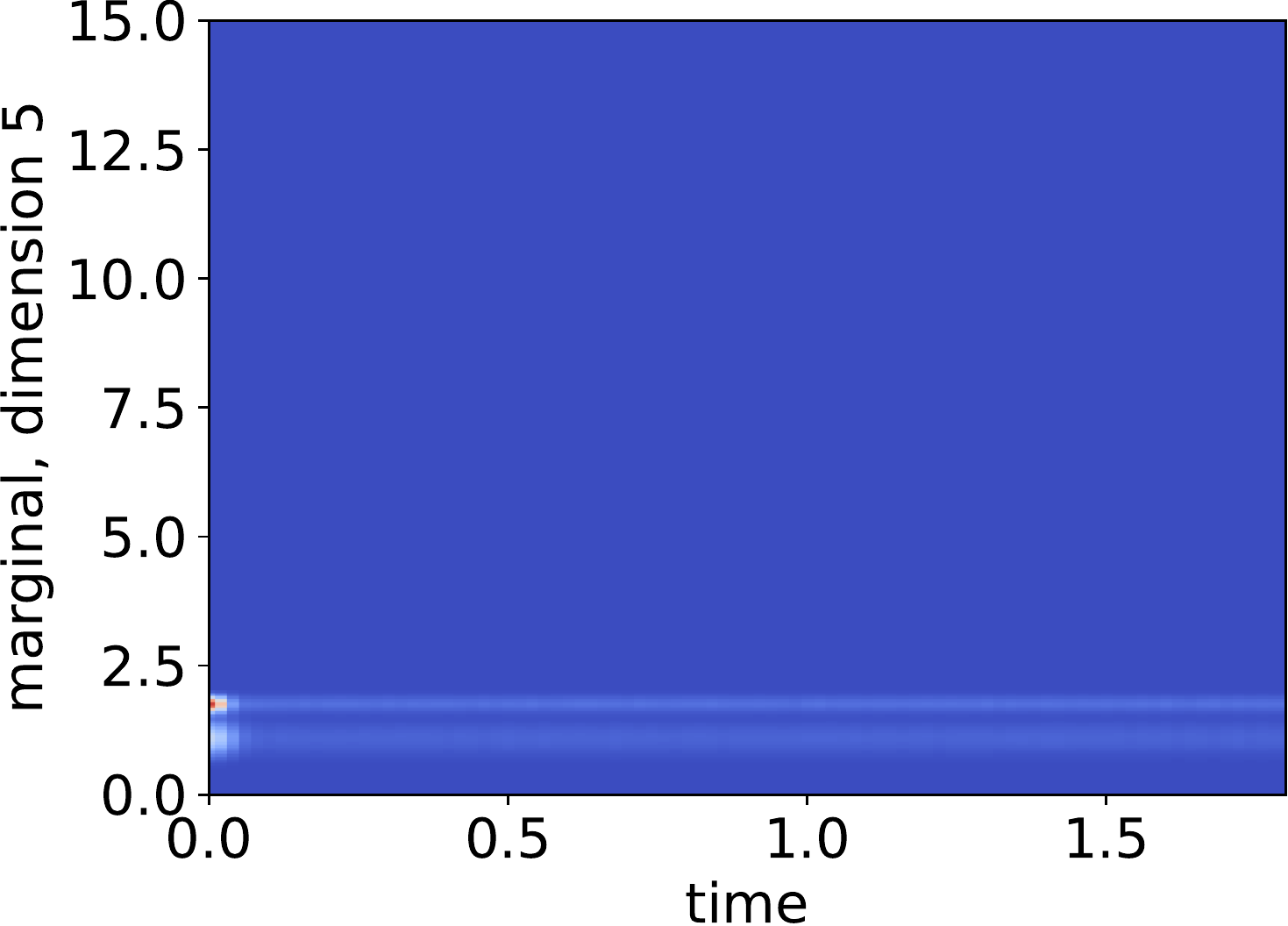} &
\hspace*{-0.25cm}\includegraphics[width=0.33\columnwidth]{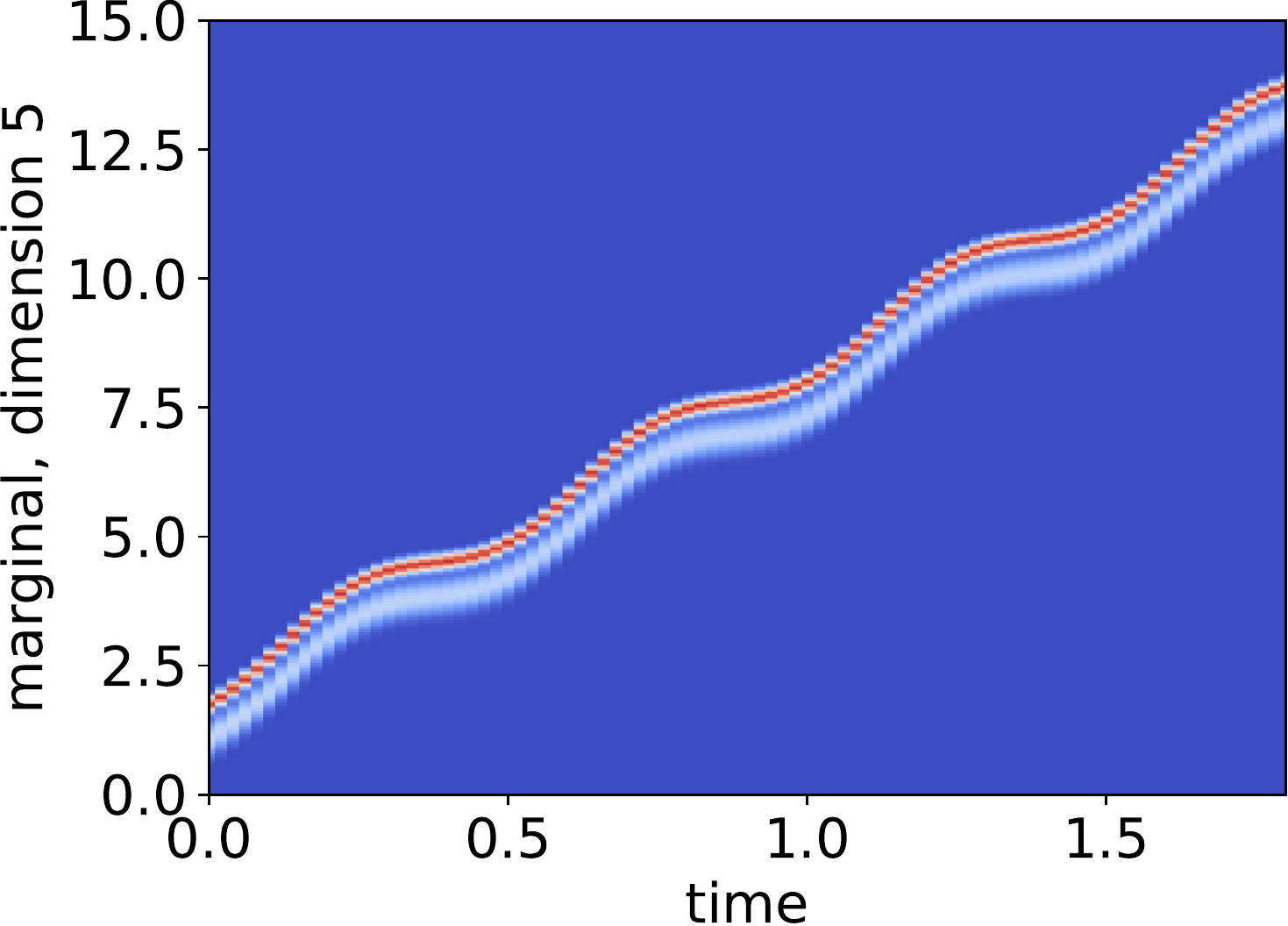}\\
\footnotesize (a) truth & \footnotesize (b) static & \footnotesize (c) Neural Galerkin\\
&\\
\multicolumn{3}{c}{\includegraphics[width=0.66\columnwidth]{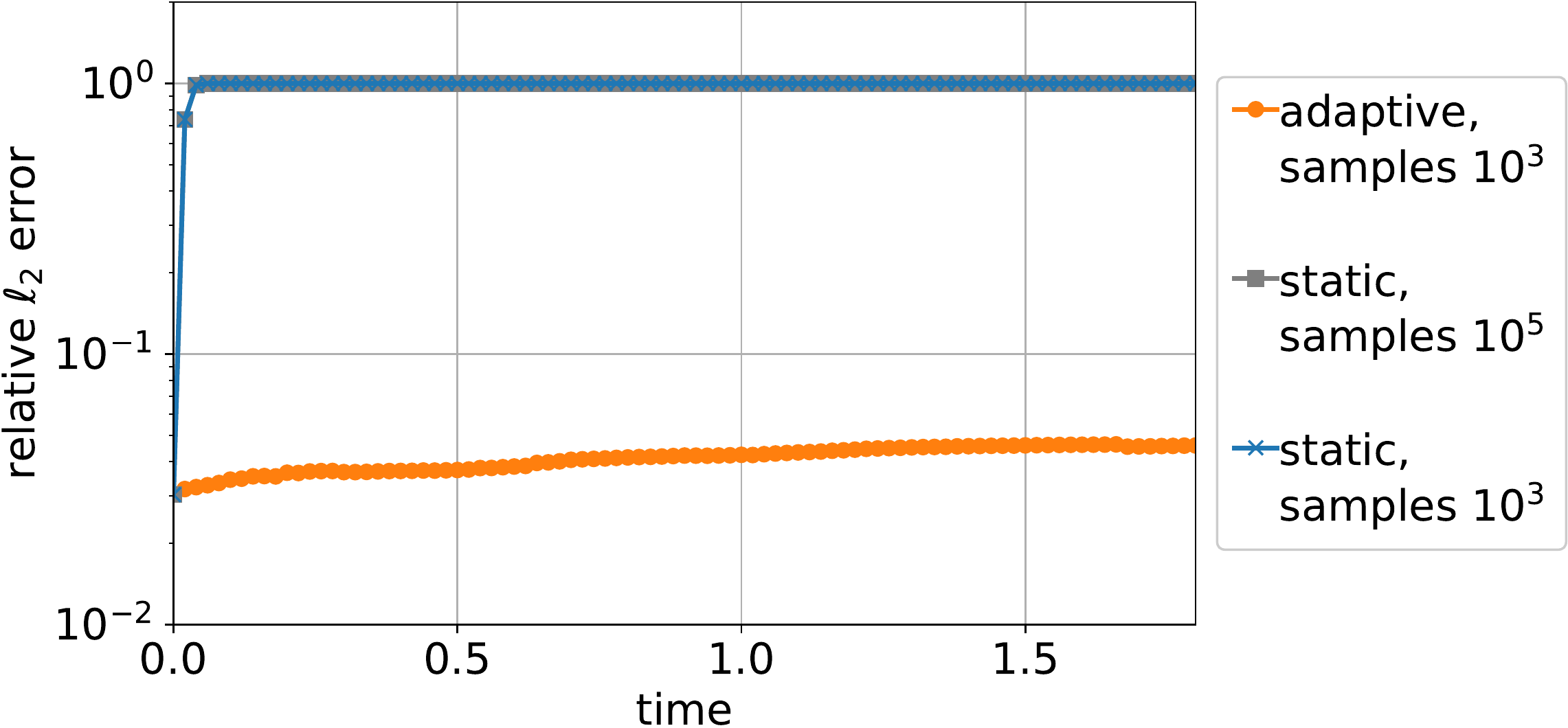}}\\
\multicolumn{3}{c}{\footnotesize (d) error over time}
\end{tabular}
\end{center}
\caption{High-dimensional advection with time-only varying advection speed: (a) marginal in dimension 5 of the analytic solution; (b) approximation obtained with static sampling; (c) result of the Neural Galerkin method that adaptively samples data over time to estimate the operators $M$ and $F$ depending on the dynamics of the problem at hand. This is in stark contrast to classical time-space collocation approaches that uniformly sample over space and time when optimizing for a solution. In this example, the adaptive sampling is key for Neural Galerkin to accurately predict the local-in-space dynamics of the solution. In contrast, an approximate solution obtained with uniform sampling that is static over time fails to lead to meaningful predictions, as shown by plot (b) and the error shown in plot (d).}
\label{fig:HighTransportOnlyTime}
\end{figure}

\subsection{Advection in unbounded, high-dimensional domains}\label{sec:NumExp:Adv}

Our next example is an advection equation in high dimension,
\begin{equation}
\partial_t u + \bfa(t, \bfx) \cdot \nabla_{\bfx} u = 0, \qquad u(0,\bfx)= u_0(\bfx),
\label{eq:Transport:Eq}
\end{equation}
where $\bfx\in \Xcal\equiv \R^d$, $\bfa:[0,\infty)\times \R^d \to \R$ is some bounded velocity field that is assumed to be differentiable in both its arguments, and $u_0:\R^d\to\R$ is some initial condition satisfying $\lim_{|\bfx|\to\infty} u_0(\bfx) = 0$. This equation can in principle be solved by the method of characteristics, i.e. by considering
\begin{equation}
    \label{eq:charact}
    \dot \bfX(t,\bfx) = \bfa(t,\bfX(t,\bfx)), \qquad \bfX(0,\bfx) = \bfx,
\end{equation}
and using
\begin{equation}
    \label{eq:sol:adv}
    u(t,\bfx) = u_0(\bfX(-t,\bfx)).
\end{equation}
In practice, however, Eq.~\eqref{eq:sol:adv} is hard to use to obtain global information about the solution, except in situations where Eq.~\eqref{eq:charact} can be solved explicitly (an instance of which we will also consider below as benchmark). Next we show how to numerically solve~\eqref{eq:Transport:Eq} directly using a Neural Galerkin scheme.

\begin{figure}[!h]
\begin{tabular}{cc}
\includegraphics[width=0.44\columnwidth]{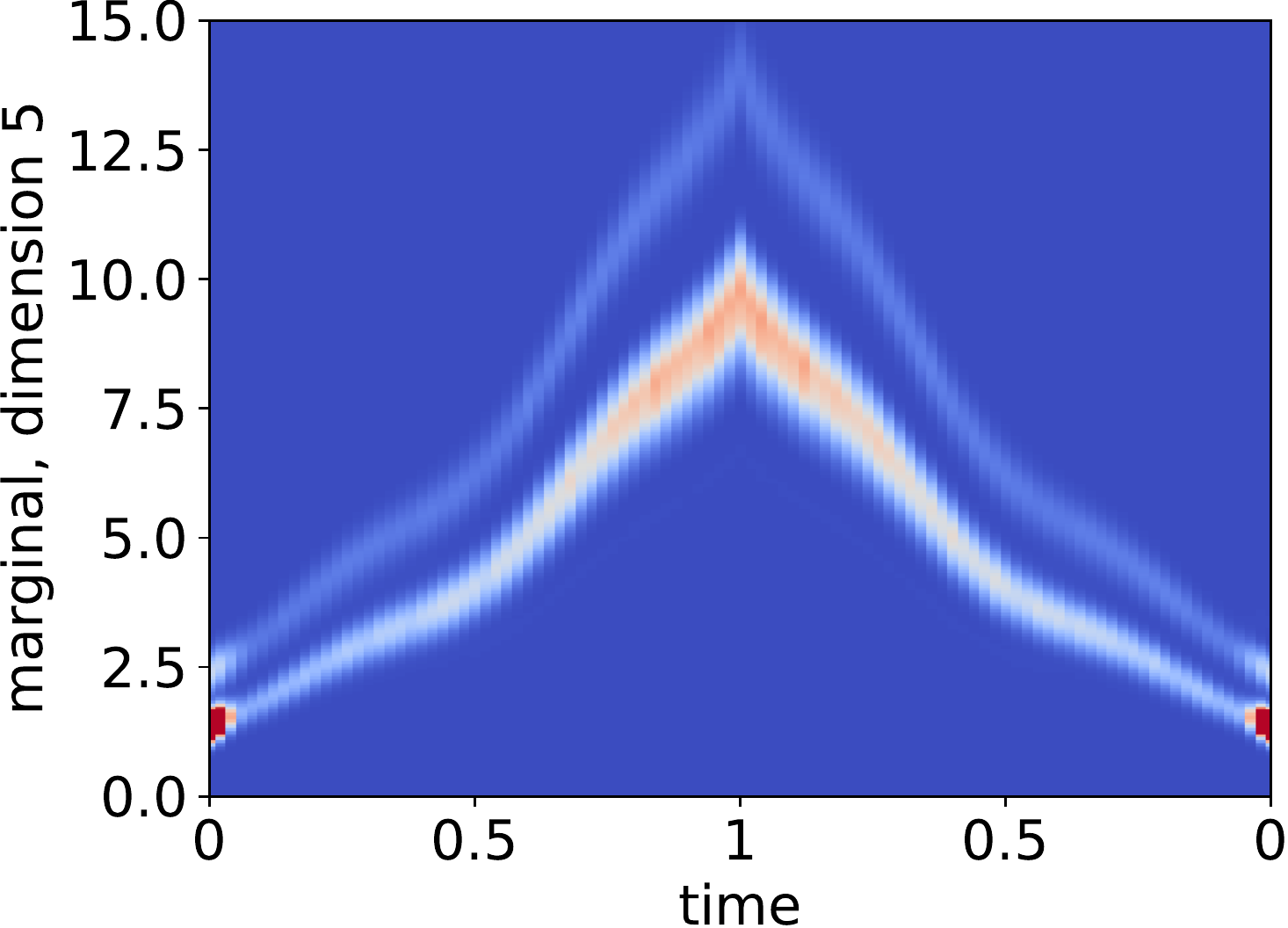} \hspace*{0.5cm}&
\includegraphics[width=0.42\columnwidth]{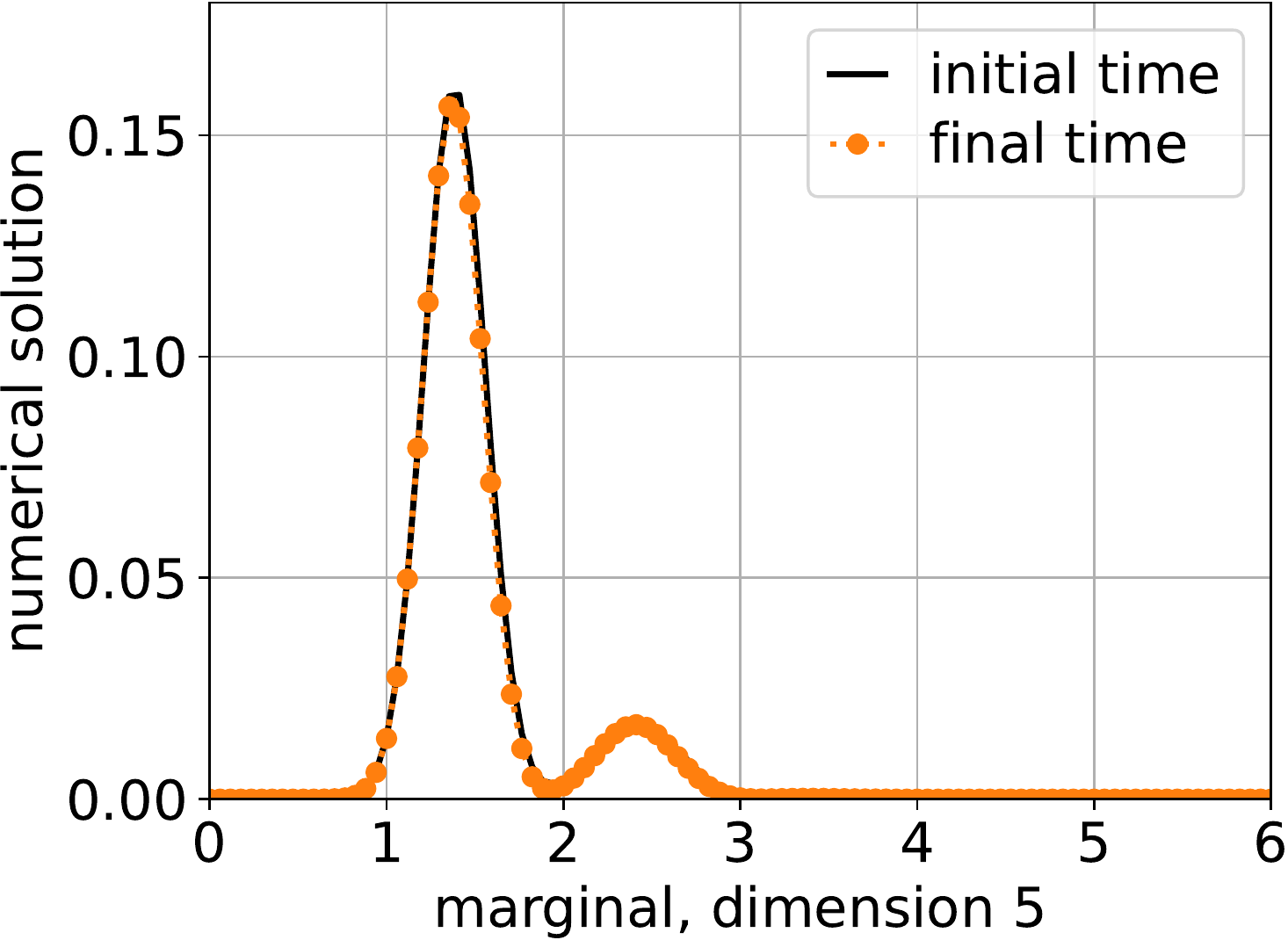}\\ 
\footnotesize (a) Neural Galerkin with adaptive measure & \footnotesize (b) initial vs.~final solution\\
&\\
\multicolumn{2}{c}{\includegraphics[width=0.66\columnwidth]{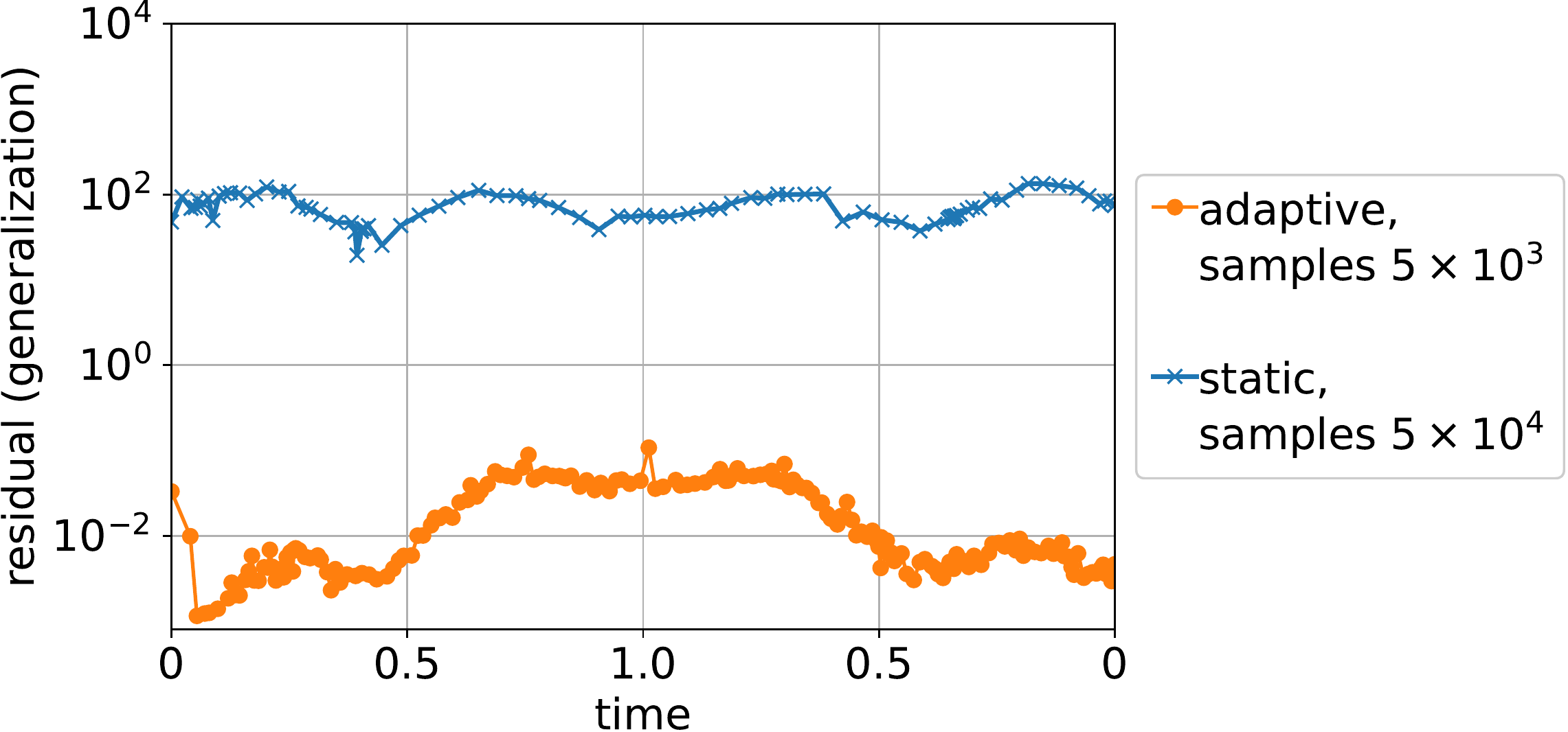}} \\
\multicolumn{2}{c}{\footnotesize (c)  residual}
\end{tabular}
\caption{High-dimensional advection with time-space varying advection speed: The equation is integrated forward in time until $t = 1$ and then integrated backwards until time $t = 0$, where the initial condition is reached again. Plot (a) shows the prediction of the Neural Galerkin method. Plot (b) shows that the Neural Galerkin solution at final time closely matches the initial condition, which means that little error is accumulated by integrating forward and then backward in time.  Adaptive sampling is key in this example because the solution is local in the high-dimensional spatial domain. Neural Galerkin with adaptive sampling tracks the local dynamics well, whereas  approximations based on static sampling fail after only a few time steps as shown by a large residual in (c).}
\label{fig:HighTransportTimeSpace}
\end{figure}

\subsubsection{Advection with a time-dependent coefficient.}
In the first experiment with the advection equation, we make the transport coefficient depend only on time using
\begin{equation}
\bfa_t(t) = \bfa_s \odot \left(\sin(\bfa_v \pi t) + 5/4\right)
\label{eq:Transport:TimeCoeff}
\end{equation}
with $\bfa_s = [1, 2, \dots, d]^T$, $\bfa_v = 2 + \frac{2}{d}[0, 1, \dots, d - 1]^T$ and the element-wise vector multiplication $\odot$. We set $d=5$  and take as  initial condition $u_0$  a mixture of two non-isotropic Gaussian packets with means
\[
\bfmu_1 = \frac{11}{10}\begin{bmatrix}
1\\
1\\
\vdots\\
1
\end{bmatrix}\,,\qquad \bfmu_2 = \frac{3}{4} \begin{bmatrix}
1.5 - (-1)^11/(d + 1)\\
1.5 - (-1)^22/(d + 1)\\
\vdots\\
1.5 - (-1)^dd/(d + 1)
\end{bmatrix}
\]
and covariance matrices
\begin{equation}
\bfSigma_1 = \frac{1}{200}\begin{bmatrix}
2 & & &\\
& 4 & &\\
& & \ddots &\\
& & & 2d
\end{bmatrix}\,,\qquad \bfSigma_2 = \frac{1}{200}\begin{bmatrix}
d & & &\\
& d-1 & &\\
& & \ddots &\\
& & & 1
\end{bmatrix}\,.
\end{equation}
In this case, the solution~\eqref{eq:sol:adv} can be derived explicitly,
\begin{equation}
    \label{eq:explicit:sol:adv}
    u(t,\bfx) = u_0\left( \bfx-{\textstyle\int_0^t}\bfa_t(s)ds\right)\,,
\end{equation}
where $\bfa_t$ is the time-varying coefficient defined in the main text. We will use the exact solution as benchmark. Marginals of the solutions are computed as follows: Let $\Xcal = \Xcal_1 \times \dots \times \Xcal_d$ be the domain of interest. We then plot for each dimension $i = 1, \dots, d$ the function
\[
(t, x_i) \mapsto \int_{\Omega_1} \cdots \int_{\Omega_{i - 1}} \int_{\Omega_{i + 1}} \cdots \int_{\Omega_d} u(t, x_1, \dots, x_{i - 1}, x_{i + 1}, \dots, x_d)\,,
\]
where $u$ is the function of which marginals are to be plotted and the integrals are numerically approximated via Monte Carlo and $8192$ samples. The samples are drawn from the analytic solution, which is available in this experiment. The marginal in dimension five is shown in Figure~\ref{fig:HighTransportOnlyTime}(b); see also Figure~\ref{supp:fig:NumExp:TransportTime} in the Appendix.

For Neural Galerkin, see Figure~\ref{fig:HighTransportOnlyTime}(b), we use a shallow network as in Eq.~\eqref{eq:NumExp:Exp:ShallowNetwork} with the exponential unit~\eqref{eq:NumExp:Exp:ExpUnit}---notice that since these units are isotropic, several of them are required to accurately approximate both the initial condition and the solution of Eq.~\eqref{eq:Transport:Eq} at time $t>0$. We take $m = 50$ nodes and use the adaptive RK45 method as time-integrator. We also use a measure $\nu_\theta$ adapted to the solution to sample $n=1000$ data points $\{\bfx_i\}_{i=1}^n$ at each integration step to estimate $M$ and $F$ via Eq.~\eqref{eq:F:M:emp:2}. Specifically, we take $\nu_\theta$ to be the Gaussian mixture with 50 nodes obtained from the current neural approximation of the solution by equating the weights $c_i$, dividing $w_i$ by a factor $\kappa > 0$, and keeping the same $\bfb_i$, i.e. if the current Neural Galerkin approximation of the solution is
\begin{equation}
U(\theta, \bfx) = \sum\nolimits_{i = 1}^{50} c_i(t) \exp\left(-w_i^2(t)|\bfx - \bfb_i(t)|^2\right)\,,
\label{eq:NumExp:Exp:ExpUnit:current}
\end{equation}
we sample the data points from
\begin{equation}
d\nu_{\theta(t)} (\bfx) = C^{-1}\sum\nolimits_{i = 1}^{50}  \exp\left(-\tfrac{1}{\kappa^2} w_i^2(t)|\bfx - \bfb_i(t)|^2\right)d\bfx\,,
\label{eq:nu:current}
\end{equation}
where $C$ is a normalization constant and $\kappa = 1$ in this experiment; see Section~\ref{sec:Active:Arch}. Using this adapted measure is key for accuracy, and leads to the small $\ell_2$ error plotted in Figure~\ref{fig:HighTransportOnlyTime}(d). If instead we estimate $M$ and $F$ by drawing the data points uniformly in $[0, 15]^d$ (to cover the domain over which the solution propagates for $t\in[0,2]$), the relative error is 100\% after only a few time steps even if we use as many as $n=10^5$ data points, see  Figure~\ref{fig:HighTransportOnlyTime}(b) and Figure~\ref{fig:HighTransportOnlyTime}(d). This emphasizes the importance of using adaptivity in both the function approximation via neural networks, and the data acquisition to estimate $M$ and $F$. This double-adaptivity is a key distinguishing feature of the proposed Neural Galerkin schemes compared to existing DNN approaches.

\subsubsection{Advection with time- and space-dependent coefficients.}
Next we consider Eq.~\eqref{eq:Transport:Eq} with the time-space varying advection speed,
\begin{equation}
\bfa_{st}(t, \bfx) = \bfa_s \odot \left(\sin(\bfa_v \pi t) + 3\right) \odot (\bfx + 1)/10\,,
\end{equation}
where $\bfa_s, \bfa_v$ are the vectors defined in Eq.~\eqref{eq:Transport:TimeCoeff}.
For the initial condition we take a sum of two Gaussian packets, $u_0(x) = p_1(x)/10 + p_2(x)/10$. The two Gaussian probability density functions $p_1, p_2$ are scaled by a factor $1/10$ to avoid scaling issues. The means are
\[
\bfmu_1 = 2 - \frac{1}{12}\begin{bmatrix}-1 & 2 & -3 & 4 & -5 \end{bmatrix}^T\,,\qquad \bfmu_2 = 1.8 - \frac{1}{12}\begin{bmatrix} 1 & -2 & 3 & -4 & 5 \end{bmatrix}^T
\]
and the covariances are
\[
\bfSigma_1 = \begin{bmatrix}
\frac{3}{50} & &\\
& \ddots & \\
& & \frac{3}{50}
\end{bmatrix}\,,\qquad \bfSigma_2 = \begin{bmatrix}
\frac{3}{100} & &\\
& \ddots & \\
& & \frac{3}{100}
\end{bmatrix}\,.
\]
The marginals for all five dimensions of the Neural Galerkin solutions are shown in the appendix in Figure~\ref{supp:fig:LinAdvSpatiallyVarying}. 

Since the analytic solution for this problem is not available in closed form, to assess the performance of our Neural Galerkin approach, we  numerically integrate the problem forward in time until time $T = 1$, then invert the advection coefficient in time, and  integrate backwards until time $0$:  we can then estimate the error by comparing the final numerical solution with the initial condition. The numerical setup is the same as in the previous experiment with the advection speed~\eqref{eq:Transport:TimeCoeff}, except that we draw $n=5000$ samples  from the adapted measure in~\eqref{eq:nu:current} at  each integration step to estimate $M$ and $F$. Figure~\ref{fig:HighTransportTimeSpace}(a) shows the time-space marginal of the Neural Galerkin solution: as  can be seen the spatially varying coefficient leads to a significantly changing solution over time.  Figure~\ref{fig:HighTransportTimeSpace}(b) compares the initial condition with the final solution obtained by integrating forward and then backward in time, which are in good agreement. In Figure~\ref{fig:HighTransportTimeSpace}(c) we show an estimate of the PDE residual defined in Eq.~\eqref{eq:NG:ContTime:MinResObj} using the adapted measure in~\eqref{eq:nu:current} and estimated using $10^5$ samples drawn independently from this measure. The mean residual is computed as follows: We draw $n = 10^5$ samples from the equally-weighted mixture of Gaussians given by the $m = 50$ nodes of the network at time $t$. The standard deviation of the Gaussians is double the standard deviation of the nodes of the network. With these samples, we compute a Monte Carlo estimate of the objective $J_t$ at time $t$. The residuals indicate that the proposed Neural Galerkin approach with adaptive sampling approximates well the local dynamics in this high-dimensional transport problem,  whereas static sampling with a uniform distribution fails again to provide meaningful predictions.

\subsection{Interacting particle systems}\label{sec:NumExp:Particle} Let $\Xcal = \mathbb{R}^d$ and consider the evolution of $d$ interacting particles with positions $X_1(t), \dots, X_d(t) \in \mathbb{R}$ governed by 
\begin{equation}
\mathrm d X_i = g(t,X_i)\mathrm dt +\sum\nolimits_{j = 1}^d K(X_i,X_j)\mathrm dt + \sqrt{2D}\,\mathrm dW_i\,,
\label{eq:ParticlesSDE:1}
\end{equation}
with $i = 1, \dots, d$ and where $g:[0,\infty)\times \R\to\R$ is a time-dependent one-body force, $K:\R\times\R\to\R$ a pairwise interaction term, $D>0$ the diffusion coefficient, and $W_i$ are independent Wiener processes. The evolution of the joint probability density of these particles, $u(t,\bfx)$ with $\bfx=[x_1,\ldots,x_d]^T$, is governed by the Fokker-Planck equation
\begin{equation}
    \label{eq:FPE:1}
    \partial_t u =\sum\nolimits_{i=1}^d \left(-\partial_{x_i} \left( u h_i(t,x_1,\ldots, x_d)\right) + D \partial^2_{x_i}  u\right),
\end{equation}
where $h_i(t,x_1,\ldots, x_d) = g(t,x_i) +\sum_{j=1}^d K(x_i,x_j)$.
In high dimension, one typically resorts to Monte Carlo methods based on integrating Eq.~\eqref{eq:ParticlesSDE:1} to obtain samples from the density $u$. In contrast, we solve the Fokker-Planck equation~\eqref{eq:FPE:1} directly with our Neural Galerkin approach to derive an approximation of $u(t,\bfx)$ via $U(\theta(t),\bfx)$ everywhere in the domain $\Xcal$. Having the density at hand allows us to efficiently compute quantities of interest that are not expressible as expectation over $u$,  such as the entropy, and therefore not directly accessible to Monte Carlo methods.
Note that Neural Galerkin schemes also involve sampling to estimate $M$ and $F$, but it is done adaptively based on the current approximation to keep low the number of samples required to obtain accurate estimates. 

\begin{SCfigure}
\centering
\includegraphics[width=0.5\columnwidth]{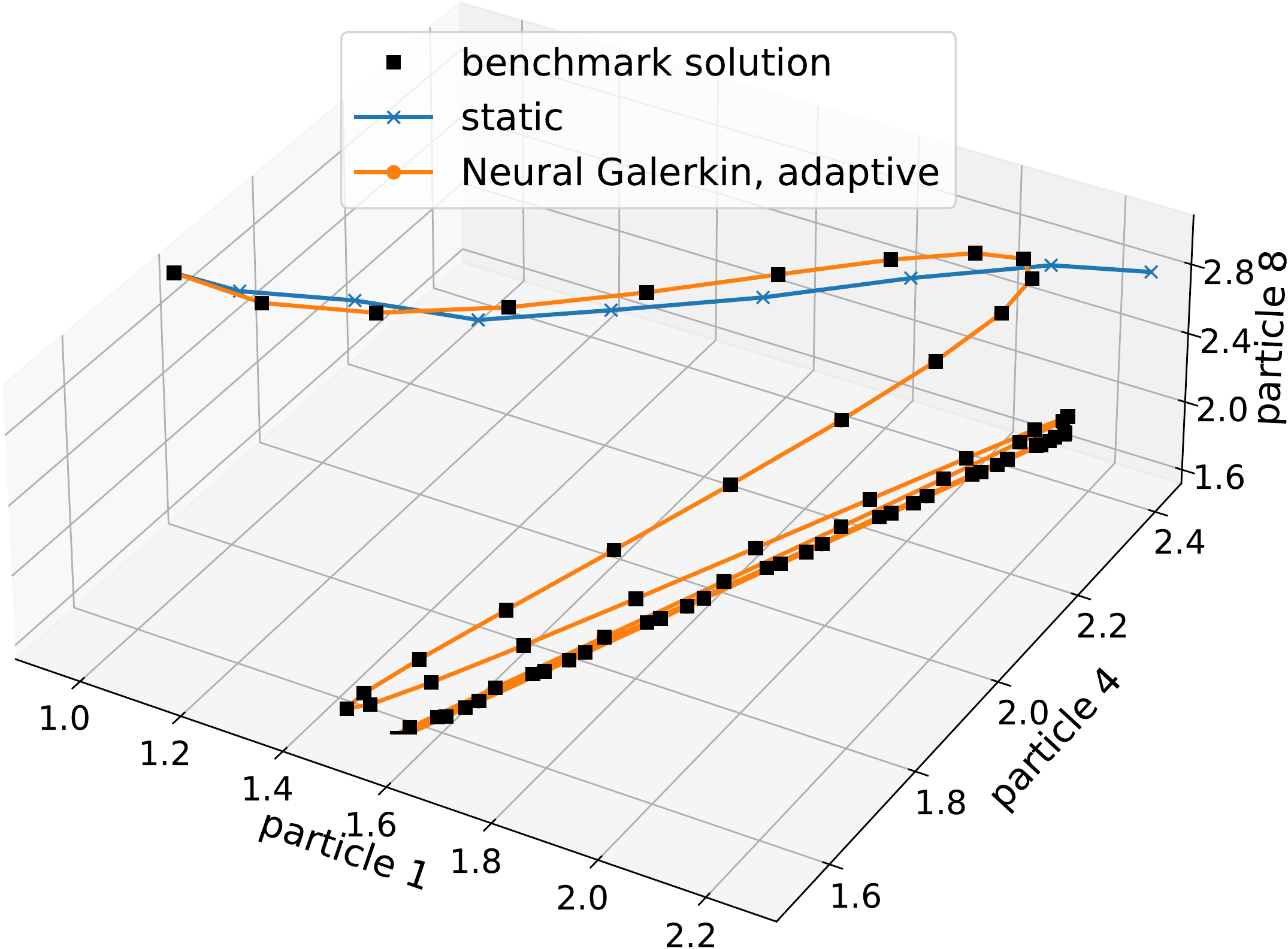}
\caption{Particles in harmonic trap: As the particles become attracted by the trap, the particle density becomes increasingly concentrated in the high-dimensional domain. The  Neural Galerkin scheme adaptively samples from the neural approximation of the density over time and thereby accurately tracks  the particles mean position (orange) in contrast to methods with static sampling done uniformly in the spatial domain (blue).}
\label{fig:ParticlesHarmonic3D}
\end{SCfigure}

\subsubsection{Particles in harmonic trap.}
In the following experiment, we set
\begin{equation}
    \label{eq:F:K}
    g(t,x) = a(t) - x, \qquad K(x, y) = \frac{\alpha}{d}(y-x)
\end{equation}
with $a:[0,\infty) \to \R$ and $\alpha>0$. This corresponds to following the evolution of $d$ particles put in a harmonic trap centered around the moving position $[a(t),\ldots,a(t)]^T \in \Xcal$ and also attracting each other harmonically. This choice has the advantage that we can derive closed ODEs for the mean and covariance of the particle positions that can be solved numerically and used as benchmark, see~\ref{sec:Appx:FokkerPlanck}. If the initial condition $u(0,\bfx)$ is a Gaussian probability density, the density $u(t, \bfx)$ remains Gaussian for all times $t$, with the mean and covariance given by the ODEs.
Accordingly, we take $u(0,\bfx)$ to be an isotropic Gaussian density with mean $9/10 + 21/(10(d - 1))[0, 1, \dots, d - 1]^T$ and variance $\sigma^2 = 0.1$. We take $d = 8$ particles and set $a(t) = 5/4(\sin(\pi t) + 3/2)$, $\alpha = 1/4$ and $D = 10^{-2}$.

\begin{figure}
\begin{tabular}{cc}
\includegraphics[width=0.45\columnwidth]{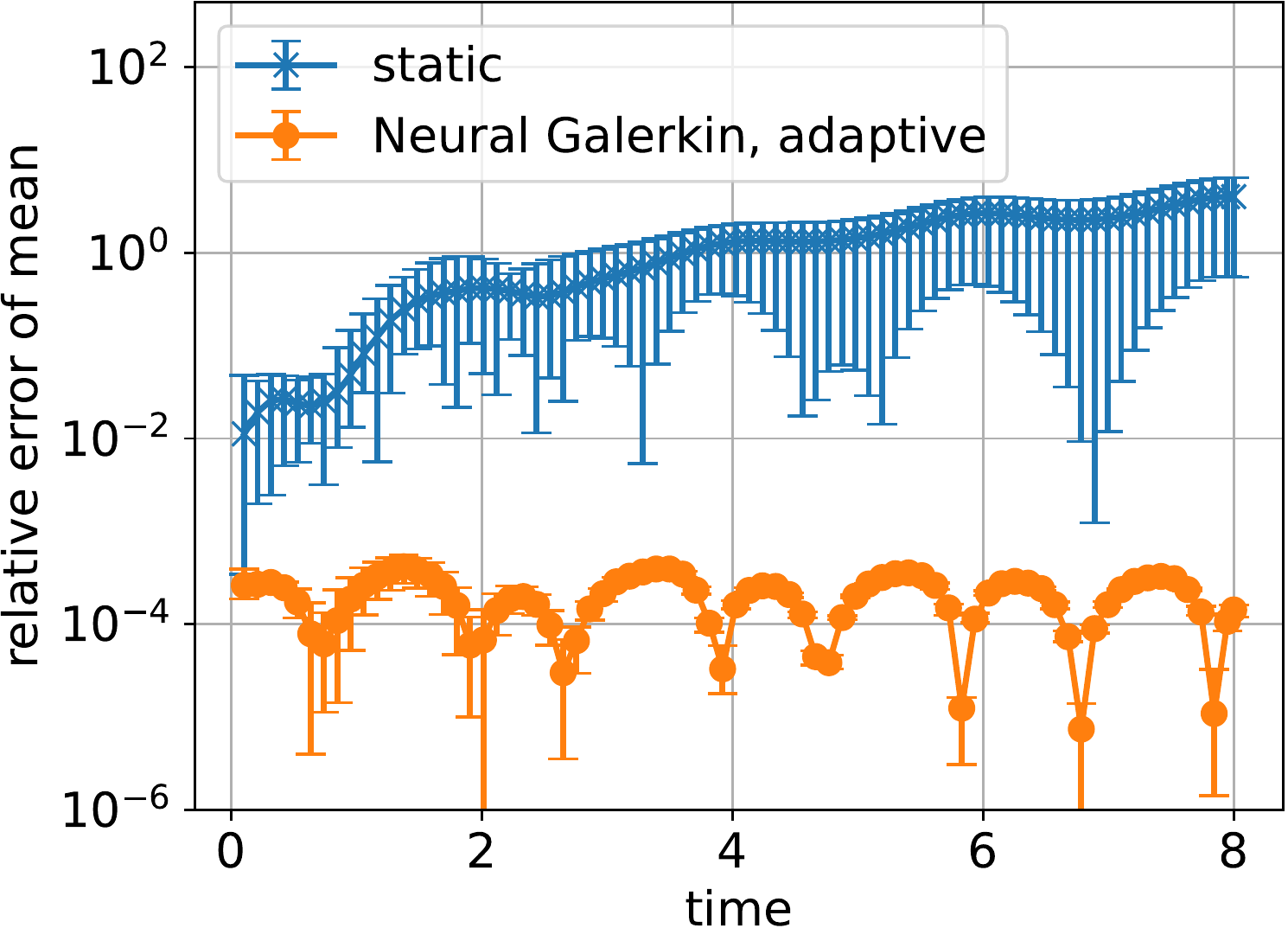} & \includegraphics[width=0.45\columnwidth]{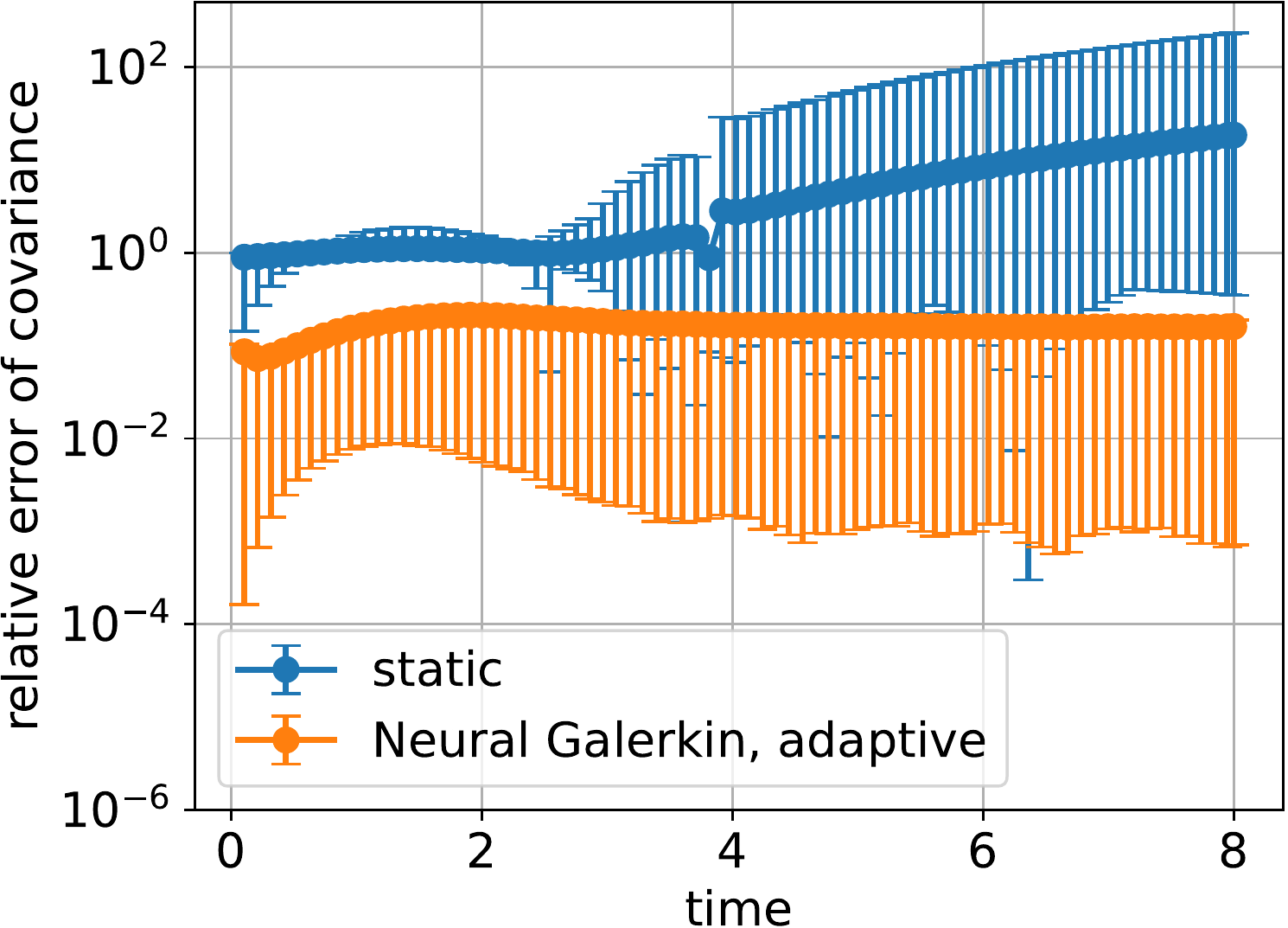}\\
\footnotesize (a) mean prediction & \footnotesize (b) covariance prediction
\end{tabular}
\caption{Particles in harmonic trap: Our Neural Galerkin schemes with adaptive sampling achieve orders of magnitude lower errors in predicting the mean (orange, panel a) and the covariance (orange, panel b) of the particle density than methods that rely on static sampling uniform over time and space (blue, panels a and b). Additionally, the Neural Galerkin scheme provides an approximation of the particle density that allows efficiently computing quantities of interest such as the entropy (Figure~\ref{fig:ParticleCubicDim8Entropy}a), in contrast to Monte Carlo that only draws samples of the particle distribution.}
\label{fig:ParticleDim8}
\end{figure}

We apply our Neural Galerkin approach to solve the Fokker-Planck equation~\eqref{eq:FPE:1} using the neural network defined in \eqref{eq:NumExp:Exp:ShallowNetwork} with $m = 30$ nodes and non-negative weights to guarantee that the Neural Galerkin approximation is non-negative, i.e., 
\begin{equation}
U(\theta, \bfx) = \sum\nolimits_{i = 1}^m c_i^2\varphi(\bfx, w_i, \bfb_i)\,.
\label{eq:SI:ShallowSquaredWeights}
\end{equation}
The normalized Neural Galerkin approximation is given by
\begin{equation}
\bar{U}(\theta, \bfx) = \frac{1}{\sum_i c_i^2} \sum\nolimits_{i = 1}^m c_i^2\varphi(\bfx, w_i, \bfb_i)\,.
\end{equation}
The scaling factor $\kappa$ in the adaptive sampling is set to $\kappa = 2$ in this experiment.

\begin{figure}
\begin{tabular}{cc}
\includegraphics[width=0.45\columnwidth]{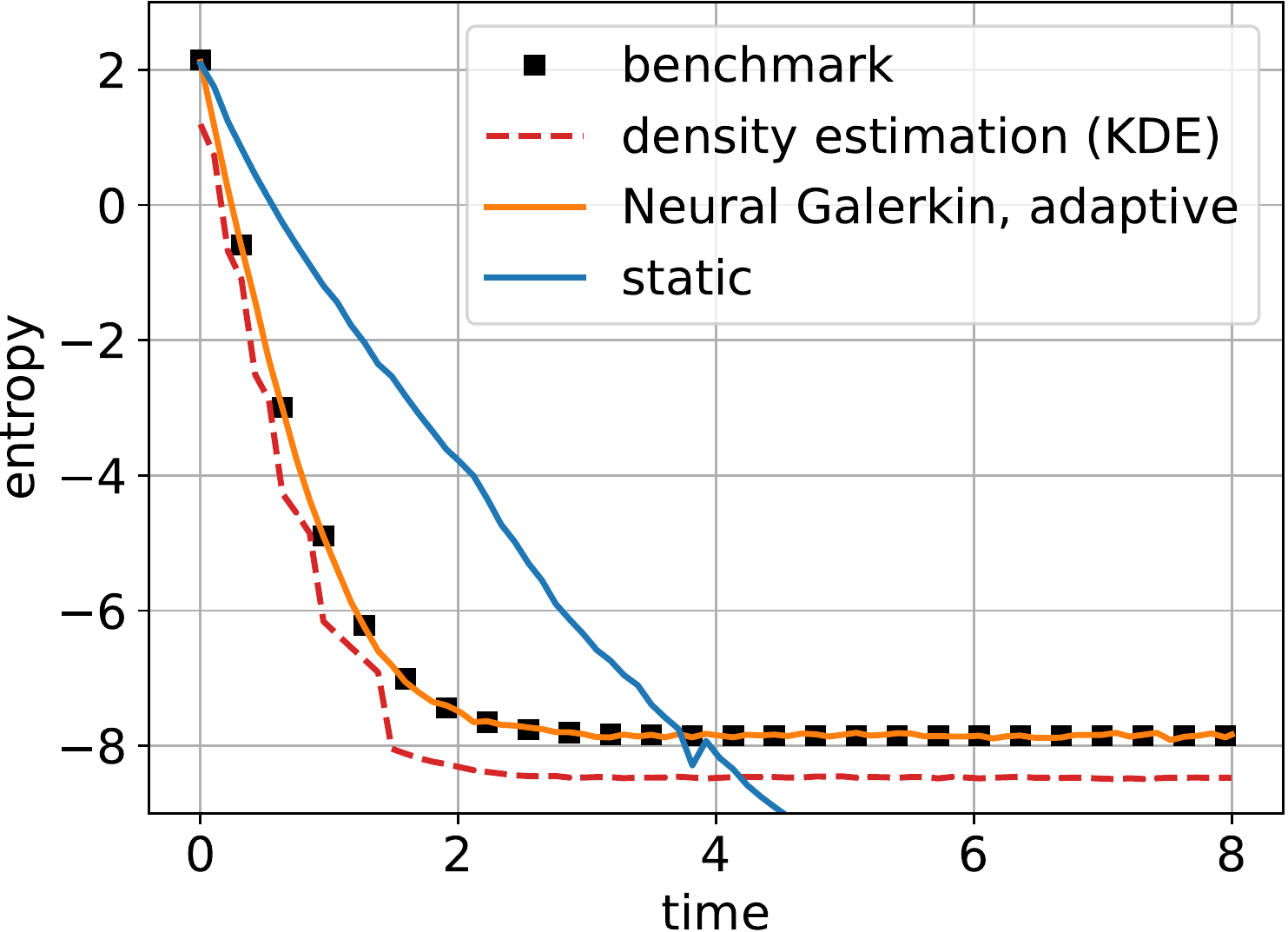} &
\includegraphics[width=0.45\columnwidth]{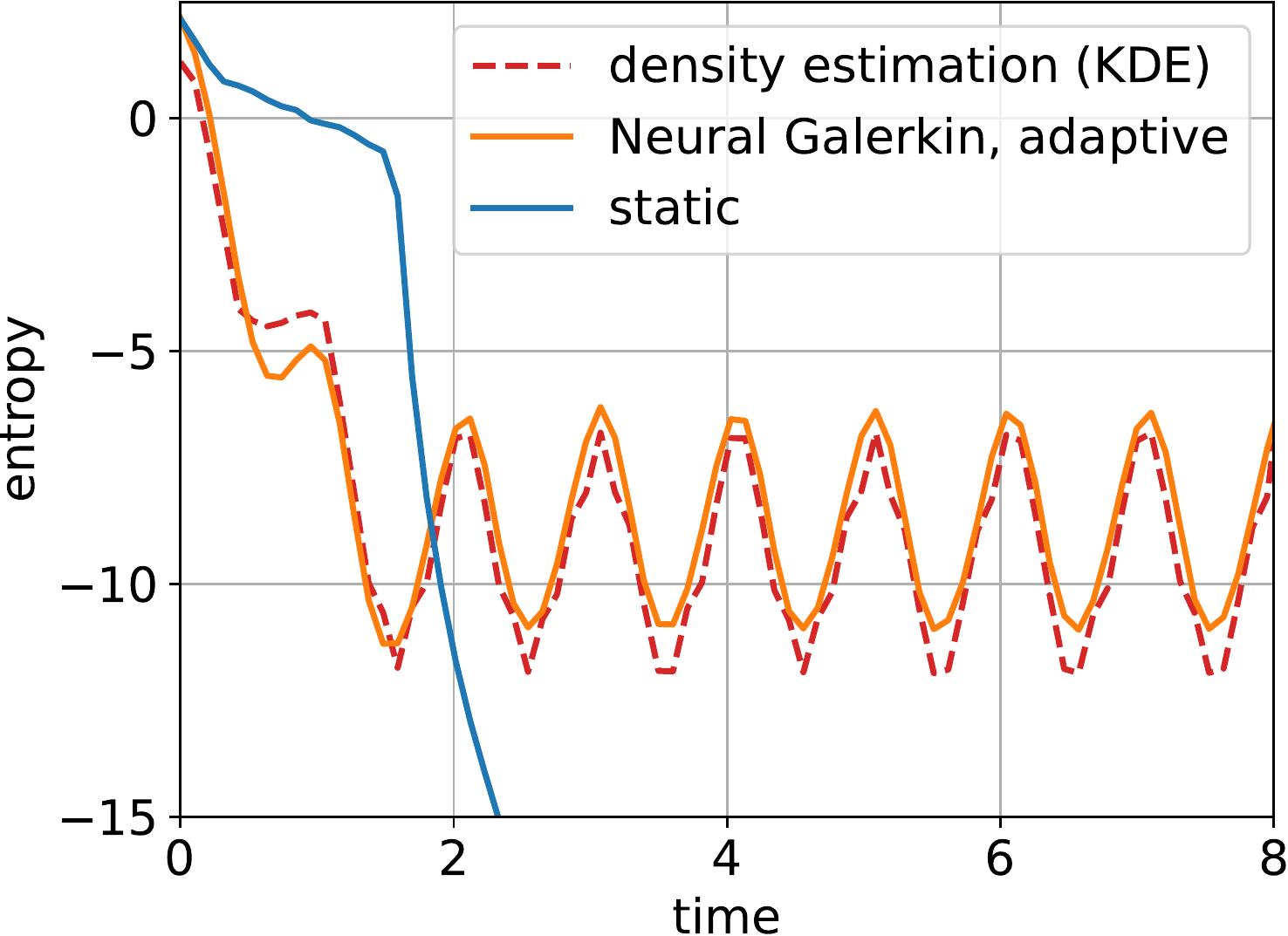}\\
\footnotesize (a) entropy, harmonic trap & \footnotesize (b) entropy, aharmonic trap\\
\end{tabular}
\caption{Particles in trap: The proposed Neural Galerkin schemes provide an approximation of the particle density that allows efficiently computing quantities of interest such as the entropy. In contrast, Monte Carlo estimation only draws samples of the particle distribution and requires subsequent density estimation with e.g., kernel density estimation (KDE) to compute quantities such as the entropy.}
\label{fig:ParticleCubicDim8Entropy}
\end{figure}

We use implicit Euler with time-step size $\delta t = 10^{-3}$. The number of SGD iterations to solve Eq.~\eqref{eq:ODE:DiscImp} in $\theta^{k+1}$ is $5 \times 10^3$.  To estimate $M$ and $F$, we draw $n = 1000$ samples at every time step, using the adaptive sampling procedure used in the previous experiments with the advection equation. Figure~\ref{fig:ParticlesHarmonic3D} shows the mean position of particles 1, 4, and 8 in the spatial domain over time. The approximation obtained with Neural Galerkin with adaptive sampling closely follows the benchmark solution. In contrast, estimation of $M$ and $F$ via uniform sampling in $[0, 5]^d$ leads again to a poor approximation. The relative error of the predicted mean with adaptive Neural Galerkin is roughly $10^{-4}$ as shown in Figure~\ref{fig:ParticleDim8}(a) and the entries of the covariance matrix are approximated with relative error of about $10^{-1}$ to $10^{-2}$ as shown in Figure~\ref{fig:ParticleDim8}(b). Because Neural Galerkin schemes provide an approximation of the density of the particles, rather than just samples of their positions as Monte Carlo methods, we can compute the entropy of the particle distribution over time; see Figure~\ref{fig:ParticleCubicDim8Entropy}(a). The entropy of the distribution described by the Neural Galerkin approximation is computed by drawing $n = 5000$ samples $\bfx_1, \dots, \bfx_n$ from the mixture given by \eqref{eq:SI:ShallowSquaredWeights} and then estimating the entropy via Monte Carlo
\[
E_n = -\frac{1}{n}\sum_{i = 1}^n \log \bar{U}(\theta(t), \bfx_i)\,.
\]
The adaptive Neural Galerkin approximation of the entropy is in close agreement with the benchmark, whereas static sampling and Monte Carlo combined with density estimation from $5,000$ samples leads to poor approximations.

\subsubsection{Particles in aharmonic trap.}
Next we consider an aharmonic trap with one-body force $g(t, X_i) = (a(t) - x)^3$ and the same pairwise interaction term $K$ as defined in Eq.~\eqref{eq:F:K} except that $\alpha$ is set to $\alpha = -0.5$, so that the particles now repel each other withn the attracting trap. Because of the nonlinear force, the density $u$ is not necessarily Gaussian anymore. A Monte Carlo estimate of the particle positions from $100,000$ samples serves as benchmark in the following.

\begin{figure}
\begin{tabular}{cc}
\includegraphics[width=0.45\columnwidth]{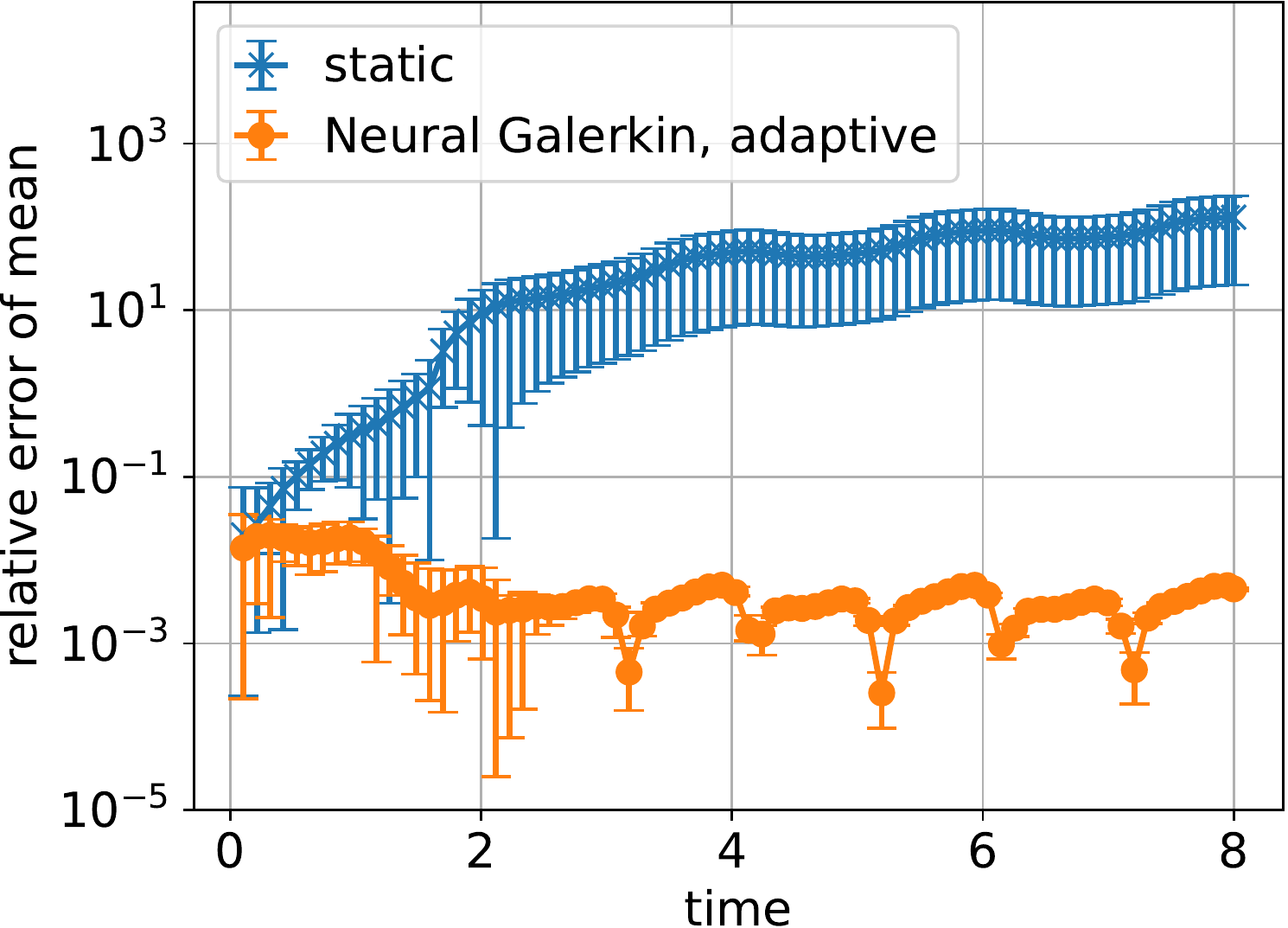} & \includegraphics[width=0.45\columnwidth]{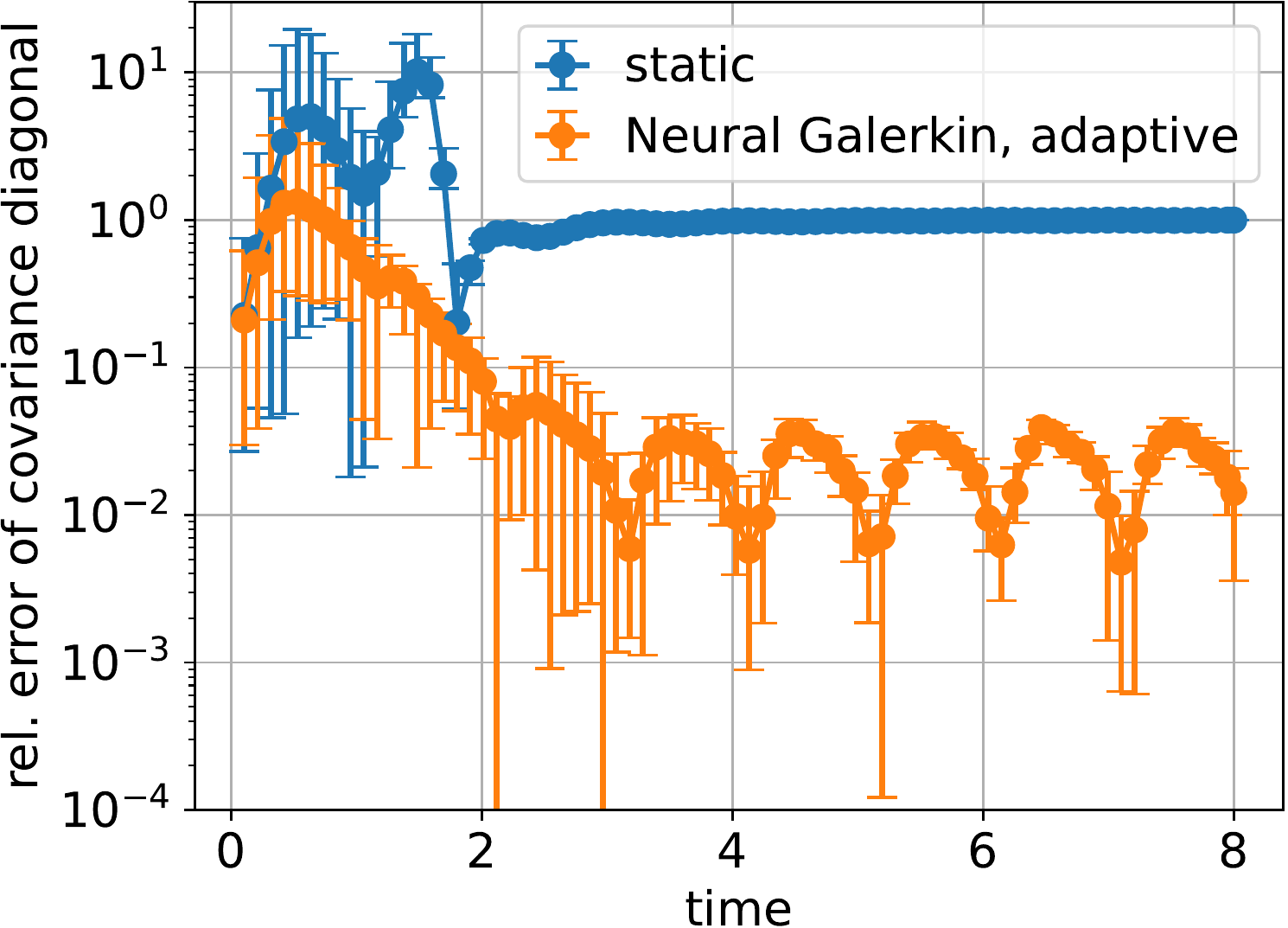}\\
\footnotesize (a) mean positions & \footnotesize (b) covariance of particles\\
\end{tabular}
\caption{Particles in aharmonic trap: We use the same plotting convention as in Figure~\ref{fig:ParticleDim8}}
\label{fig:ParticleCubicDim8}
\end{figure}

To apply a Neural Galerkin scheme in this experiment, we take the neural network defined in Eq.~\eqref{eq:NumExp:Exp:ShallowNetwork} with $m = 40$ nodes and non-negative weights. Time is discretized with implicit Euler and $\delta t = 10^{-3}$. We take $n = 5000$ sample and $1000$ SGD iterations in each time step. The rest of the setup is the same as in the example with the harmonic trap. Figure~\ref{fig:ParticleCubicDim8}(a) and (b) show the relative error of the mean and the relative error of the covariance diagonal, respectively. With adaptive sampling, the used Neural Galerkin scheme achieves an error between $10^{-3}$ and $10^{-2}$. In contrast, static sampling is insufficient and leads to large relative errors. Note that the off-diagonal elements of the covariance of the particle density converge to zero over time in this experiment and thus the relative error is not informative and not plotted here. Because Neural Galerkin schemes provide an approximation of the density of the particles, we can efficiently compute the entropy: As  shown in Figure~\ref{fig:ParticleCubicDim8Entropy}(b), it is in agreement with the prediction obtained with Monte Carlo and density estimation from $100,000$ samples.

\section{Concluding remarks}

As machine-learning techniques via deep learning play an ever more important role in enabling realistic simulations of phenomena and processes in science and engineering, there is an increasing need for training on data that are informative about the underlying physics. This problem of data generation becomes especially important in dynamical systems that evolve over time where \emph{a~priori} data is sparse or non-existent, and static data collection in a pre-processing step is insufficient. The proposed Neural Galerkin schemes intertwine data acquisition via sampling and numerical simulations in an adaptive way. The schemes learn PDE solutions by evolving a non-linear function representation, which is adjusted from data points obtained via importance sampling that is informed by the solution itself.
Our numerical experiments demonstrate that the adaptivity in function approximation as well as data collection is key to enable accurate prediction of local features in high-dimensional domains.

These results suggest that the proposed Neural Galerkin schemes can be applied to simulate many other high-dimensional equations relevant in science and engineering, such as kinetic equations, non-linear Fokker-Planck equations, Boltzmann equations, etc.
This would be especially useful in situations where Monte-Carlo methods are not directly applicable because the solution does not admit a representation via the Feynman-Kac formula. The usage of Neural Galerkin schemes in such applications also opens the door to several interesting questions that we leave as future research avenues: First, we should understand better which structural properties of the neural architectures will guarantee that the function approximation from \eqref{eq:rep:sol} $u(t, \bfx) \approx U(\theta(t), \bfx)$ will hold uniformly in time. For some initial efforts in the context of shallow ReLU architectures we refer the reader to~\cite{wojtowytsch2020some}. Second, we should quantify when the integral defining the instantaneous loss \eqref{eq:NG:ContTime:MinResObj} admits an efficient estimation, beyond what non-adaptive Monte-Carlo schemes can provide, by exploiting sparsity of the residual in a certain domain, not necessarily the spatial domain. For example, if the solution is represented in the frequency domain, then locality (often referred to as sparsity) means that only few frequencies are relevant, which in principle can be exploited numerically with Neural Galerkin schemes by using suitable nonlinear parametrizations. Third, it would be interesting to explore and systematize importance sampling strategies that adapt to the transient regularity structures appearing in the solution and thereby will allow us to maintain an adaptive neuron budget. Such strategies could for example use birth-death processes \cite{rotskoff2019global}, which would allow us to vary the width of the neural network in time, or methods that couple the PDE to some evolution equation for the input data points used to sample the residual.

\bibliography{ng}
\bibliographystyle{abbrv}

\appendix
\section{Learning in the Banach space $\mathcal{F}_1$}
\label{sec:F1}
Let us consider the Neural Galerkin representation $u(t, \bfx) = U(\theta(t), \bfx)$ 
in a situation where $\theta$ is infinite dimensional, namely when
\begin{equation}
    \label{eq:uF1}
    u(t,\bfx) = \int_\Zcal \phi(\bfx,z) d\gamma_t(z)
\end{equation}
where $\phi: \Xcal\times\Zcal \to \R$ is some nonlinear unit (like e.g. the ReLU), and $\gamma_t$ is a Radon measure with finite total variation defined on $\Zcal$: we will denote the space of such Radon measures by $\Mcal(\Zcal)$. Functions of the type defined in Eq.~\eqref{eq:uF1} form a Banach space in which the norm of $u(t)$ is the minimum total variation of $\gamma_t$ among all Radon measures consistent with Eq.~\eqref{eq:uF1}; these function can also be approximated by shallow (two-layer) neural networks, which converge to functions as defined in Eq.~\eqref{eq:uF1} in the limit of infinite width. In this context, one can specify the evolution of $\gamma_t$ by requiring that
\begin{equation}
    \label{eq:gammadot:var}
    \dot \gamma_t \in \argmin_{\dot\gamma\in\Mcal (\Zcal)} \int_\Xcal \left| \int_\Zcal \phi(\bfx,z) d\dot \gamma(z) - f\left(t,\bfx, {\textstyle\int_\Zcal} \phi(\cdot,z) d\gamma_t(z)\right) \right|^2 d\nu(\bfx)
\end{equation}
where $\nu$ is some positive measure on $\Xcal$. The Euler-Lagrange equation associated with this minimization problem gives the evolution equation for $\gamma_t$:
\begin{equation}
    \label{eq:gamma:evol}
    \int_{\Zcal} M(z,z') d\dot \gamma_t(z') = F(t,z,\gamma_t) 
\end{equation}
where
\begin{equation}
    \label{eq:MF:gam:def}
    \left\{
    \begin{aligned}
    M(z,z') & = \int_\Xcal \phi(\bfx,z)\phi(\bfx,z') d\nu(\bfx),\\
    F(t,z,\gamma) &= \int_\Xcal \phi(\bfx,z) f\left(t,\bfx, {\textstyle\int_\Zcal} \phi(\cdot,z') d\gamma(z')\right)d\nu(\bfx)\,.
    \end{aligned}
    \right.
\end{equation}

\section{Another derivation of the Neural Galerkin equation $M(\theta)\dot{\theta} = F(t, \theta)$}
\label{sec:deriv:NGODE}

Let $u(t)$ be the solution of the PDE at time $t$ and $U(\theta(t))$ be its parametric representation. In general we do not have access to~$u(t)$ (except at initial time through the prescribed initial condition), but if we had an oracle giving us this solution one option to determine the parameters $\theta(t)$ would be to minimize an objective like
\begin{equation}
    \label{eq:loss:instant}
    E(\theta(t)) =  R\left(u(t),U(\theta(t))\right) 
\end{equation}
where $R:\Ucal\times\Ucal \to \R$ is such that (i) $R(u,v)\ge 0$ for all $u,v\in\Ucal$; (ii) $R(u,u) = 0$ for all $u\in \Ucal$ and (iii) $R(u,v)$ is strictly convex in $v$ for all $u\in \Ucal$. 
Assuming that $R(u,v)$ is twice differentiable in both its arguments, this requires that
\begin{equation}
    \label{eq:convex}
    \forall u \in \Ucal \quad :\quad  R(u,u) =0, \quad D_v R(u,v)|_{v=u} =0 \quad \text{and}\quad D^2_{v,v} R(u,v)|_{v=u} \quad \text{is positive-definite}
\end{equation}
where $D$ denotes derivative in $\Ucal$. Together with $R(u,u) = 0$, note that this also implies that $D_u R(u,v)|_{v=u} =0$ for all $u\in\Ucal$. 

Suppose that the current solution $u(t)$ is representable exactly, i.e. there exists $\theta(t)\in\Theta$ such that $U(\theta(t))=u(t)$. This implies that $E(\theta(t)) =0$, and it is natural to require that this objective remains zero, or as close to zero as possible, as time evolves. With this in mind, we can use the following result:

\begin{proposition}
\label{th:loss}
Assume that (i) $U(\theta(t))=u(t)$  and (ii) the objective $R(u,v)$ satisfies the conditions~\eqref{eq:convex}, so that $E(\theta(t)) =0$. Then
\begin{equation}
    \label{eq:limit:loss}
     \lim_{h\to0} h^{-2} E(\theta(t+h)) = \tfrac12 \dot \theta^T M(\theta(t) )\dot \theta + \dot \theta \cdot F(t,\theta(t)) + \tfrac12 b(t,\theta(t)) =: L(\dot \theta),
\end{equation}
where
\begin{equation}
    \label{eq:MF2}
    \begin{aligned}
    M(\theta) &= \<\nabla_\theta U(\theta), D^2_{v,v} R(U(\theta),U(\theta))\nabla_\theta U(\theta)\>\\
    F(t,\theta) & = \<f(t,U(\theta)), D^2_{u,v} R(U(\theta),U(\theta))\nabla_\theta U(\theta)\>\\
    b(t,\theta) & = \<f(t,U(\theta)), D^2_{u,u} R(U(\theta),U(\theta))f(t,U(\theta))\>
    \end{aligned}
\end{equation}
\end{proposition}

\noindent\textit{Proof:} A straightforward calculation using the conditions in~\eqref{eq:convex} as well as $\partial_t u(t)=f(t,u(t)) = f(t,U(\theta(t))$ indicates that
\begin{equation}
    \label{eq:loss:taylor}
    \begin{aligned}
    E(\theta(t+h)) &=  R\left(u(t+h),U(\theta(t+h))\right) \\
    &= \tfrac12 h^2 \< f(t,U(\theta(t)), D^2_{u,u} R(U(\theta(t)),U(\theta(t)))f(t,U(\theta(t))\>\\
    &+\tfrac12 h^2 \< \nabla_\theta U(\theta(t))\cdot \dot \theta, D^2_{v,v} R(U(\theta(t)),U(\theta(t)))\nabla_\theta U(\theta(t))\cdot \dot \theta\>\\
    &\quad+ h^2 \< f(t,U(\theta(t)), D^2_{u,v} R(U(\theta(t)),U(\theta(t)))\nabla_\theta U(\theta(t))\cdot \dot \theta\> +o(h^2)
    \end{aligned}
\end{equation}
Dividing by $h^2$ and taking the limit as $h\to0$ establishes~\eqref{eq:limit:loss}.\hfill$\square$

Minimizing the objective defined in Eq.~\eqref{eq:limit:loss} over $\dot \theta$ given $\theta(t)$ leads us back to an evolution equation for $\theta(t)$ that is structurally similar the Neural Galerkin equation derived in the main text, 
 that is
\begin{equation}
    \label{eq:ODE2}
    M(\theta(t)) \dot \theta(t) = F(t,\theta(t)),
\end{equation}
and is to be solved with the initial condition
\begin{equation}
    \label{eq:init:ODE}
    \theta(0) \in \argmin_\theta R(u_0,U(\theta)).
\end{equation} 
In particular it is easy to see that if we take
\begin{equation}
    \label{eq:loss:1}
    R(u,v) = \int_\Xcal |u(t,\bfx)-v(t,\bfx)|^2 d\nu(x),
\end{equation}
then the system of ODEs defined in Eq.~\eqref{eq:ODE2} is exactly the same system as the system of ODEs given in the main text with $M$ and $F$ given in Eq.~\eqref{eq:M:F:def}  with $\nu_{\theta}=\nu$. 
The formulation above is however more general, and offers a principled way to specify the measure $\nu_{\theta}$.  
For example, for equations like the Fokker-Planck equation that evolve a probability density~$u(t)>0$ such that $\int_\Xcal u(t,\bfx) d\bfx =1$, it common objective is the Kullback-Leibler divergence
\begin{equation}
    \label{eq:DKL}
    R(u,v) = \int_{\Xcal}\log \left(\frac{u(\bfx)}{v(\bfx)}\right) u(\bfx) d\bfx\,.
\end{equation}
With this choice  we arrive at Eq.~\eqref{eq:ODE2} with
\begin{equation}
    \label{eq:MF3}
    \begin{aligned}
    M(\theta) &= \int_\Xcal \frac{\nabla_\theta U(\theta,\bfx)\otimes \nabla_\theta U(\theta,\bfx)}{U(\theta,\bfx)} d\bfx\\
    F(t,\theta) & = \int_\Xcal \frac{\nabla_\theta U(\theta,\bfx)f(t,\bfx, U(\theta))}{U(\theta,\bfx)} d\bfx
    \end{aligned}
\end{equation}
These correspond to using $d\nu_{\theta}(\bfx) = |U(\theta,\bfx)|^{-1} d\bfx $.

\section{Comparison with global methods}
\label{sec:comp:coll}

Let us compare the Neural Galerkin equation we use with the equations one obtains by treating the problem globally in time, i.e. by putting the spatial and the tempral variables on the same footing as in the Deep Ritz method~\cite{doi:10.1002/cnm.1640100303,SIRIGNANO20181339,RAISSI2019686,BERG201828}. 
An optimal trajectory (in the $L_2([0,T] \times \mathcal{X}, \nu_0)$ sense) in parameter space satisfies 
\begin{equation}
    \label{eq:loss:global}
    \min_{\{\theta(t) \in \Theta; t\in [0, T]\}} \int_0^T J_t(\theta, \dot \theta(t)) dt,
\end{equation}
subject to $\theta(0)=\theta_0$.
The Euler-Lagrange equations associated with this minimization problem are
\begin{equation}
    \label{eq:EL}
    \frac{d}{dt} \partial_\eta J_t(\theta,\dot \theta) = \partial_\theta J_t(\theta,\dot \theta), \qquad \theta(0) = \theta_0, \qquad \partial_\eta J_T(\theta(T),\dot \theta(T)) =0.
\end{equation}
Explicitly, these read
\begin{equation}
\label{eq:EL:ex}
    \frac{d}{dt} \left( M(t,\theta) \dot \theta - F(t,\theta)\right) = H(t,\theta,\dot \theta), \qquad \theta(0) = \theta_0, \qquad M(T,\theta(T)) \dot \theta(T) = F(T,\theta(T)),
\end{equation}
where we defined
\begin{equation}
    \label{eq:Hdef}
    H(t,\theta,\eta) = \int_\Xcal \left(\nabla_\theta\nabla_\theta U(\theta,\bfx) \eta - \nabla_\theta U(\theta,\bfx) D_u f(t,\bfx,U(\theta))\right) \left(\nabla_\theta U(\theta,\bfx)\cdot \eta -  f(t,\bfx,U(\theta)) \right) d\nu_{t,\theta}(\bfx)
\end{equation}
Eqs.~\eqref{eq:EL:ex} are a boundary value problem which is harder to solve than the Neural Galerkin Eq.~\eqref{eq:NG:ODE}.

\section{Advection in unbounded, high-dimensional domains}
\label{sec:Appx:AdvEq}
\begin{figure}
\includegraphics[width=1.0\columnwidth]{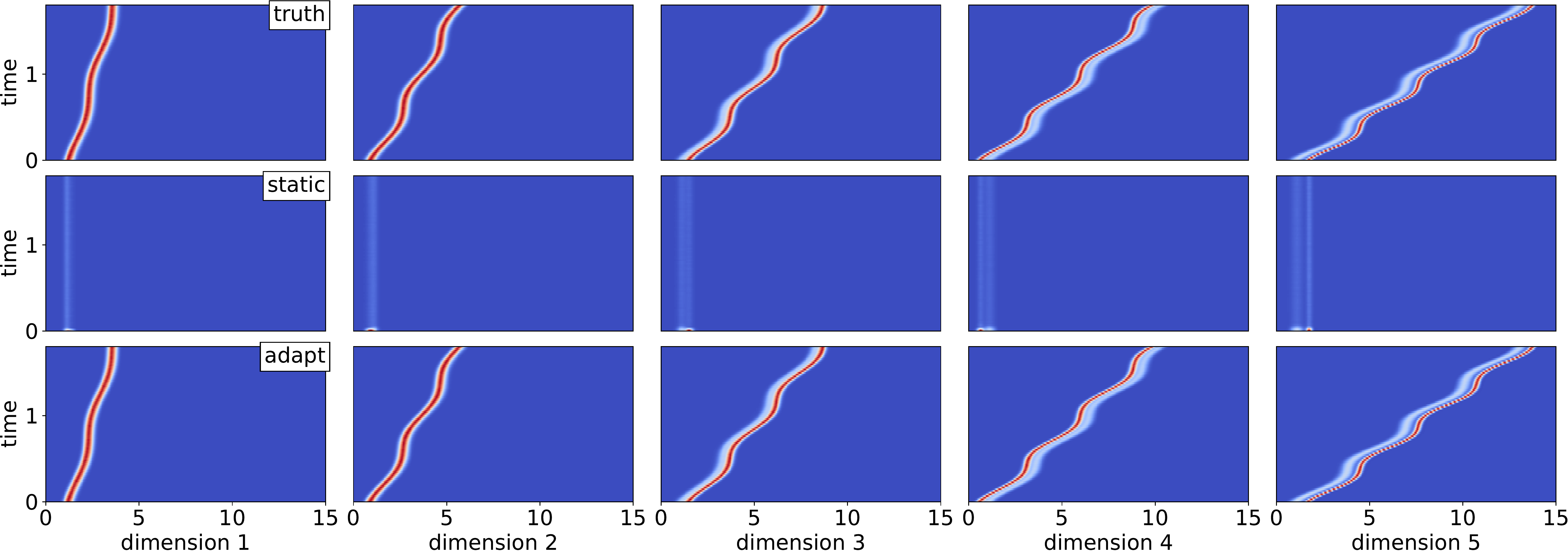}
\caption{High-dimensional transport with time-varying coefficient: The plots labeled with ``truth'' show the marginals of the analytic solution. The plots labeled with ``uniform'' are the Neural Galerkin approximations with a static uniform sampling in $[0, 15]^d$. The plots labeled with ``adapt'' show the Neural Galerkin solution with adaptive sampling, which is in close agreement with the analytic solution.}
\label{supp:fig:NumExp:TransportTime}
\end{figure}

\begin{figure}
\centering
\includegraphics[width=1.0\columnwidth]{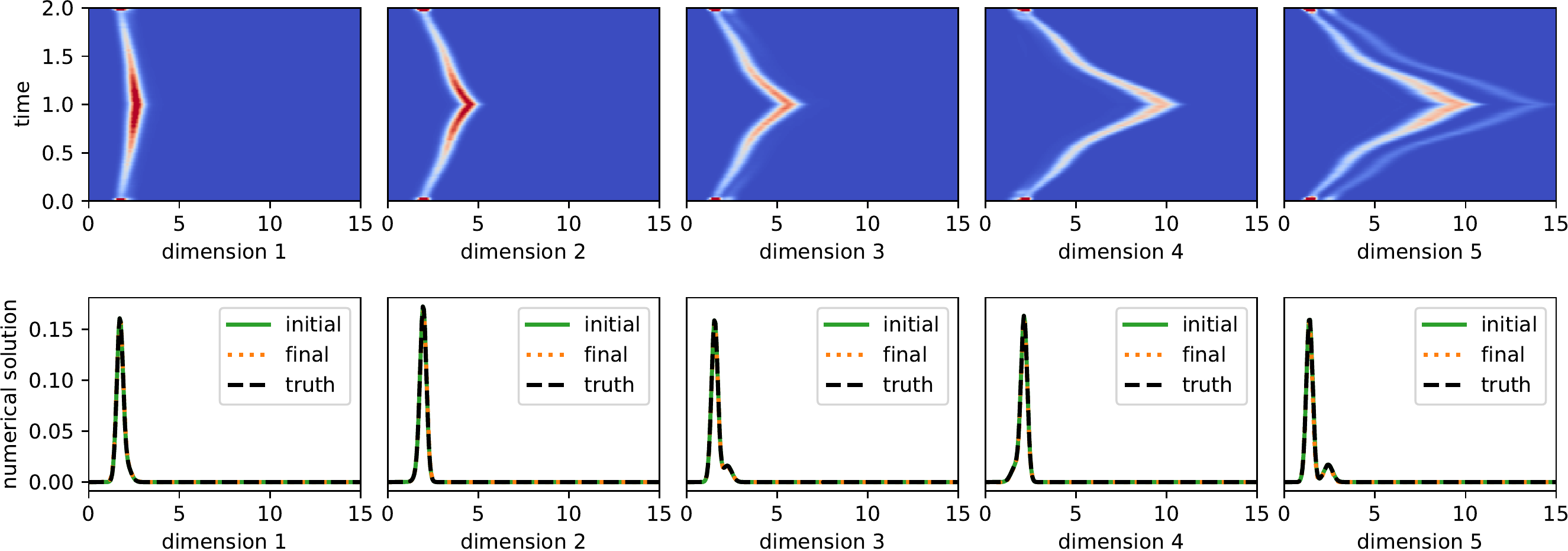}
\caption{High-dimensional transport with time and spatially varying coefficients: The first row shows the marginals of the approximation obtained with Neural Galerkin with adaptive sampling and adaptive time integration. Because of the spatially varying coefficient, the coherent structure that is propagated over time changes shape, which is different from the time-only varying case shown in Figure~\ref{supp:fig:NumExp:TransportTime}. The transport equation is integrated forward and backward in time to compare the approximation at initial and end time, which is shown in the second row. The approximations at initial and end time match well, which indicates that the proposed adaptive approach approximates well the local, high-dimensional dynamics in this problem.}
\label{supp:fig:LinAdvSpatiallyVarying}
\end{figure}

All five marginals for the time-varying and the time-and-space varying case are shown in Figure~\ref{supp:fig:NumExp:TransportTime} and Figure~\ref{supp:fig:LinAdvSpatiallyVarying}, respectively.

\section{Fokker-Planck equation in high-dimensions} \label{sec:Appx:FokkerPlanck}

\begin{figure}
\begin{tabular}{ccc}
\includegraphics[width=0.30\columnwidth]{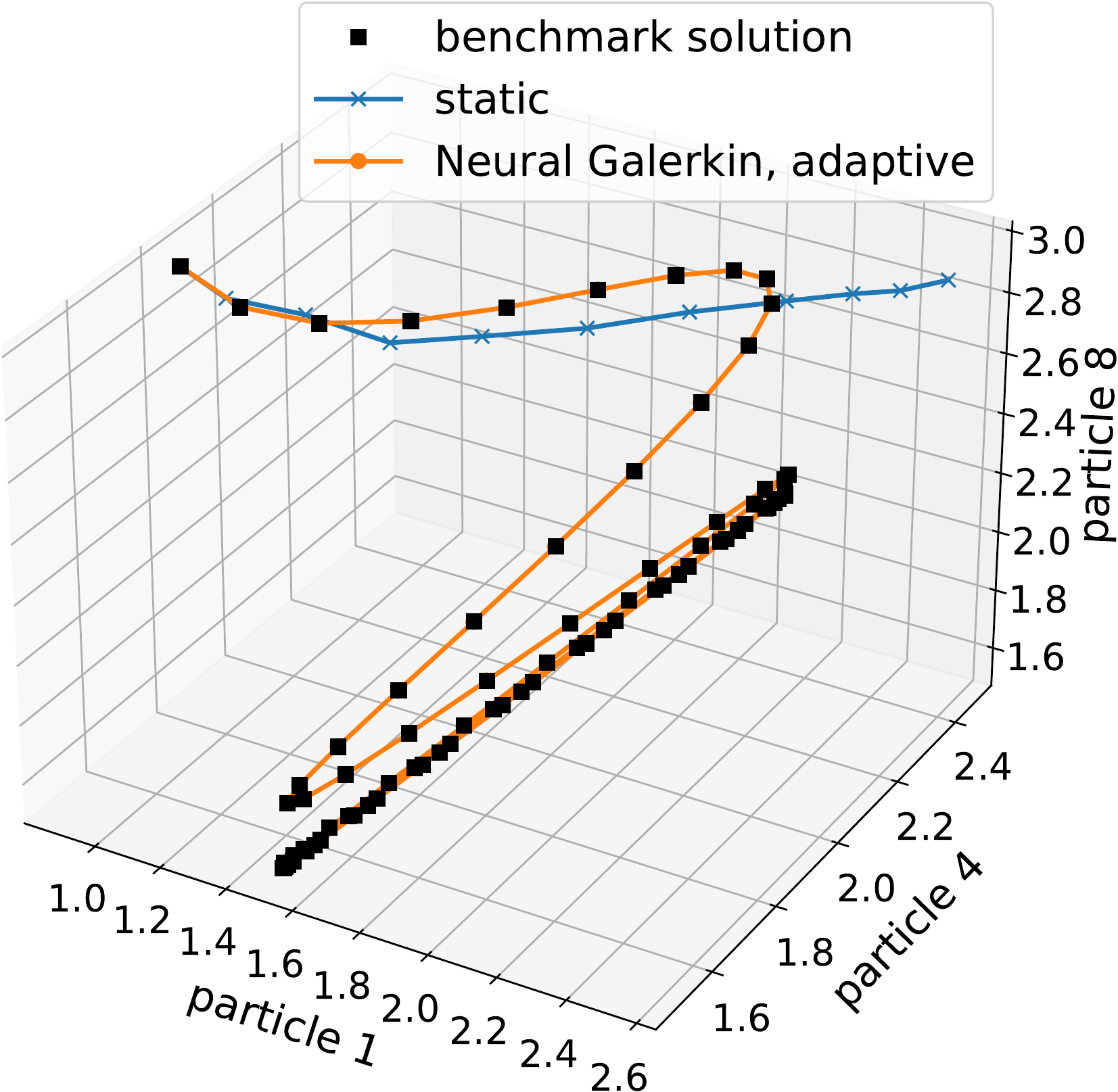} & \includegraphics[width=0.30\columnwidth]{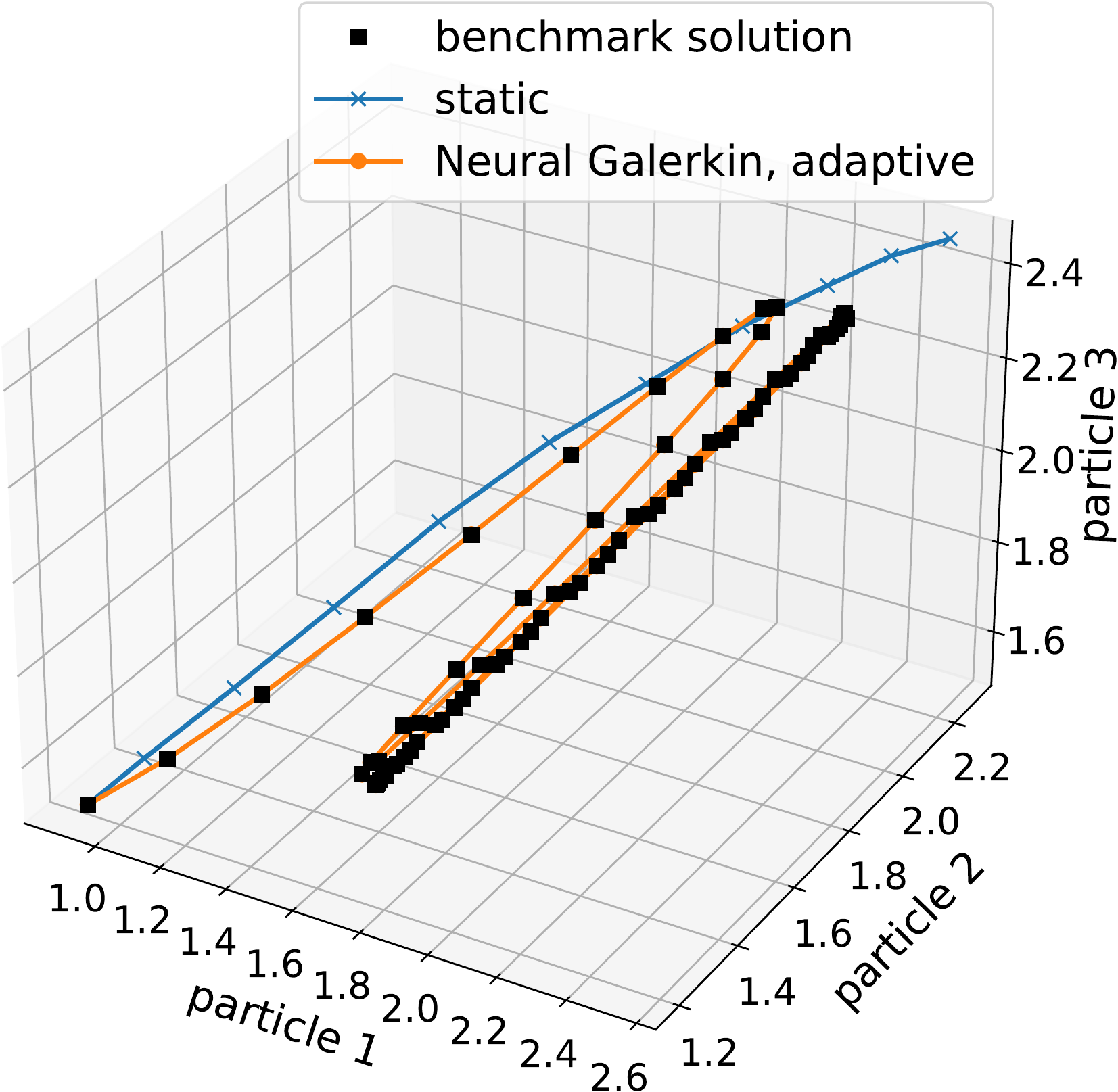} & \includegraphics[width=0.30\columnwidth]{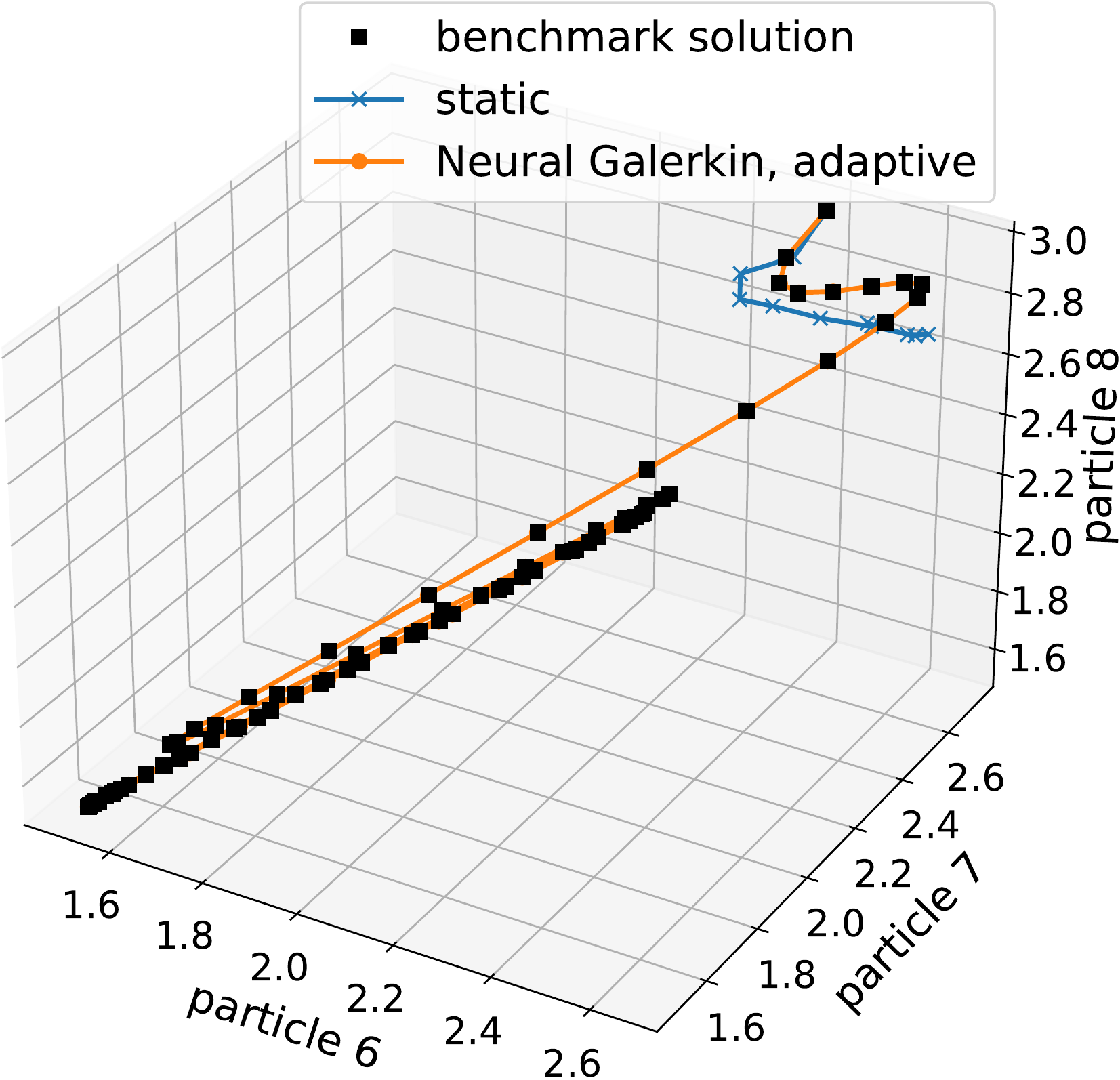}
\end{tabular}
\caption{Particles in harmonic trap: Visualization of particle positions over time.}
\label{supp:fig:Particle3DPlots}
\end{figure}

In case of the harmonic trap, the following system of ODEs describe the mean $\bar{X}_i(t)$ and covariance $C_{ij}(t)$ of particles $X_i$ and $X_j$ at time $t$ for $i , j = 1, \dots, d$ are 
\begin{align}
\dot{\bar{X}}_i(t) = & -(1 + \alpha)\bar{X}_i(t) + a(t) + \frac{\alpha}{d}\sum_{j = 1}^d \bar{X}_j(t)\,,\qquad i = 1, \dots, d\,,\\
\dot{C}_{ij} = & -2(1 + \alpha)C_{ij} + \frac{\alpha}{d}\sum_{k = 1}^d (C_{kj} + C_{ki}) + 2\beta^{-1}\delta_{ij}\,,\qquad i,j = 1, \dots, d\,,\label{eq:SI:FPEODECov}
\end{align}
where $\delta_{ij}$ is the Dirac delta. The benchmark solution is obtained by integrating these ODEs with Runge-Kutta 45. Figure~\ref{supp:fig:Particle3DPlots} compares the particle positions predicted by Neural Galerkin to the benchmark solution for various dimensions.

\begin{figure}
\begin{tabular}{ccc}
\includegraphics[width=0.30\columnwidth]{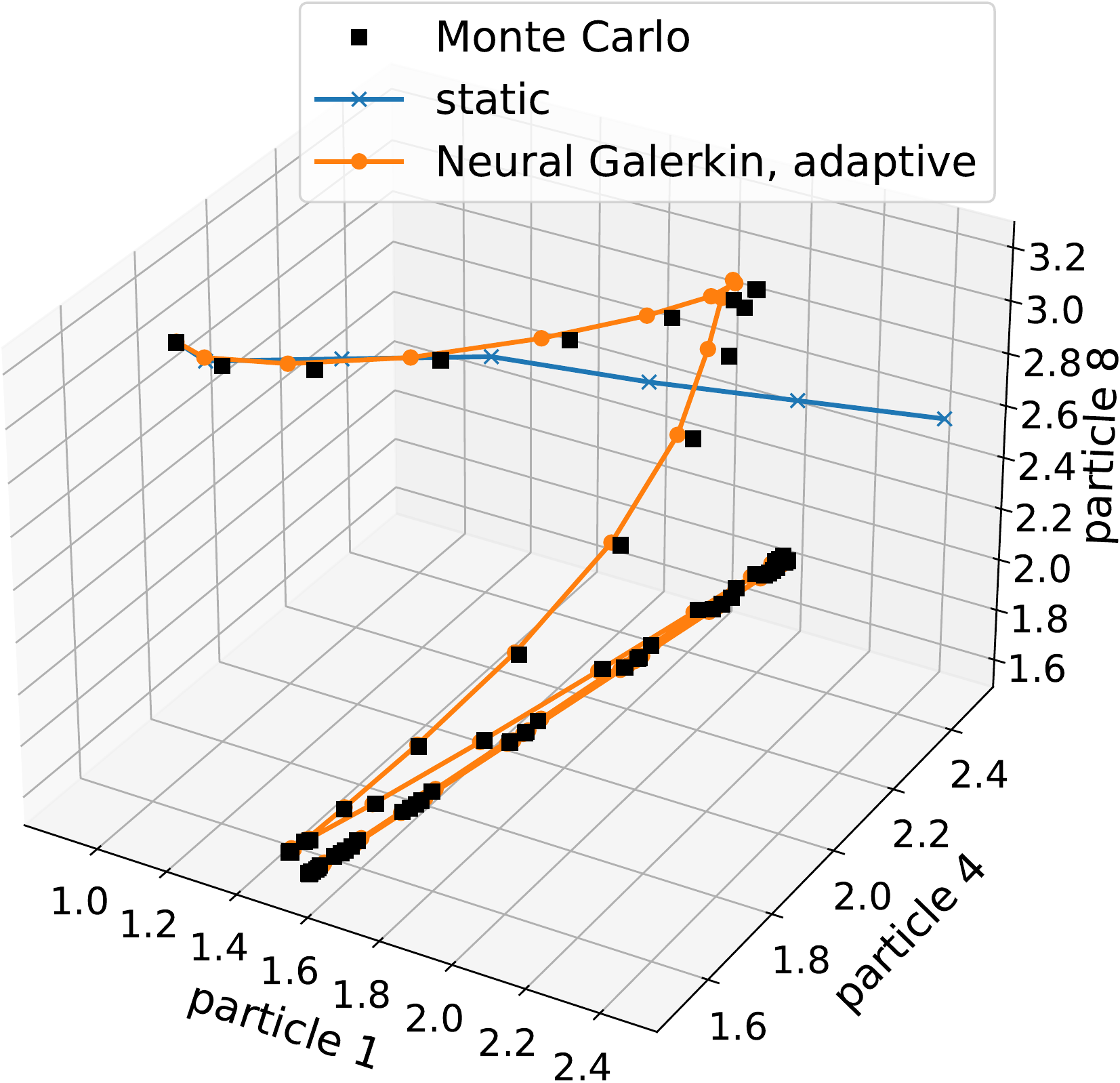} & \includegraphics[width=0.30\columnwidth]{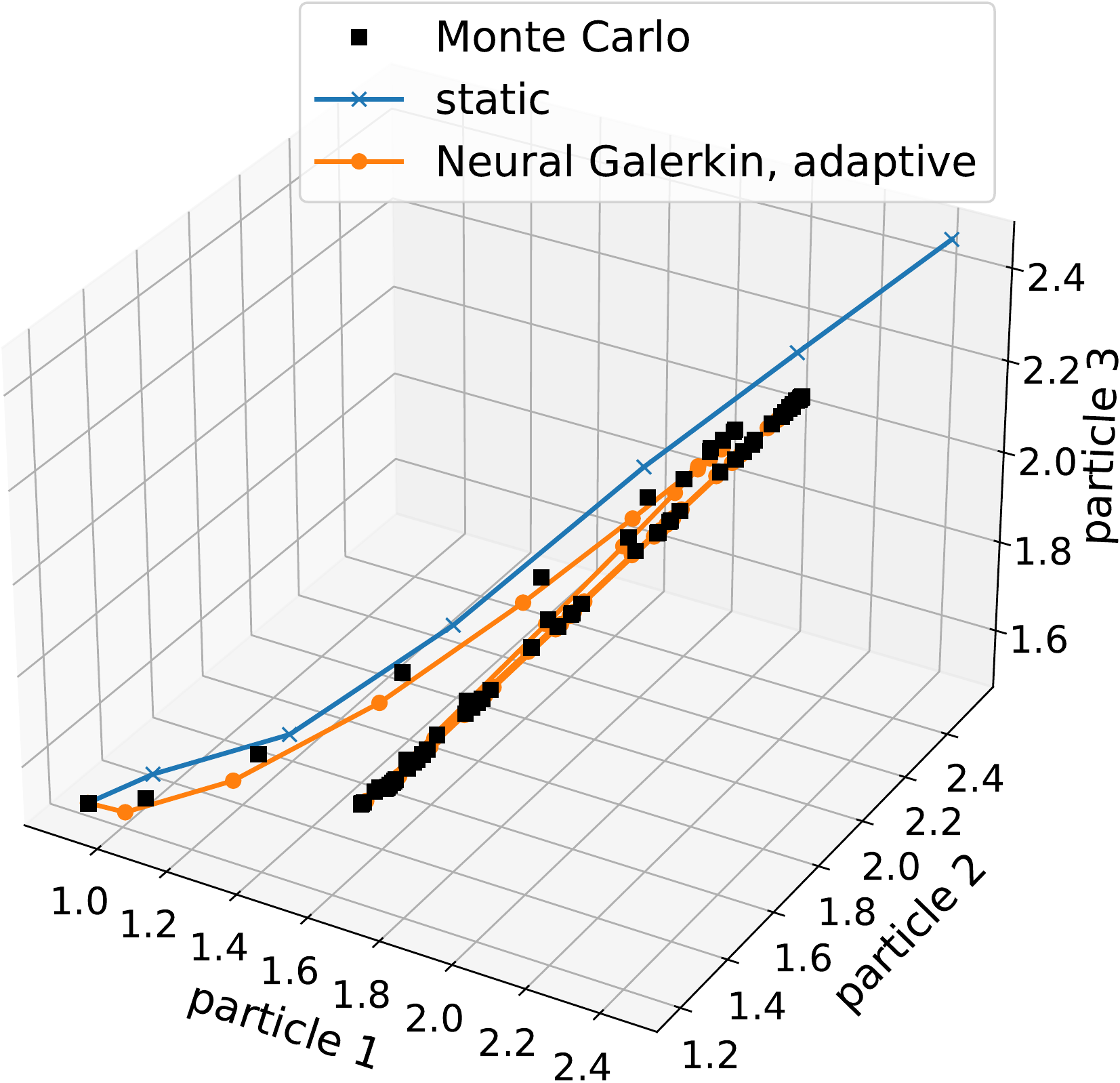} & \includegraphics[width=0.30\columnwidth]{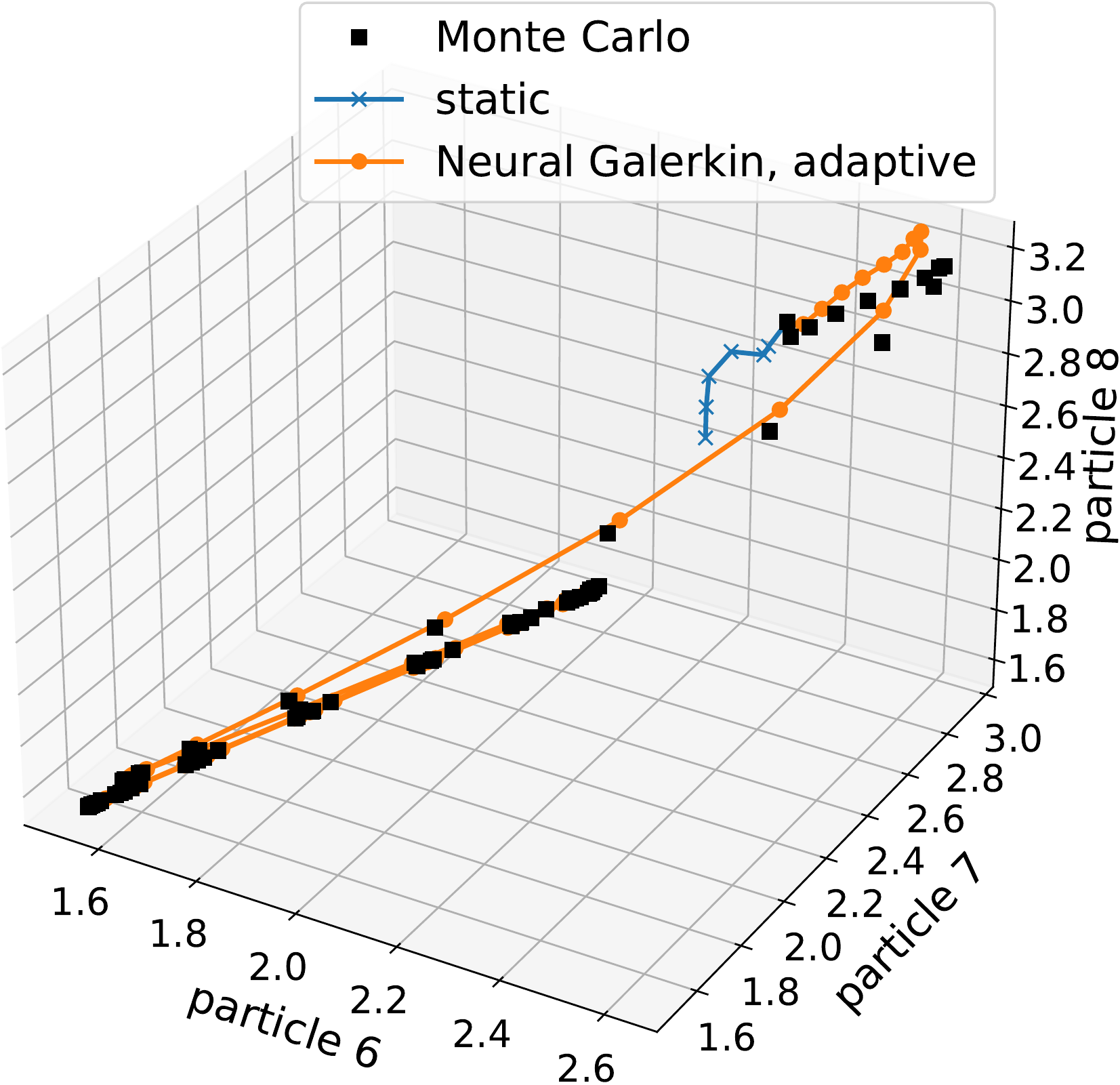}
\end{tabular}
\caption{Particles in aharmonic trap: Visualization of particle positions over time. Note that a Monte Carlo estimate of mean serves as benchmark here.}
\label{supp:fig:ParticleAHarmonic3DPlots}
\end{figure}

In case of the aharmonic trap, a Monte Carlo estimate serves as benchmark. To that end, we discretize the stochastic differential equation of the particle positions 
via the Euler–Maruyama method with time-step size $\delta t = 10^{-4}$ and draw $100,000$ paths. From these $100,000$ paths we estimate the mean and covariance, which serve as the benchmark for the Neural Galerkin approximation. The particle positions predicted by Neural Galerkin versus the Monte Carlo benchmark are shown in Figure~\ref{supp:fig:ParticleAHarmonic3DPlots}.

\end{document}